\documentclass[11pt]{article}
 \usepackage{benstyle}
 
\usepackage{cite}

\newcommand {\lb} {\label}
\newcommand {\rf[1]} {(\ref{#1})}
\usepackage{url}

\title{{Fast and stable approximation of analytic functions from equispaced samples via polynomial frames}}
\author{Ben Adcock\thanks{Department of Mathematics, Simon Fraser University, Burnaby BC, V5A 1S6, Canada, \texttt{ben\_adcock@sfu.ca}, +1-778-782-4819}
\and
Alexei Shadrin\thanks{DAMTP, University of Cambridge, Cambridge CB3 0WA, UK, \texttt{a.shadrin@damtp.cam.ac.uk}, +44-1223-766-887}
}

\begin{document}

\maketitle

\begin{abstract}
We consider approximating analytic functions on the interval $[-1,1]$ from their values at a set of $m+1$ equispaced nodes. A result of Platte, Trefethen \& Kuijlaars states that fast and stable approximation from equispaced samples is generally impossible. In particular, any method that converges exponentially fast must also be exponentially ill-conditioned. We prove a positive counterpart to this `impossibility' theorem. Our `possibility' theorem shows that there is a well-conditioned method that provides exponential decay of the error down to a finite, but user-controlled tolerance $\epsilon > 0$, which in practice can be chosen close to machine epsilon. The method is known as \textit{polynomial frame} approximation or \textit{polynomial extensions}. It uses algebraic polynomials of degree $n$ on an extended interval $[-\gamma,\gamma]$, $\gamma > 1$, to construct an approximation on $[-1,1]$ via a SVD-regularized least-squares fit. A key step in the proof of our main theorem is a new result on the maximal behaviour of a polynomial of degree $n$ on $[-1,1]$ that is simultaneously bounded by one at a set of $m+1$ equispaced nodes in $[-1,1]$ and $1/\epsilon$ on the extended interval $[-\gamma,\gamma]$. We show that linear oversampling, i.e., $m = c n \log(1/\epsilon) / \sqrt{\gamma^2-1}$, is sufficient for uniform boundedness of any such polynomial on $[-1,1]$. This result aside, we also prove an extended impossibility theorem, which shows that such a possibility theorem (and consequently the method of polynomial frame approximation) is essentially optimal.
\end{abstract}


\pbk

\section{Introduction}

In this paper, we consider the problem of approximating an analytic function $f : [-1,1] \rightarrow \bbC$ from its values at $m+1$ equispaced points in $[-1,1]$. Several years ago, Platte, Trefethen \& Kuijlaars \cite{platte2011impossibility} demonstrated that this problem is intrinsically difficult. They proved an `impossibility' theorem which states that any method that offers exponential rates of convergence in $m$ for all functions analytic in a fixed, but arbitrary region of the complex plane must necessarily be exponentially ill-conditioned. Furthermore, the best rate of convergence achievable by a stable method is necessarily subexponential -- specifically, root-exponential in $m$.

This result generalizes what has long been known for polynomial interpolation at equispaced nodes: namely, Runge's phenomenon. Polynomial interpolation is divergent for functions that not analytic in a sufficiently large complex region (the Runge region). And while it converges exponentially fast for functions that are analytic in this region, its condition number is also exponentially large. Such ill-conditioning means that, when computed in floating point arithmetic, the error of polynomial interpolation eventually increases, even for entire functions, due to the accumulation of round-off error.

Many methods have been proposed to overcome Runge's phenomenon by stably and accurately approximating analytic functions from equispaced nodes (see, e.g., \cite{adcock2016mapped,boyd2009exponentially,platte2011impossibility} and references therein). Several such methods appear to offer fast convergence in practice{, which, for some functions at least, appears faster than root-exponential.} Yet full mathematical explanations for this phenomenon have hereto been lacking.

The purpose of this paper is to fully analyze one such method in view of the impossibility theorem. This method is termed \textit{polynomial frame approximation} or \textit{polynomial extensions} \cite{adcock2020approximating}, and is closely related to so-called \textit{Fourier extensions} \cite{matthysen2017function,lyon2012fast,huybrechs2010fourier,bruno2007accurate,boyd2002fourier,pasquetti1996spectral}. It approximates a function $f$ on $[-1,1]$ using orthogonal polynomials on an \textit{extended} interval $[-\gamma,\gamma]$ for some fixed $\gamma > 1$. The approximation is then computed by solving a regularized least-squares problem, with user-controlled regularization parameter $\epsilon > 0$.

Our main contribution is to show that this method offers a positive counterpart to the impossibility theorem of \cite{platte2011impossibility}. We prove a `possibility' theorem, which asserts that the polynomial frame approximation method is well-conditioned and, for all functions that are analytic in a sufficiently large region (related to the parameter $\gamma$), its error decreases exponentially fast down to roughly $\ordu{m^{3/2} \epsilon}$, where $\epsilon$ is the user-determined tolerance. This tolerance may be taken to be on the order of machine epsilon without impacting the conditioning of the method, thus rendering the approach suitable for practical purposes. But it may also be taken larger (with benefits in terms of stability) if fewer digits of accuracy are required. 
While the impossibility theorem dictates that exponential convergence to zero cannot be achieved by a well-conditioned method, we show that exponential decrease of the error down to an arbitrary tolerance, multiplied by the slowly growing factor proportional to $m^{3/2}$, is indeed possible. Additionally, we establish an `extended' impossibility theorem, which relates conditioning and error decay for approximation methods that achieve only a finite accuracy $\epsilon$. This theorem explains how the method we consider is essentially optimal, in a suitable sense.

Our main result hinges on a new bound for the maximal behaviour of polynomials that are bounded simultaneously at a set of $m+1$ equispaced nodes on $[-1,1]$ and on the extended interval $[-\gamma,\gamma]$. The impossibility theorem of \cite{platte2011impossibility} uses a classical result of Coppersmith and Rivlin \cite{coppersmith1992growth}, which states that a polynomial of degree $n$ that is bounded by one at $m \geq n$ equispaced nodes can grow as large as $\alpha^{n^2/m}$ on $[-1,1]$ outside of these nodes, where $\alpha > 1$ is a constant. In particular, $m = c n^2$ equispaced points are both sufficient and necessary to ensure boundedness of such a polynomial on $[-1,1]$. We consider a nontrivial variation of this setting, where the polynomial is also assumed to be no larger than $1/\epsilon$ on $[-\gamma,\gamma]$. Our key result shows that $m = c n \log(1/\epsilon) / \sqrt{\gamma^2-1}$ equispaced points suffice for ensuring boundedness of any such polynomial on $[-1,1]$. While we use this bound to analyze polynomial frames, we expect it to be of independent interest from a pure approximation-theoretic perspective.

{
More broadly, our work has some connections with optimal recovery.  We note that optimal recovery of functions in general, and from equispaced samples 
in particular has been actively studied in approximation theory, with emphasis on
optimal approximation methods without taking stability issues into considerations. See, for instance, \cite{temlyakov1993approximate,michelli1985lecture,temlyakov2018multivariate}.
}

The outline of the remainder of this paper is as follows. In \S \ref{s:overview} we present an overview of the main components of the paper. We introduce notation, describe the impossibility theorem of \cite{platte2011impossibility} and then present polynomial frame approximation. We next state our main results, then conclude with a discussion of related work. In \S \ref{s:acc-stab-poly-frame} we analyze the accuracy and conditioning of polynomial frame approximation. Then in \S \ref{s:maximal-behaviour-proof} we give the proof of the aforementioned result on the maximal behaviour of polynomials bounded at equispaced nodes and on the extended interval $[-\gamma,\gamma]$. With this in hand, in \S \ref{s:main-thm-proof} we give the proof of the main results: the possibility theorem for polynomial frame approximation and the extended impossibility theorem. We then present several numerical examples in \S \ref{s:numerical}, before offering some concluding remarks in \S \ref{s:conclusions}.

\section{Overview of the paper}\label{s:overview}

\subsection{Notation}

Throughout, $\bbP_n$ denotes the space of polynomials of degree at most $n$.
For $m \geq 1$, we let $\{ x_i \}^{m}_{i=0}$ be a set of $m+1$ equispaced points in $[-1,1]$ including endpoints, i.e.\ $x_i = -1 + 2 i /m$. Given an interval $I$, we let $C(I)$ be the space of continuous functions on $I$ and
\bes{
\nm{g}_{I,\infty} = \sup_{x \in I} |g(x)|,\quad g \in C(I),
}
be the uniform norm over $I$.
We also let
\bes{
\ip{f}{g}_{I,2} = \int_{I} f(x) \overline{g(x)} \D x,\quad f , g \in C(I),
}
be the usual $L^2$-inner product over $I$ and $\nm{\cdot}_{I,2} = \sqrt{\ip{\cdot}{\cdot}_{I,2}}$ be the corresponding $L^2$-norm. 

Next, we define several discrete semi-norms and semi-inner products. We define
\bes{
\nm{g}_{m,\infty} = \max_{i = 0,\ldots,m} |g(x_i) |,\quad g \in C(I),
}
where $\{ x_i \}^{m}_{i=0}$ is the equispaced grid on $[-1,1]$, and
\bes{
\ip{f}{g}_{m,2} = \frac{2}{m+1} \sum^{m}_{i=0} f(x_i) \overline{g(x_i)},\quad f , g \in C(I).
}
We also let $\nm{\cdot}_{m,2} = \sqrt{\ip{\cdot}{\cdot}_{m,2}}$ be the corresponding discrete semi-norm. Note that $\ip{\cdot}{\cdot}_{m,2}$ is an inner product on $\bbP_n$ for any $m \geq n$, since a polynomial of degree $n$ cannot vanish at $m+1 \geq n+1$ distinct points unless it is the zero polynomial. Observe also that
\be{
\label{discnormbds}
\sqrt{2/(m+1)} \nm{f}_{m,\infty} \leq \nm{f}_{m,2} \leq \sqrt{2} \nm{f}_{m,\infty} \leq \sqrt{2} \nm{f}_{[-1,1],\infty},
}
for any $f \in C(I)$.

Finally, given a compact set $E \subset \bbC$, we write $B(E)$ for the set of functions that are continuous on $E$ and analytic in its interior. We also define $\nm{f}_{E,\infty} = \sup_{z \in E} |f(z)|$.

\subsection{The impossibility theorem}

Throughout this paper, we consider families of mappings
\bes{
\cR_m : C([-1,1]) \rightarrow C([-1,1]),
}
where, for each $m \geq 1$ and $f \in C([-1,1])$, $\cR_m(f)$ depends only on the values $\{f(x_i)\}^{m}_{i=0}$ of $f$ on the equispaced grid $\{x_i\}^{m}_{i=0}$. We define the (absolute) condition number of $\cR_m$ (in terms of the continuous and discrete uniform norms) as
\be{
\label{kappa-def}
\kappa(\cR_m) = \sup_{f \in C([-1,1])} \lim_{\delta \rightarrow 0^+} \sup_{\substack{h \in C([-1,1]) \\ 0 < \nm{h}_{m,\infty} \leq \delta }} \frac{\nmu{\cR_{m}(f+h) - \cR_{m}(f) }_{[-1,1],\infty}}{\nm{h}_{m,\infty}}.
}
We are now ready to state the impossibility theorem:

\thm{
[The impossibility theorem, \cite{platte2011impossibility}]
\label{t:impossibility}
 Let $E \subset \bbC$ be a compact set containing $[-1,1]$ in its interior and $\{ \cR_m \}^{\infty}_{m=1}$ be an approximation procedure based on equispaced grids of $m+1$ points such that, for some $C,\rho > 1$ and $1/2 < \tau \leq 1$, we have
\bes{
\nmu{f - \cR_m(f) }_{[-1,1],\infty} \leq C \rho^{-m^{\tau}} \nm{f}_{E,\infty},\quad \forall m \in \bbN,\ f \in B(E).
}
Then the condition numbers \R{kappa-def} satisfy
\bes{
\kappa(\cR_m) \geq \sigma^{m^{2 \tau - 1}}
}
for some $\sigma > 1$ and all sufficiently large $m$.
}

When $\tau = 1$, this implies that any approximation procedure that achieves exponential convergence at a geometric rate must also be exponentially ill-conditioned at a geometric rate. Furthermore, the best possible (and achievable \cite{adcock2019optimal}) rate of convergence of a stable method is root-exponential in $m$, i.e.\ the error decays like $\rho^{-\sqrt{m}}$ for some $\rho > 1$.

\subsection{Polynomial frame approximation}

We now describe polynomial frame approximation. As observed, this method was formalized in \cite{adcock2020approximating}, and is related to earlier works on Fourier extensions \cite{adcock2014parameter,adcock2014numerical,matthysen2017function,lyon2012fast,huybrechs2010fourier,bruno2007accurate,boyd2002fourier,pasquetti1996spectral}, and more generally, numerical approximations with frames \cite{adcock2019frames,adcock2020frames}.

Let $\gamma > 1$. Polynomial frame approximation uses orthogonal polynomials on an extended interval $[-\gamma,\gamma]$ to construct an approximation to a function over $[-1,1]$. In this paper, we use orthonormal Legendre polynomials, although we remark in passing that other orthogonal polynomials such as Chebyshev polynomials could also be employed. Note that an orthonormal basis on $[-\gamma,\gamma]$ fails to constitute a basis when restricted to the smaller interval $[-1,1]$. It forms a so-called \textit{frame} \cite{adcock2019frames,christensen2016introduction}, hence the name polynomial `frame' approximation.

Let $P_i(x)$ be the classical Legendre polynomial on $[-1,1]$, normalized so that $P_i(1) = 1$. Since $\int^{1}_{-1} |P_i(x)|^2 \D x = (i+1/2)^{-1}$, we define Legendre polynomial frame on $[-1,1]$ as $\{ \psi_i \}^{\infty}_{i=0}$, where the $i$th such function is given by
\be{
\label{psii-def}
\psi_i(x) = \sqrt{i+1/2} P_i  (x/\gamma  ) /\sqrt{\gamma},\quad x \in [-1,1].
}
Let $m,n \geq 0$ and consider a function $f \in C([-1,1])$. Our aim is to compute a polynomial approximation to $f$ of the form
\bes{
f \approx \hat{f} = \sum^{n}_{i=0} \hat{c}_i \psi_i \in \bbP_n
}
for suitable coefficients $\hat{c}_i$. It is natural to strive to do this via a least-squares fit, i.e.
\be{
\label{LScoeff}
\hat{c} = (\hat{c}_i)^{n}_{i=0}  \in \argmin{c \in \bbC^{n+1}} \nm{A c - b}_{2},
}
where
\be{
\label{Adef}
A =  \sqrt{2/(m+1)} (\psi_j(x_i))^{m,n}_{i,j=0} \in \bbC^{m \times n},\qquad b = \sqrt{2/(m+1)} (f(x_i))^{m}_{i=0} \in \bbC^m.
}
Note that $\sqrt{2/(m+1)}$ is simply a normalization factor that is included for convenience. Unfortunately, as described in \cite{adcock2020approximating}, this least-squares problem is ill-conditioned for large $n$ (even when $m \gg n$) due to the use of a frame rather than a basis \cite{adcock2019frames}. Therefore, we instead solve a suitably regularized least-squares problem. There are a number of different ways to do this, but, following previous works we consider an $\epsilon$-truncated Singular Value Decomposition (SVD).

Suppose that the least-squares matrix \R{Adef} has SVD $A = U \Sigma V^*$, where $\Sigma = \diag(\sigma_0,\ldots,\sigma_n) \in \bbR^{m \times n}$ is the diagonal matrix of singular values. Recall that the minimal $2$-norm solution $\hat{c}$ of \R{LScoeff} is given by
\bes{
\hat{c} = A^{\dag} b = V \Sigma^{\dag} U^* b,
}
where $\dag$ denotes the pseudoinverse. Given $\epsilon > 0$, we define $\Sigma^{\epsilon}$ as $\epsilon$-regularized version of $\Sigma$ as
\bes{
(\Sigma^{\epsilon} )_{ii} =  \begin{cases} \sigma_i & \sigma_i > \epsilon \\ 0 & \mbox{otherwise} \end{cases}  ,
}
and let $\Sigma^{\epsilon,\dag}$ be its pseudoinverse, i.e.
\bes{
(\Sigma^{\epsilon,\dag} )_{ii} = \begin{cases} 1/\sigma_i & \sigma_i > \epsilon \\ 0 & \mbox{otherwise} \end{cases}  .
}
Then we define the $\epsilon$-regularized approximation of \R{LScoeff} as $\hat{c}^{\epsilon} = V \Sigma^{\epsilon,\dag} U^* b$ and the corresponding approximation to $f$ as
\bes{
\hat{f}^{\epsilon} = \sum^{n}_{i=0} \hat{c}^{\epsilon}_i \psi_i.
}
With this in hand, we define the overall approximation procedure as the mapping
\be{
\label{poly-frame-approx-1}
\cP^{\epsilon,\gamma}_{m,n} : C([-1,1]) \rightarrow C([-1,1]),\ f \mapsto \hat{f}^{\epsilon} = \sum^{n}_{i=0} \hat{c}^{\epsilon}_{i} \psi_i,
}
where
\be{
\label{poly-frame-approx-2}
\hat{c}^{\epsilon} = (\hat{c}^{\epsilon}_i)^{n}_{i=0} = V \Sigma^{\epsilon,\dag} U^* b,\qquad b  = \sqrt{2/(m+1)} (f(x_i))^{m}_{i=0} .
}

\rem{
[Why not use orthogonal polynomials on the original interval]
\label{r:why-not-unit-interval}
The use of orthogonal polynomials on an extended interval may at first seem bizarre, since, as noted, the infinite collection of such functions no longer forms a basis when restricted to $[-1,1]$, but a frame. And even though the first $n+1$ such functions $\psi_0,\ldots,\psi_n$ constitute a basis for $\bbP_n$, they are extremely ill-conditioned as $n \rightarrow \infty$. To be precise, the condition number of their Gram matrix $G = (\ip{\psi_j}{\psi_i}_{[-1,1],2} )^{n}_{i,j=0}$ is exponentially-large in $n$. In turn, the matrix $A$ is also exponentially ill-conditioned in $n$, even when $m \gg n$. To understand why, observe that $A^* A = (\ip{\psi_j}{\psi_i}_{m,2} )^{n}_{i,j=0}$ is simply a discrete approximation to $G$.

Why then, do we not consider Legendre polynomials in $[-1,1]$? These constitute a perfectly well-conditioned basis, which means that the corresponding matrix $A$ would also be well-conditioned for $m \gg n$. Thus there is no need for regularization, and we may simply compute the least-squares fit by solving \R{LScoeff}. Note that this simply corresponds to $\cP^{\epsilon,\gamma}_{m,n}$ with $\epsilon = 0$ and $\gamma = 1$. The problem is that such an approximation, which we term \textit{polynomial least-squares approximation}, {behaves exactly as the} impossibility theorem (Theorem \ref{t:impossibility}) {predicts in the best case}. It is ill-conditioned if $m = o(n^2)$ as $n \rightarrow \infty$, and in particular, if $m \sim c n$ as $n \rightarrow \infty$, then the condition number of the mapping grows exponentially fast. On the other hand, with \textit{quadratic oversampling}, i.e.\ $m \sim c n^2$ as $n \rightarrow \infty$, the approximation is well-conditioned and its convergence rate is root-exponential in $m$ for all analytic functions. See \cite{adcock2019optimal} for a discussion on polynomial least-squares and the impossibility theorem. See also Remark \ref{r:least-squares-cond} below.
}

\subsection{Maximal behaviour of polynomials bounded at equispaced nodes}

As mentioned, both the impossibility theorem and the subsequent possibility theorem rely on estimates for the maximal behaviour of polynomials that are bounded at equispaced nodes. The former is based on a classical result due to Coppersmith and Rivlin concerning the maximal growth of a polynomial $p \in \bbP_n$ that is at most one at the $m+1$ equispaced nodes $\{ x_i \}^{m}_{i=0}$. Specifically, in \cite{coppersmith1992growth} they showed that there exist constants $\beta \geq \alpha > 1$ such that, if
\bes{
B(m,n) = \sup \{ \nm{p}_{[-1,1],\infty} : p \in \bbP_n,\ \nm{p}_{m,\infty} \leq 1 \},
}
then
\be{
\label{coppersmith-rivlin}
\alpha^{n^2/m} \leq B(m,n) \leq \beta^{n^2/m},\quad \forall 1 \leq n \leq m.
}

\rem{
[The condition number of polynomial least-squares approximation]
\label{r:least-squares-cond}
It is not difficult to show that the condition number of polynomial least-squares approximation $\cP_{m,n} = \cP^{0,1}_{m,n}$ satisfies
\bes{
B(m,n) \leq \kappa(\cP_{m,n}) \leq \sqrt{m+1} B(m,n),
}
where $B(m,n)$ is as in \R{coppersmith-rivlin} (this follows by setting $\gamma = 1$ and $\epsilon = 0$ in a result we show later, Lemma \ref{l:poly-frame-accuracy-stability}). Thus, \R{coppersmith-rivlin} immediately explains why this approximation is ill-conditioned when $m = o(n^2)$ as $n \rightarrow \infty$, and only well-conditioned when $m \sim c n^2$ (or faster) as $n \rightarrow \infty$.
}

As described above, the polynomial frame approximation is constructed via an $\epsilon$-truncated SVD. We will see later in \S \ref{s:acc-stab-poly-frame} that such truncation means that the approximation $\cP^{\epsilon,\gamma}_{m,n}(f)$ of a function $f$ belongs to (a subspace of) the set of polynomials
\bes{
P^{\epsilon,\gamma}_{m,n} = \left \{ p \in \bbP_n : \nm{p}_{[-\gamma,\gamma],2} \leq \nm{p}_{m,2} / \epsilon \right \} \subseteq \bbP_n
}
whose $L^2$-norm over the extended interval $[-\gamma,\gamma]$ is at most $1/\epsilon$ times larger than their discrete $2$-norm over the equispaced grid. In other words, the effect of regularization via the truncated SVD is to restrict the type of polynomial the approximation can take to one that does not grow too large on the extended interval.

After interchanging the $2$-norms for uniform norms, this observation motivates the study of the quantity
\be{
\label{C-def-general}
C(m,n,\gamma,\epsilon) = \sup \{ \nm{p}_{[-1,1],\infty} : p \in \bbP_n,\ \| p \|_{m,\infty} \leq 1,\ \| p \|_{[-\gamma,\gamma],\infty} \leq 1/\epsilon \}.
}
In \S \ref{s:acc-stab-poly-frame}, we show that the condition number of the polynomial frame approximation $\cP^{\epsilon,\gamma}_{m,n}$ satisfies
\bes{
\kappa(\cP^{\epsilon,\gamma}_{m,n}) \leq \sqrt{m+1} C(m,n,\gamma,\epsilon).
}
Notice that $C(m,n,\gamma,0) = B(m,n)$ and, in general,
\bes{
C(m,n,\gamma,\epsilon) \leq B(m,n),
}
where $B(m,n)$ is the as in \R{coppersmith-rivlin}. However, whereas $B(m,n)$ is only bounded as $m,n\rightarrow\infty$ in the quadratic oversampling regime (i.e.\ $m \sim c n^2$ for some $c>0$), for $C(m,n,\gamma,\epsilon)$ we show that linear oversampling is sufficient. Specifically:

\thm{
[Maximal behaviour of polynomials bounded at equispaced nodes and extended intervals]
\label{t:polynomial-inequality-main}
Let $0 < \epsilon \leq 1/\E$, $\gamma > 1$ and $n \geq \sqrt{\gamma^2-1} \log(1/\epsilon)$, and consider the quantity $C(m,n,\gamma,\epsilon)$ defined in \R{C-def-general}. Suppose that
\be{
\label{m-n-eps-poly-growth}
m \geq 36 n \log(1/\epsilon) / \sqrt{\gamma^2-1}.
}
Then
\bes{
C(m,n,\gamma,\epsilon) \leq c,
}
for some numerical constant $c > 0$. Specifically, $c$ can be taken as $c = 4 \beta + 3$ for $\beta$ as in \R{coppersmith-rivlin}.
}

This theorem is a direct consequence of a more general result (Theorem \ref{t:poly-inequality-general}) that we state and prove in \S \ref{s:maximal-behaviour-proof}. Note that the factor $36$ appearing in \R{m-n-eps-poly-growth} was chosen for convenience. Theorem \ref{t:poly-inequality-general} describes in general how a condition roughly of the form $m \geq c_1 n  \log(1/\epsilon) / \sqrt{\gamma^2-1}$ leads to a bound $C(m,n,\gamma,\epsilon) \leq h(c_1)$ for some function $h(c_1)$ depending on $c_1$.

\subsection{The motivations for considering polynomials on an extended interval}

Above we asserted that, by considering orthogonal polynomials on an extended interval and using regularization, the polynomial frame approximation constructs an approximation in a space within which the polynomials cannot grow too large on $[-\gamma,\gamma]$. As Theorem \ref{t:polynomial-inequality-main} makes clear, this prohibits such polynomials from behaving too wildly on $[-1,1]$ away from the equispaced grid, whenever $m$ scales linearly with $n$. Thus, using orthogonal polynomials on an extended interval allows for a reduction in the oversampling rate from quadratic in $n$ (as is the case for standard polynomial-least squares approximation -- see Remark \ref{r:least-squares-cond}) to linear in $n$.

On the other hand, by restricting the approximation space in this way, we potentially limit the ability of the scheme to accurately approximate analytic functions. In Theorem \ref{t:acc-cond-poly-frame}, we establish an error bound for the polynomial frame approximation that takes the form
\bes{
\nmu{f - \cP^{\epsilon',\gamma}_{m,n}(f)}_{[-1,1],\infty} \leq 2 c \sqrt{m+1} \inf_{p \in \bbP_n} \left \{ \nm{f - p}_{[-1,1],\infty} + (n+1)\epsilon \nm{p}_{[-\gamma,\gamma],\infty} \right \},
}
where $\epsilon' = \epsilon (n+1) / \sqrt{\gamma}$ (the reasons for this choice of $\epsilon'$ are discussed further below). This holds for any $c > 1$ and $f \in C([-1,1])$, provided $m$ and $n$ are chosen so that
\bes{
C(m,n,\gamma,\epsilon) \leq c.
}
Such a bound is very similar to those shown previously for both Fourier extensions \cite{adcock2014numerical} and polynomial frame approximations \cite{adcock2020approximating}, the main difference being the use of the $L^{\infty}$-norm instead of the $L^2$-norm. The key component of it is the best approximation term 
\bes{
\inf_{p \in \bbP_n} \left \{ \nm{f - p}_{[-1,1],\infty} + (n+1)\epsilon \nm{p}_{[-\gamma,\gamma],\infty} \right \}.
}
In other words, the effect of the truncated SVD regularization is to replace the classical best approximation error term
\bes{
\inf_{p \in \bbP_n} \nm{f - p}_{[-1,1],\infty},
}
(which arises in the case $\epsilon = 0$, i.e.\ standard polynomial least-squares approximation)
by one that also involves a term depending on $\epsilon$ multiplied by the norm of $p$ over the extended interval. As a result, the overall approximation error depends on how well $f$ can be approximated by a polynomial $p \in \bbP_n$ uniformly on $[-1,1]$ (the term $\nm{f - p}_{[-1,1],\infty}$) that does not grow too large on the extended interval $[-\gamma,\gamma]$ (the term $ \nm{p}_{[-\gamma,\gamma],\infty} $). 

Our main result, stated next, arises by bounding this best approximation term for functions that are analytic in sufficiently large complex regions.

\subsection{Main result: the possibility theorem}

We now state our main result (see \S \ref{s:main-thm-proof} for the proof). For this, we first recall the definition of the Bernstein ellipse with parameter $\theta > 1$:
\be{
\label{Bernstein-ellipse}
E_{\theta} = \left \{ \frac{z+z^{-1}}{2} : z \in \bbC,\ 1 \leq | z | \leq \theta \right \}.
}
A classical theorem in approximation theory states that any function $f \in B(E_{\theta})$ is approximated to exponential accuracy by polynomials. Specifically,
\be{
\label{poly-BA}
\inf_{p \in \bbP_n} \nm{f - p}_{[-1,1],\infty} \leq \frac{2}{\theta-1} \nm{f}_{E_{\theta},\infty} \theta^{-n},
}
(see also Lemma \ref{poly-approx-bounds} later). Our main results asserts that polynomial frame approximation can achieve a similar rate of decay in $n$, subject to linear oversampling in $m$. Specifically:

\thm{
[The possibility theorem]
\label{t:possibility-thm}
Let $0 < \epsilon \leq 1/\E$, $\gamma > 1$ and $n \geq \sqrt{\gamma^2-1} \log(1/\epsilon)$, and consider the polynomial frame approximation $\cP^{\epsilon',\gamma}_{m,n}$ defined in \R{poly-frame-approx-1}--\R{poly-frame-approx-2}, where
\be{
\label{m-n-relation}
m = 
\left \lceil 36 n \log(1/\epsilon) \big / \sqrt{\gamma^2-1} \right \rceil ,\qquad \epsilon' = \frac{\epsilon(n+1)}{\sqrt{\gamma}}.
}
Then the condition number of the mapping $\cP^{\epsilon',\gamma}_{m,n}$ satisfies
\bes{
 \kappa(\cP^{\epsilon',\gamma}_{m,n}) \leq c \sqrt{m+1} ,
}
where $c$ is as in Theorem \ref{t:polynomial-inequality-main}.
Moreover, if $E_{\theta}$ is a Bernstein ellipse with parameter
\be{
\label{theta-cond}
\theta > \gamma + \sqrt{\gamma^2-1}
}
then, for all $f \in B(E_{\theta})$,
\be{
\label{main-err-bd}
\begin{split}
\nmu{f - \cP^{\epsilon',\gamma}_{m,n}(f)}_{[-1,1],\infty} & \leq  c g(\theta,\gamma) \sqrt{m} \left ( \theta^{-n} + n \epsilon \right ) \nm{f}_{E_{\theta},\infty}
\\
& \leq  c g(\theta,\gamma) \sqrt{m} \left ( \rho^{1-m} + m \epsilon \right ) \nm{f}_{E_{\theta},\infty},
\end{split}
}
where $g(\theta,\gamma)$ depends on $\theta$ and $\gamma$ only and
\be{
\label{rho-main-thm}
\rho = \theta^{c_*},\qquad c_* =  \frac{\sqrt{\gamma^2-1} }{36 \log(1/\epsilon)} .
}
}

This result shows that polynomial frame approximation is well-conditioned when $n$ scales linearly with $m$ (specifically, \R{m-n-relation} holds), with its condition number being at worst $\ordu{\sqrt{m}}$ as $m \rightarrow \infty$. Moreover, for functions that are analytic in $E_{\theta}$ (note that this region contains the extended interval $[-\gamma,\gamma]$ in its interior, due to the condition \R{theta-cond}) its error decreases exponentially fast in $m$ down to the level $\ordu{m^{3/2} \epsilon}$. Recall that the rate $\theta^{-n}$ in \R{main-err-bd} is the same as in \R{poly-BA} for the best polynomial approximation of a function in $B(E_{\theta})$. Thus, one can achieve a near-optimal error decay rate in $n$ with only linear oversampling in $m$.

Overall, Theorem \ref{t:possibility-thm} provides a positive counterpart to the impossibility theorem (Theorem \ref{t:impossibility}). It is important emphasize that it does not violate Theorem \ref{t:impossibility}: indeed, exponential decay of the error is only guaranteed down to a finite accuracy. 

\rem{
[Varying $\epsilon$ with $n$]
\label{rem:varying-eps}
It is worth observing how the theorem changes if one strives to scale $\epsilon$ with $n$ so as to achieve exponential convergence of the error down to zero{, rather than exponential decay down to a finite level of accuracy}. If $\epsilon$ is chosen as $\epsilon = \theta^{-n}$ then the scaling between $m$ and $n$ becomes quadratic, since $\log(1/\epsilon) = n \log(\theta)$ in this case. When substituted into Theorem \ref{t:polynomial-inequality-main}, this implies  quadratic scaling of $m$ with $n$, and therefore root-exponential convergence {(to zero)} of the approximation with respect to $m$. {This is in precise agreement with the impossibility theorem.}
}
Another key aspect of Theorem \ref{t:possibility-thm} is the dependence on $\epsilon$ in \R{m-n-relation} and, in turn, the exponential rate \R{rho-main-thm}. Since $\epsilon > 0$ dictates the limiting accuracy of the approximation scheme, it is often desirable to choose $\epsilon$ close to machine epsilon $\epsilon_{\mathrm{mach}}$, which in IEEE double precision is roughly $\epsilon_{\mathrm{mach}} \approx 1.1 \times 10^{-16}$. Thus, the scaling $\log(1/\epsilon)$ -- which is proportional to the number of digits of accuracy desired -- is highly appealing. A scaling of, for example, $1/\epsilon$, would be meaningless for practical purposes.

Note that Theorem \ref{t:possibility-thm} does not say anything about the rate of the decay of the error for functions that are not analytic in a Bernstein ellipse $E = E_{\theta}$ that is large enough to contain the extended interval $[-\gamma,\gamma]$. We discuss the behaviour of the error for such functions in \S \ref{ss:lower-regularity}. On the other hand, this theorem also offers some insight into the effect of the choice of $\gamma$ on the approximation. Specifically, choosing a smaller $\gamma$ means that \R{theta-cond} holds for smaller values of $\theta$, thus the analyticity requirement $f \in B(E_{\theta})$ becomes less stringent. However, this also leads to a slower rate of exponential convergence in $m$, since $\rho$ is an increasing function of both $\gamma$ and $\theta$. We discuss this matter further in \S \ref{s:numerical}.

Finally, we remark that this theorem actually considers a polynomial frame approximation with parameter $\epsilon' = \epsilon(n+1)/\sqrt{\gamma}$ that grows linearly in $n$. The reason for this can be traced to the need to switch between the $L^2$-norm (or corresponding discrete seminorm) and the $L^{\infty}$-norm (or corresponding discrete seminorm) at various stages in the proof. See the proofs of Lemma \ref{C1-C2-as-C} and Theorem \ref{t:acc-cond-poly-frame} for the precise details. This choice of scaling is made to ensure the first term in the error bound decreases exponentially fast in $m$, which in turn follows from the linear relationship $m \approx 24 n \log(1/\epsilon)/\sqrt{\gamma^2-1}$. It is also possible to use $\epsilon$ as the truncation parameter rather than $\epsilon'$. Following much the same arguments, one can show that this choice results in a log-linear relationship between $m$ and $n$, which in turn leads to subexponential convergence of the form $\sigma^{-m / \log(m)}$ for large $m$, where, like $\rho$, $\sigma > 1$ depends on $\theta$, $\gamma$ and $\epsilon$.

\subsection{Error decay rates for functions of lower regularity}\label{ss:lower-regularity}

Theorem \ref{t:possibility-thm} only asserts exponential decay of the error for functions that are analytic in complex regions containing the extended interval $[-\gamma,\gamma]$. We now consider arbitrary analytic functions.  The following result shows that the error for such functions decays exponentially fast with the same rate $\theta^{-n}$, but only down to a larger tolerance.

\thm{
[Error decay for arbitrary analytic functions]
\label{t:possibility-slower-exp}
Consider the setup of Theorem \ref{t:possibility-thm}. Let $E_{\theta}$ be the Bernstein ellipse with parameter
\bes{
1 < \theta < \tau : =\gamma + \sqrt{\gamma^2-1}.
}
Then, for all $f \in B(E_{\theta})$,
\be{
\label{main-err-bd2}
\begin{split}
\nmu{f - \cP^{\epsilon',\gamma}_{m,n}(f)}_{[-1,1],\infty} & \leq  c g(\theta,\gamma) \sqrt{m} \left ( \theta^{-n} + n \epsilon^{\frac{\log(\theta)}{\log(\tau)}} \right ) \nm{f}_{E_{\theta},\infty}
\\
& \leq  c g(\theta,\gamma) \sqrt{m} \left ( \rho^{1-m} + m \epsilon^{\frac{\log(\theta)}{\log(\tau)}} \right ) \nm{f}_{E_{\theta},\infty},
\end{split}
}
where $g(\theta,\gamma)$ and $\rho$ are also as in Theorem \ref{t:possibility-thm}.
}

This result shows decrease of the error exponentially fast down to roughly $\epsilon^{\frac{\log(\theta)}{\log(\tau)}}$, i.e.\ some fractional power of $\epsilon$ depending on the relative sizes of $\theta$ and $\tau = \gamma + \sqrt{\gamma^2-1}$. 
It raises an immediate question: what happens after such an accuracy level is reached? As we show later through an extended impossibility theorem, we cannot expect exponential decay down to $\epsilon$ in general. However, we now show that superalgebraic decay -- i.e.\ faster than any power of $m^{-k}$ -- is indeed possible down to this level. 

We do this by first noting that an analytic function is infinitely smooth, and then by establishing an error bound for functions that are $k$-times continuously differentiable.
Specifically, in the following result we consider the space $C^{k}([-1,1])$ of functions that are $k$-times continuously differentiable on $[-1,1]$. We define the norm on this space as
\bes{
\nm{f}_{C^k([-1,1])} = \max_{j = 0,\ldots,k} \nmu{f^{(j)}}_{[-1,1],\infty}.
}

\thm{
[Error decay for $C^k$ functions]
\label{t:possibility-algebraic}
Consider the setup of Theorem \ref{t:possibility-thm}. Then, for all $k \in \bbN$ and $f \in C^k([-1,1])$,
\bes{
\nmu{f - \cP^{\epsilon',\gamma}_{m,n}(f)}_{[-1,1],\infty} \leq c g(k,\gamma) \sqrt{m} \left ( n^{-k} + n \epsilon \right ) \nm{f}_{C^k([-1,1])},
}
where  $g(k,\gamma)$ depends on $k$ and $\gamma$ only.
}

Proofs of Theorems \ref{t:possibility-slower-exp} and \ref{t:possibility-algebraic} can be found in \S \ref{s:main-thm-proof}.
Note that all these observations about the rate of error decay are seen in practice in numerical examples. We present a series of experiments confirming these results in \S \ref{s:numerical}.

We remark in passing that superalgebraic decay is slower than root-exponential decay, which is the best possible stipulated by the impossibility theorem for analytic function approximation. Whether or not polynomial frame approximation exhibits root-exponential decay in $m$ after the breakpoint $\epsilon^{\frac{\log(\theta)}{\log(\tau)}}$ is an open problem.

\subsection{An extended impossibility theorem}

Theorems \ref{t:possibility-thm} and \ref{t:possibility-slower-exp} show that polynomial frame approximation can achieve roughly the same exponential rate $\theta^{-n}$ as the best polynomial approximation for functions in $B(E_{\theta})$ when subject to linear oversampling. However, it only maintains this rate down to roughly $\epsilon$ for sufficiently large $\theta$. We now ask whether or not there exists an approximation scheme that can perform better than this: namely, whether an exponential rate $\theta^{-n}$ down to $\epsilon$ can be attained for any $\theta > 1$ in the linear oversampling regime. The following extended impossibility theorem shows that the answer to this question is no. 

\thm{[Extended impossibility theorem]
\label{t:impossibility-extended}
Let $\{ \cR_m \}^{\infty}_{m=1}$ be an approximation procedure based on equispaced grids of $m+1$ points such that, for some $c > 0$, $C,\theta^* > 1$, $0 < \epsilon \leq (4C)^{-2}$ and $1/2 < \tau \leq 1$, we have
\be{
\label{extended-error-bound}
\nmu{f - \cR_m(f) }_{[-1,1],\infty} \leq C \left ( \theta^{-c m^{\tau}} + \epsilon \right ) \nm{f}_{E_{\theta},\infty},\quad 
}
for all $m \in \bbN$, $f \in B(E_{\theta})$ and $1 < \theta \leq \theta^*$.
Then the condition numbers \R{kappa-def} satisfy
\bes{
\kappa(\cR_m) \geq \sigma^{m^{2 \tau - 1}}
}
for some $\sigma > 1$ and all sufficiently large $m$.
}

The proof of this theorem can be found in \S \ref{s:numerical}. It follows essentially the same steps as that of the impossibility theorem.

This theorem has several consequences. First, it extends the impossibility theorem by showing that exactly the same relationship between fast error decay and conditioning holds even when the overall error decreases only down to a constant tolerance $\epsilon$. To do this, it makes the stronger assumption that the scheme yields exponential decay for all analytic functions -- including those with singularities arbitrarily close to $[-1,1]$ -- with the rate of exponential decay being dependent on the size of the region of analyticity. Specifically, $f \in B(E_{\theta})$ implies a term of the form $\theta^{-c m^{\tau}}$ for \textit{all} $1 < \theta \leq \theta^*$.

Second, it implies the following. Any well-conditioned method must either (i) yield root-exponential decrease of the error down to $\epsilon$, i.e. $\theta^{-c \sqrt{m}} + \epsilon$; or (ii) fail to yield exponential convergence for all analytic functions, i.e.\ \R{extended-error-bound} holds with $\tau = 1$ only for $\theta^{**} \leq \theta \leq \theta^*$ for some $1 < \theta^{**} \leq \theta^*$. As discussed previously, (ii) is exactly how polynomial frame approximation behaves, up to small algebraic factors in $m$. In other words, it is nearly optimal.

\rem{
[Varying $\gamma$ with $n$]

Recall that Theorem \ref{t:possibility-thm} only asserts exponential decay down to $\epsilon$ for functions that are analytic Bernstein ellipses containing the interval $[-\gamma,\gamma]$. A natural idea is therefore to decrease $\gamma$ with $n$ so that, for any fixed, analytic function, the interval $[-\gamma,\gamma]$ is included in its region of analyticity for all large $n$.

To determine a suitable scaling for $\gamma$, one can use similar ideas to those employed in so-called \textit{mapped polynomial} spectral methods \cite{adcock2016mapped,don1997accuracy,boyd1989chebyshev,hale2008new,kosloff1993modified}: namely, choose $\gamma$ to match the terms in the error bound \R{main-err-bd}. Ignoring the term $n$ (for simplicity), we therefore solve the equation
$(\gamma + \sqrt{\gamma^2-1})^{-n} = \epsilon$
with respect to $\gamma$, which yields
\be{
\label{gamma-scale}
\gamma = \gamma(n,\epsilon) = \frac{\epsilon^{1/n} + \epsilon^{-1/n}}{2}.
}
Observe that $\gamma \rightarrow 1^+$ as $n \rightarrow \infty$ for fixed $n$. By combining Theorems \ref{t:possibility-thm} (for large $n$, since $\gamma < \theta$ for all large $n$) and \ref{t:possibility-slower-exp} (for small $n$, noting that $\epsilon^{\frac{\log(\tau)}{\log(\theta)}} = \theta^{-n}$ with this choice of $\gamma$), one can show that
\be{
\label{gamma-scale-err-bound}
\nmu{f - \cP^{\epsilon',\gamma}_{m,n}(f)}_{[-1,1],\infty} \leq  c g(\theta,\gamma) \sqrt{m} \left ( \theta^{-n} + n \epsilon \right ) \nm{f}_{E_{\theta},\infty},
}
for $f \in B(E_{\theta})$ and \textit{any} $\theta > 1$. Yet, scaling $\gamma$ as in \R{gamma-scale} causes the relation between $m$ and $n$ to become quadratic. Indeed,
$\sqrt{\gamma^2 - 1} = \log(1/\epsilon)/n + \ord{(\log(1/\epsilon)/n)^3}$ as $n \rightarrow \infty$,
and therefore the condition \R{m-n-eps-poly-growth} becomes $m > 36 n^2$
for all large $n$. Of course, this is to be expected in view of Theorem \ref{t:impossibility-extended}, since the error bound \R{gamma-scale-err-bound} holds for any $\theta > 1$.
}

\section{Accuracy and conditioning of $\cP^{\epsilon,\gamma}_{m,n}$}\label{s:acc-stab-poly-frame}

The next three sections develop the proofs of Theorems \ref{t:polynomial-inequality-main}--\ref{t:impossibility-extended}. We commence in this section by analyzing the error and condition number of the polynomial frame approximation $\cP^{\epsilon,\gamma}_{m,n}$. Our analysis follows similar lines to that of \cite{adcock2020frames} and, in particular, \cite{adcock2020approximating}, the main difference being that we work in the $L^{\infty}$-norm, rather than the $L^2$-norm.

\subsection{Reformulation in terms of singular vectors}

Let $v_0,\ldots,v_{n} \in \bbC^n$ be the right singular vectors of the matrix $A$ defined in \R{Adef}. Then, to each such vector, we can associate a polynomial in $\bbP_n$:
\bes{
\xi_i = \sum^{n}_{j=0} (v_i)_j \psi_j \in \bbP_n,\quad i =0,\ldots,n.
}
Here the $\psi_j$ are as in \R{psii-def}.
Since these functions are orthonormal on $[-\gamma,\gamma]$ and the $v_i$ are orthonormal vectors, the functions $\xi_i$ are themselves orthonormal on $[-\gamma,\gamma]$:
\bes{
\ip{\xi_i}{\xi_j}_{[-\gamma,\gamma],2} = v^*_j v_i = \delta_{ij},\quad i,j = 0,\ldots,n.
}
However, the functions $\xi_i$ are also orthogonal with respect to the discrete inner product $\ip{\cdot}{\cdot}_{m,2}$. Indeed, since $v_j$ and $v_k$ are singular vectors, we have
\eas{
\ip{\xi_j}{\xi_k}_{m,2} = \frac{2}{m+1} \sum^{m}_{i=0} \xi_j(x_i) \overline{\xi_k(x_i)} 
& = \frac{2}{m+1} \sum^{m}_{i=0} \sum^{n}_{s=0} \sum^{n}_{t=0} (v_j)_s \psi_s(x_i) \overline{(v_k)_t} \overline{\psi_t(x_i)}
\\
& = v^*_k A^* A v_j 
\\
& = \sigma^2_j \delta_{jk}.
}
With this in hand, we now define the subspace
\bes{
\bbP^{\epsilon,\gamma}_{m,n} = \spn \{ \xi_i : \sigma_i > \epsilon \} \subseteq \bbP_n.
}
We note in passing that this space coincides with $\bbP_n$ whenever $m \geq n$ and $\sigma_{\min} = \sigma_n > \epsilon$. In particular, $\bbP^{0,\gamma}_{m,n} = \bbP_n$ for $m \geq n$. On the other hand, when $\epsilon > 0$ we have $\bbP^{\epsilon,\gamma}_{m,n} \subseteq P^{\epsilon,\gamma}_{m,n}$,
where
\bes{
P^{\epsilon,\gamma}_{m,n} = \left \{ p \in \bbP_n : \nm{p}_{[-\gamma,\gamma],2} \leq \nm{p}_{m,2} / \epsilon \right \} \subseteq \bbP_n,
}
is, as before, the set of polynomials of degree at most $n$ whose $L^2$-norm over the extended interval $[-\gamma,\gamma]$ is at most $1/\epsilon$ times larger than their discrete $2$-norm over the equispaced grid. To see why this holds, let $p = \sum_{\sigma_i > \epsilon} c_i \xi_i \in \bbP^{\epsilon,\gamma}_{m,n}$. Then, by the double orthogonality of the $\xi_i$,
\bes{
\nm{p}^2_{[-\gamma,\gamma,2]} = \sum_{ \sigma_i > \epsilon} |c_i|^2 \leq \frac{1}{\epsilon^2} \sum_{ \sigma_i > \epsilon} \sigma^2_i |c_i|^2 = \frac{\nm{p}^2_{m,2}}{\epsilon^2}.
}
Hence $p \in P^{\epsilon,\gamma}_{m,n}$, as required.

Finally, observe that the polynomial frame approximation $\cP^{\epsilon,\gamma}_{m,n}(f)$ of $f \in C([-1,1])$ belongs to the space 
$\bbP^{\epsilon,\gamma}_{m,n}$. In fact, it is the orthogonal projection onto this space with respect to the discrete inner product $\ip{\cdot}{\cdot}_{m,2}$. Therefore, by orthogonality, we may write
\be{
\label{Feps-in-terms-of-SVs}
\cP^{\epsilon,\gamma}_{m,n}(f) = \sum_{\sigma_i > \epsilon} \frac{\ip{f}{\xi_i}_{m,2}}{\sigma^2_i} \xi_i,\qquad f \in C([-1,1]).
}

\subsection{Accuracy and conditioning up to constants}

We can now examine the accuracy and conditioning of $\cP^{\epsilon,\gamma}_{m,n}$. We first require the following:

\lem{
\label{l:poly-frame-accuracy-stability}
Let $\epsilon \geq 0$, $\gamma \geq 1$ and $\cP^{\epsilon,\gamma}_{m,n}$ the corresponding polynomial frame approximation. Then, for any $f \in C([-1,1])$,
\be{
\label{acc-bound}
\nmu{f - \cP^{\epsilon,\gamma}_{m,n}(f)}_{[-1,1],\infty} \leq \left ( 1 + \sqrt{m+1} C_1 \right )\nmu{f - p}_{[-1,1],\infty} + C_2 \epsilon \nmu{p}_{[-\gamma,\gamma],\infty},\quad \forall p \in \bbP_n,
}
where
\be{
\label{C1-C2-def}
\begin{split}
C_1 & = C_1(m,n,\gamma,\epsilon) = \sup \{ \nmu{p}_{[-1,1],\infty} : p \in \bbP^{\epsilon,\gamma}_{m,n},\ \nm{p}_{m,\infty} \leq 1 \},
\\
C_2 & = C_2(m,n,\gamma,\epsilon) =  \{ \nmu{p -\cP^{\epsilon,\gamma}_{m,n}(p)}_{[-1,1],\infty} : p \in \bbP_{n},\  \nm{p}_{[-\gamma,\gamma],\infty} \leq \epsilon^{-1} \}.
\end{split}
}
Moreover, the condition number $\kappa(\cP^{\epsilon,\gamma}_{m,n})$ satisfies
\be{
\label{stab-bound}
C_1 \leq \kappa(\cP^{\epsilon,\gamma}_{m,n}) \leq \sqrt{m+1} C_1.
}
}

These constants are interpreted as follows. The first, $C_1$, measures how large an element of the space $\bbP^{\epsilon,\gamma}_{m,n}$ can be uniformly on $[-1,1]$ relative to its discrete uniform norm over the grid. The second, $C_2$, examines the effect of the $\epsilon$-truncation. Specifically, it measures the error in reconstructing a polynomial $p \in \bbP_n$ relative to its size on the extended interval. Note, in particular, that $C_2 = 0$ when $\epsilon = 0$ and $m \geq n$. We also have $C_1(m,n,\gamma,0) = B(m,n)$ in this case, where $B(m,n)$ is as in \R{coppersmith-rivlin}.

\prf{
Let $f \in C([-1,1])$ and $p \in \bbP_n$. Then
\eas{
\| f  -  \cP^{\epsilon,\gamma}_{m,n}(f) \|_{[-1,1],\infty} 
& \leq \nmu{f - p}_{[-1,1],\infty} + \nmu{\cP^{\epsilon,\gamma}_{m,n}(f) - \cP^{\epsilon,\gamma}_{m,n}(p)}_{[-1,1],\infty} + \nmu{p - \cP^{\epsilon,\gamma}_{m,n}(p)}_{[-1,1],\infty}
\\
& \leq \nmu{f - p}_{[-1,1],\infty} + C_1 \nmu{\cP^{\epsilon,\gamma}_{m,n}(f) - \cP^{\epsilon,\gamma}_{m,n}(p)}_{m,\infty} + C_2 \epsilon \nmu{p}_{[-\gamma,\gamma],\infty}.
}
Here, in the second step, we used the fact that $\cP^{\epsilon,\gamma}_{m,n}$ is linear and its range is the space $\bbP^{\epsilon,\gamma}_{m,n}$. Now recall that $\cP^{\epsilon,\gamma}_{m,n}$ is an orthogonal projection with respect to $\ip{\cdot}{\cdot}_{m,2}$. Thus, using \R{discnormbds},
\eas{
\nmu{\cP^{\epsilon,\gamma}_{m,n}(f) - \cP^{\epsilon,\gamma}_{m,n}(p)}_{m,\infty} 
& \leq \sqrt{(m+1)/2} \nmu{\cP^{\epsilon,\gamma}_{m,n}(f) - \cP^{\epsilon,\gamma}_{m,n}(p)}_{m,2} 
\\
& \leq \sqrt{(m+1)/2} \nmu{f - p}_{m,2}
\\
& \leq \sqrt{m+1} \nmu{f-p}_{m,\infty} 
 \leq \sqrt{m+1} \nmu{f-p}_{[-1,1],\infty}.
}
Substituting this into the previous expression now gives \R{acc-bound}.

For the second result, we use the fact that $\cP^{\epsilon,\gamma}_{m,n}$ is linear once more to write
\bes{
\kappa(\cP^{\epsilon,\gamma}_{m,n}) = \sup_{\substack{f \in C([-1,1]) \\ \nm{f}_{m,\infty} \neq 0}} \frac{\nmu{\cP^{\epsilon,\gamma}_{m,n}(f) }_{[-1,1],\infty} }{\nm{f}_{m,\infty}}.
}
Notice that $\cP^{\epsilon,\gamma}_{m,n}(p) = p$ for all $p \in \bbP^{\epsilon,\gamma}_{m,n}$. Hence
\bes{
\kappa(\cP^{\epsilon,\gamma}_{m,n}) \geq \sup_{\substack{p \in \bbP^{\epsilon,\gamma}_{m,n} \\ \nm{p}_{m,\infty} \neq 0}} \frac{\nm{p}_{[-1,1],\infty} }{\nm{p}_{m,\infty}} = C_1,
}
which gives the lower bound in \R{stab-bound}. On the other hand, using \R{discnormbds} and the fact that $\cP^{\epsilon,\gamma}_{m,n}$ is an orthogonal projection once more, we have
\eas{
\nmu{\cP^{\epsilon,\gamma}_{m,n}(f) }_{[-1,1],\infty} &\leq C_1 \nmu{\cP^{\epsilon,\gamma}_{m,n}(f) }_{m,\infty} 
 \leq C_1 \sqrt{(m+1)/2} \nmu{\cP^{\epsilon,\gamma}_{m,n}(f) }_{m,2} 
\\
& \leq  C_1 \sqrt{(m+1)/2} \nmu{f}_{m,2}
\\
& \leq C_1 \sqrt{m+1} \nm{f}_{m,\infty}.
}
This gives the upper bound in \R{stab-bound}.
}

\subsection{Bounding the constants}

The next step is to estimate the constants $C_1$ and $C_2$ appearing in this lemma.

\lem{
\label{C1-C2-as-C}
Consider the setup of the previous lemma. Then the constants $C_1$ and $C_2$ defined in \R{C1-C2-def} satisfy
\bes{
C_1 \leq C \left(m,n,\gamma,\epsilon / (\sqrt{2} c_{n,\gamma})  \right ),
}
and
\bes{
C_2 \leq \sqrt{\gamma(m+1)} \cdot C \left ( m , n , \gamma , \sqrt{m+1} \epsilon / (\sqrt{2} c_{n,\gamma}) \right ),
}
where $C$ is as in \R{C-def-general} and $c_{n,\gamma} = \frac{n+1}{\sqrt{2\gamma}}$.
}
\prf{
Consider $C_1$ first. Let $p \in \bbP^{\epsilon,\gamma}_{m,n}$ with $\nm{p}_{m,\infty} \leq 1$. Then we can write $p = \sum_{\sigma_i > \epsilon} c_i \xi_i$ and, using the orthogonality of the $\xi_i$, we see that
\bes{
\nm{p}^2_{m,2} = \sum_{\sigma_i > \epsilon} \sigma^2_i |c_i|^2 ,\qquad \nm{p}^2_{[-\gamma,\gamma],2} = \sum_{\sigma_i > \epsilon} |c_i|^2.
}
Therefore
\bes{
 \nm{p}_{[-\gamma,\gamma],2} \leq \epsilon^{-1} \nm{p}_{m,2} \leq \sqrt{2} \epsilon^{-1} \nm{p}_{m,\infty} \leq \sqrt{2} \epsilon^{-1}.
}
We also recall the following inequality over $\bbP_n$:
\be{
\label{poly-inf-2}
\nm{q}_{[-\gamma,\gamma],\infty} \leq c_{n,\gamma} \nm{q}_{[-\gamma,\gamma],2},\quad \forall q \in \bbP_n.
}
This can be show directly by recalling that the classical Legendre polynomial $P_i(x)$ attains its maximum at $x = 1$ and takes value $P_i(1) = 1$. Hence, writing $q \in \bbP_n$ as $q = \sum^{n}_{i=0} c_i \psi_i$ and recalling \R{psii-def}, we get
\eas{
\nm{q}_{[-\gamma,\gamma],\infty} & \leq \sum^{n}_{i=0} |c_i| \sqrt{\frac{i+1/2}{\gamma}} 
\leq \left ( \sum^{n}_{i=0} |c_i|^2 \right )^{1/2} \left ( \sum^{n}_{i=0} \frac{i+1/2}{\gamma} \right )^{1/2} 
 = c_{n,\gamma} \nm{q}_{[-\gamma,\gamma],2},
}
which establishes \R{poly-inf-2}.
Therefore, since $p \in \bbP^{\epsilon,\gamma}_{m,n} \subseteq \bbP_n$, we get
\bes{
 \nm{p}_{[-\gamma,\gamma],\infty} \leq c_{n,\gamma} \nm{p}_{[-\gamma,\gamma],2} \leq \sqrt{2} c_{n,\gamma} \epsilon^{-1}.
}
We deduce that
\eas{
C_1 &\leq \sup \{ \| p \|_{[-1,1],\infty} : p \in \bbP^{\epsilon,\gamma}_{m,n},\ \nm{p}_{m,\infty} \leq 1,\ \nm{p}_{[-\gamma,\gamma],\infty} \leq \sqrt{2} c_{n,\gamma} \epsilon^{-1} \} 
\\
& = C(m,n,\gamma,\epsilon / (\sqrt{2} c_{n,\gamma} ) ),
}
which gives the first result.

We now consider $C_2$. Let $p \in \bbP_n$ with $\nm{p}_{[-\gamma,\gamma],\infty} \leq \epsilon^{-1}$. Since $p \in \bbP_n$, we may write
\bes{
p = \sum^{n}_{i = 0} \frac{\ip{p}{\xi_i}_{m,2}}{\sigma^2_i} \xi_i,
}
and using \R{Feps-in-terms-of-SVs}, we may also write
\bes{
\cP^{\epsilon,\gamma}_{m,n}(p) = \sum_{\sigma_i > \epsilon} \frac{\ip{p}{\xi_i}_{m,2}}{\sigma^2_i} \xi_i.
}
Therefore
\bes{
p - \cP^{\epsilon,\gamma}_{m,n}(p) = \sum_{\sigma_i \leq \epsilon}  \frac{\ip{p}{\xi_i}_{m,2}}{\sigma^2_i} \xi_i.
}
Using the fact that the $\xi_i$ are orthonormal over $[-\gamma,\gamma]$, we deduce that
\be{
\label{p-minus-Feps-p-1}
\nmu{p - \cP^{\epsilon,\gamma}_{m,n}(p)}^2_{[-\gamma,\gamma],2} =  \sum_{\sigma_i \leq \epsilon} \frac{|\ip{p}{\xi_i}_{m,2} |^2}{\sigma^4_i} \leq \sum^{n}_{i = 0} \frac{|\ip{p}{\xi_i}_{m,2} |^2}{\sigma^4_i} = \nm{p}^2_{[-\gamma,\gamma],2},
}
and, using the fact that the $\xi_i$ are orthogonal with respect to the discrete semi-inner product $\ip{\cdot}{\cdot}_{m,2}$, we see that
\be{
\label{p-minus-Feps-p-2}
\nmu{p - \cP^{\epsilon,\gamma}_{m,n}(p)}^2_{m,2} =  \sum_{\sigma_i \leq \epsilon} \frac{|\ip{p}{\xi_i}_{m,2} |^2}{\sigma^2_i} \leq \epsilon^2  \sum_{\sigma_i \leq \epsilon} \frac{|\ip{p}{\xi_i}_m |^2}{\sigma^4_i} \leq \epsilon^2 \nm{p}^2_{[-\gamma,\gamma],2}.
}
Now observe that we can write
\bes{
C_2 = \max \{ \nmu{q}_{[-1,1],\infty} : q \in \cA \},\qquad \cA = \left \{ q : q = p - \cP^{\epsilon,\gamma}_{m,n}(p),\ p \in \bbP_n,\ \nmu{p}_{[-\gamma,\gamma],\infty} \leq \epsilon^{-1} \right \}.
}
Let $q = p - \cP^{\epsilon,\gamma}_{m,n}(p) \in \cA$. Then $q \in \bbP_n$ and, due to \R{poly-inf-2} and \R{p-minus-Feps-p-1},
\bes{
\nmu{q}_{[-\gamma,\gamma],\infty} \leq c_{n,\gamma} \nmu{q}_{[-\gamma,\gamma],2} \leq c_{n,\gamma}  \nmu{p}_{[-\gamma,\gamma],2} \leq \sqrt{2 \gamma} c_{n,\gamma}  \nmu{p}_{[-\gamma,\gamma],\infty} \leq \sqrt{2 \gamma} c_{n,\gamma}  /\epsilon.
}
Also, by \R{discnormbds}, \R{p-minus-Feps-p-2} and the fact that $\nm{p}_{[-\gamma,\gamma],\infty} \leq \epsilon^{-1}$,
\eas{
\nmu{q}_{m,\infty}& \leq \sqrt{(m+1)/2} \nmu{q}_{m,2}  \leq \sqrt{(m+1)/2} \epsilon \nm{p}_{[-\gamma,\gamma],2}  \leq \sqrt{(m+1)/2} \sqrt{2 \gamma} \epsilon \nm{p}_{[-\gamma,\gamma],\infty},
}
and therefore $ \nmu{q}_{m,\infty} \leq \sqrt{(m+1)/2} \sqrt{2 \gamma}$.
Hence
\bes{
q \in \cB : = \left \{ q \in \bbP_n : \nmu{q}_{[-\gamma,\gamma],\infty} \leq \sqrt{2 \gamma} c_{n,\gamma} / \epsilon,\ \nmu{q}_{m,\infty} \leq \sqrt{(m+1)/2} \sqrt{2 \gamma} \right \},
}
which implies that $C_2 \leq \max \{ \nmu{q}_{[-1,1],\infty} : q \in \cB \}$ and, after renormalizing,
\bes{
C_2 \leq \sqrt{(m+1)/2} \sqrt{2 \gamma}  \max \left \{ \nmu{p}_{[-\gamma,\gamma],\infty} : p \in \bbP_{n},\ \nm{p}_{m,\infty} \leq 1,\ \nm{p}_{[-\gamma,\gamma],\infty} \leq \frac{c_{n,\gamma} }{\sqrt{(m+1)/2} \epsilon} \right \}.
}
This gives the second result.
}

\subsection{Main result on accuracy and conditioning}

We now summarize these two lemmas with the following theorem:

\thm{
[Accuracy and conditioning of polynomial frame approximation]
\label{t:acc-cond-poly-frame}
Let $\epsilon > 0$, $\gamma \geq 1$, $c > 1$ and $m , n \geq 1$ be such that
\be{
\label{C-cond-main-acc-stab}
C(m,n,\gamma,\epsilon ) \leq c.
}
Then the polynomial frame approximation $\cP^{\epsilon',\gamma}_{m,n}$ with $\epsilon' = \epsilon (n+1) / \sqrt{\gamma}$ satisfies
\bes{
\kappa(\cP^{\epsilon',\gamma}_{m,n}) \leq c\sqrt{m+1},
}
and, for any $f \in C([-1,1])$,
\bes{
\nmu{f - \cP^{\epsilon',\gamma}_{m,n}(f)}_{[-1,1],\infty} \leq 2 c \sqrt{m+1} \inf_{p \in \bbP_n} \left \{ \nm{f - p}_{[-1,1],\infty} + (n+1)\epsilon \nm{p}_{[-\gamma,\gamma],\infty} \right \}.
}
}
\prf{
Observe that $C(m,n,\gamma,\epsilon)$ is a decreasing function of $\epsilon$. Hence, by this, Lemma \ref{C1-C2-as-C} and the fact that $\epsilon' = \sqrt{2} c_{n,\gamma} \epsilon$,
\bes{
C_1(m,n,\gamma,\epsilon') \leq C(m,n,\gamma,\epsilon' / (\sqrt{2} c_{n,\gamma}) ) = C(m,n,\gamma,\epsilon),
}
and
\bes{
C_{2}(m,n,\gamma,\epsilon') \leq \sqrt{\gamma(m+1)}C(m,n,\gamma,\sqrt{m+1} \epsilon' / (\sqrt{2} c_{n,\gamma}) ) \leq \sqrt{\gamma(m+1)} C(m,n,\gamma,\epsilon).
}
Therefore, the condition \R{C-cond-main-acc-stab} implies that
\bes{
C_1(m,n,\gamma,\epsilon' ) \leq c,\quad C_2(m,n,\gamma,\epsilon' ) \leq c \sqrt{\gamma(m+1)}.
}
We now apply Lemma \ref{l:poly-frame-accuracy-stability} to get that $\kappa(\cP^{\epsilon',\gamma}_{m,n}) \leq c \sqrt{m+1}$. For the error bound, we have
\eas{
\nmu{f - \cP^{\epsilon',\gamma}_{m,n}(f)}_{[-1,1],\infty} \leq & \left ( 1 + c\sqrt{m+1}  \right ) \nm{f-p}_{[-1,1],\infty} + c \sqrt{\gamma(m+1)} \epsilon' \nm{p}_{[-\gamma,\gamma],\infty}
\\
\leq & \left ( 1 + c\sqrt{m+1}  \right )  \left (  \nm{f-p}_{[-1,1],\infty} + (n+1)\nm{p}_{[-\gamma,\gamma],\infty} \right ),
}
where in the second step we used the definition of $\epsilon'$.
The result now follows, since $1 + c\sqrt{m+1} \leq 2 c \sqrt{m+1}$.
}

\section{Maximal behaviour of polynomials bounded at equispaced nodes}\label{s:maximal-behaviour-proof}

In this section, we prove Theorem \ref{t:polynomial-inequality-main}.

\subsection{Pointwise Markov inequality }

We first require a pointwise Markov inequality. The following lemma may be viewed as a generalization of the pointwise Bernstein inequality for algebraic polynomials $p \in \bbP_n$:
\be{\lb{B}
         |p'(x)| \le \frac{n}{\sqrt{1-x^2}} \|p\|_{[-1,1],\infty}, \qquad |x| < 1,
}
to higher derivatives. It appeared in \cite{konyagin2021stable} in a slightly less general form. We provide essentially the same proof for completeness.

\begin{lemma} \lb{M}
For any $k,n \in \bbN$ and $\delta \in (0,1)$ such that  $k < n \sqrt{(1-\delta^2)/2}$, 
and for any polynomial $p \in \bbP_n$, we have
\be{\lb{k}
     |p^{(k)}(x)| \le \frac{1.251 n^k}{(1-x^2)^{k/2}} \|p\|_{[-1,1],\infty}, \qquad |x| \le \delta.
}
\end{lemma}

Note that we cannot get \rf[k] just by iterating the Bernstein inequality \rf[B]. Such iterated use produces a much weaker result
$$
|p^{(k)}(x)| \le \frac{n^k k^{k/2}}{(1-x^2)^{k/2}} \|p\|_{[-1,1],\infty}, \qquad |x| < 1.
$$

\prf{
We will use the following known results. First, Shaeffer and Duffin \cite{schaeffer1938some} proved that, for any $k,n \in \bbN$ and $p \in \bbP_n$,
\be{\label{sd1}
          |p^{(k)}(x)| \le D_{n,k}(x) \|p\|_{[-1,1],\infty}, \quad |x| < 1,
}
where
\be{\lb{D}
      D_{n,k}(x) = \left|(\cos n \arccos x)^{(k)} + i (\sin n \arccos x)^{(k)}\right|.
}
Second, Shadrin \cite{shadrin2004twelve} derived the explicit expression for $D_{n,k}(\cdot)$. Specifically,
\be{\lb{Df}
   \frac{1}{n^2}\big(D_{n,k}(x)\big)^2
= \sum_{m=0}^{k-1} \frac{b_{m,n}}{(1-x^2)^{k+m}},
}
where
\ea{
     b_{m,n}
 &= c_{m,k} (n^2-(m+1)^2)\cdots (n^2-(k-1)^2), \lb{b} \\
     c_{m,k}
 &:= \begin{cases} 1 & m=0, \\ 
      {k-1+m \choose 2m} (2m-1)!!^2 & m\geq 1
         \end{cases}. \lb{c}
}
In particular,
\eas{
    \frac{1}{n^2}(D_{n,1}(x))^2
& = 
       \frac{1}{1-x^2} , \\
    \frac{1}{n^2}(D_{n,2}(x))^2
& = 
      \frac{(n^2-1)}{(1-x^2)^2} + \frac{1}{(1-x^2)^3},\\
    \frac{1}{n^2} (D_{n,3}(x))^2
& = 
      \frac{(n^2-1)(n^2-4)}{(1-x^2)^3} + \frac{3(n^2-4)}{(1-x^2)^4}
       + \frac{9}{(1-x^2)^5}.
}
Now, using \rf[Df] and \rf[b], we see that
\eas{
       (D_{n,k}(x))^2
&= n^2 \sum_{m=0}^{k-1}
    \frac{c_{m,k}}{(1-x^2)^{k+m}}
       (n^2-(m+1)^2)\cdots(n^2-(k-1)^2) \\
&= \frac{n^2(n^2\!-\!1^2)\cdots(n^2\!-\!(k\!-\!1)^2)} {(1-x^2)^k}
       \left(1 + \sum_{m=1}^{k-1}  \frac{c_{m,k}}{(1-x^2)^m}
         \frac{1}{(n^2\!-\!1^2)\cdots(n^2\!-\!m^2)} \right) \\
&=: (A_{n,k}(x))^2  \left( 1 + B_{n,k}(x) \right).
}
Clearly
$$
     (A_{n,k}(x))^2 \le \frac{n^{2k}}{(1-x^2)^k}.
$$
We now find an upper bound for the sum $B_{n,k}(x)$. Using \rf[c] we expand $c_{m,k}$ as
$$
          c_{m,k} 
=   {k-1+m \choose 2m}(2m-1)!!^2
  = \frac{(2m-1)!!^2}{(2m)!} \frac{(k-1+m)!}{(k-1-m)!}.
$$
In the latter expression, we estimate the first factor using Wallis' inequality (see, e.g., \cite[Chpt.\ 22]{bullen2015dictionary}):
$$
     \frac{(2m-1)!!^2}{(2m)!} = \frac{(2m-1)!!}{(2m)!!} < \frac{1}{\sqrt{\pi m}} \le \frac{1}{\sqrt{\pi}}.
$$
For the second factor, we observe that
$$
     \frac{(k-1+m)!}{(k-1-m)!}
= \frac{k}{k+m} (k^2-1^2)\cdots(k^2-m^2) \\
< (k^2-1^2)\cdots(k^2-m^2).
$$
Therefore
$$
    B_{n,k}(x)
\le \frac{1}{\sqrt{\pi}} \sum_{m=1}^{k-1} \frac{1}{(1-x^2)^m}
     \frac{(k^2-1^2)\cdots(k^2-m^2)}{(n^2\!-\!1^2)\cdots(n^2\!-\!m^2)}
\le \frac{1}{\sqrt{\pi}} \sum_{m=1}^{k-1} \left( \frac{1}{1-x^2} \frac{k^2}{n^2}\right)^m,
$$
where in the second step we used the inequality $\frac{k^2-s^2}{n^2-s^2} < \frac{k^2}{n^2}$. 
Finally, if $k \le n \sqrt{(1-\delta^2)/2}$ and $|x| \le \delta$, then $\frac{1}{1-x^2} \frac{k^2}{n^2} \le \frac{1}{2}$, and 
$$
     B_{n,k}(x)
\le \frac{1}{\sqrt{\pi}} \sum_{m=1}^{k-1} \frac{1}{2^m} < \frac{1}{\sqrt{\pi}}.
$$
Altogether, this gives
$$
D_{n,k}(x) < \frac{c_1 n^k}{(1-x^2)^{k/2}}, \quad |x| \le \delta, \qquad c_1 = (1+ 1/\sqrt{\pi})^{1/2} < 1.251,
$$
as required.
}

\begin{corollary} \lb{1/2}
For any $k,n \in \bbN$ and $\delta \in (0,1)$ such that  $k < n \sqrt{(1-\delta^2)/2}$,  
and for any polynomial $p \in \bbP_n$, we have
\be{\lb{Markov}
     \|p^{(k)}\|_{[-\delta,\delta],\infty} 
\le  \frac{1.251 n^k}{(1-\delta^2)^{k/2}} \|p\|_{[-1,1],\infty}.
}
\end{corollary}


\subsection{A lemma on best approximation}


Next, we require the following lemma on best approximation by polynomials. Although well known, we give a proof for completeness.

\begin{lemma} \lb{T}
Let $f \in C^r([a,b])$. Then 
$$
     E_{r-1}(f) := \inf_{p \in \bbP_{r-1}} \|f - p\|_{[a,b],\infty} \le \frac{2}{r!} \left( \frac{b-a}{4} \right)^r \|f^{(r)}\|_{[a,b],\infty}.
$$
\end{lemma}

\prf{
Let $p_\Delta \in \bbP_{r-1}$ be the polynomial that interpolates $f$ at the $r$ points of the set $\Delta = (x_i)^{r}_{i=1}$.  By the Lagrange interpolation formula, for any $x \in [a,b]$ we have
$$ 
      f(x) - p_\Delta(x) = \frac{1}{r!} \omega_\Delta(x) f^{(r)}(\xi), \qquad   \omega_\Delta(x) := \prod_{i=1}^r (x-x_i),
$$
for some $\xi = \xi_x \in [a,b]$. It follows that
$$
      E_{r-1}(f) 
\le \inf_{\substack{\Delta\subset[a,b] \\ |\Delta| = r}} \|f - p_\Delta\|_{[a,b],\infty}
\le \frac{1}{r!} \inf_{\substack{\Delta\subset[a,b] \\ |\Delta| = r}} \|\omega_\Delta\|_{[a,b],\infty} \cdot \|f^{(r)}\|_{[a,b],\infty}.
$$
It is well-known that the infimum of $\|\omega_\Delta\|_{[a,b],\infty}$ is attained for the set $\Delta_*$ 
of zeros of the Chebyshev polynomial $T_r^*$ on $[a,b]$, and that for this set we have 
$$
     \|\omega_{\Delta_*}\|_{[a,b],\infty} = \frac{1}{2^{r-1}} \left( \frac{b-a}{2} \right)^r.
$$
This completes the proof.
}


\subsection{Bounding the norm on a subinterval}


We are now in a position to state and prove the main result of this section. Note that, for convenience, we formulate this result in terms of intervals $[-\delta,\delta]$ and $[-1,1]$ as opposed to intervals $[-1,1]$ and $[-\gamma,\gamma]$. The analogous result for the latter is obtained by setting $\delta = 1/\gamma$.

\begin{theorem}\label{t:poly-inequality-general}
Given $\epsilon> 0$, let $p \in \bbP_n$ satisfy 
$$
     \|p\|_{[-1,1],\infty} \le 1/\epsilon,
$$
and assume that, for some $\delta \in (0,1)$ and $m \in \bbN$, we also have
$$
      |p(x_i)| \le 1, \qquad x_i =  - \delta +  2\delta i/m,\quad i = 0,\ldots, m.
$$
If 
\be{\lb{m}
       m \ge \max \left\{ 12 c_1 n \log (1/\epsilon)\frac{\delta}{\sqrt{1-\delta^2}}, c_1  \log^2(1/\epsilon) \right\},
}
for some $0 < \epsilon \leq 1/\E$ and $c_1 \ge \max \{1/\lfloor \log(1/\epsilon) \rfloor, 3/n\}$, 
 then
$$
     \|p\|_{[-\delta,\delta]} \le C_0 := \beta^{3/c_1}(1+8\epsilon) +8 \epsilon,
$$
where $\beta > 1$ is the upper constant in the Coppersmith and Rivlin bound \R{coppersmith-rivlin}.
\end{theorem}

This theorem is essentially a more general version of Theorem \ref{t:polynomial-inequality-main} in which the dependence of the bound on $C_0$ in the constant factor $c_1$ in \R{m} is made explicit. We now show how it implies Theorem \ref{t:polynomial-inequality-main}.

\prf{
[Proof of Theorem \ref{t:polynomial-inequality-main}]
We use Theorem \ref{t:poly-inequality-general} with $\delta = 1/\gamma$ and $c_1 = 3$. Recall that $0 < \epsilon \leq 1/\E$ and $n \geq  \sqrt{\gamma^2-1} \log(1/\epsilon) > 0$ by assumption. In particular, $n \geq 1$ since it is an integer. This implies that
\bes{
\max \{1/\lfloor\log(1/\epsilon) \rfloor, 3/n\} \leq \max \{ 1/\lfloor\log(\E) \rfloor, 3 \} = 3.
}
Hence the value $c_1 = 3$ is permitted. Observe now that, since $\delta = 1/\gamma$,
\eas{
\max \left\{ 12 c_1 n\log (1/\epsilon)  \frac{\delta}{\sqrt{1-\delta^2}}, c_1 \log^2(1/\epsilon) \right\}  &\leq 12 c_1 n \log (1/\epsilon) \frac{\delta}{\sqrt{1-\delta^2}} 
\\
& = 36 n \log(1/\epsilon) \frac{1}{\sqrt{\gamma^2-1}}.
}
Therefore, \R{m-n-eps-poly-growth} implies that \R{m} holds. We deduce that
\bes{
C(m,n,\gamma,\epsilon) \leq C_0 = \beta (1+\epsilon) +\epsilon \leq \beta \left ( 1 + 8/\E \right ) +8/\E < 4 \beta + 3,
}
as required.
}

\prf{
[Proof of Theorem \ref{t:poly-inequality-general}]
For a given $p \in \bbP_n$ that satisfies assumptions of the theorem, set 
$$
       k := \lfloor \log (1/\epsilon) \rfloor.
$$
If $k > n \sqrt{(1-\delta^2)/2}$, then the condition \rf[m] implies that
$$
      m >  c_1 n^2 \max  \left \{3\sqrt{2} \delta, (1-\delta^2)/2 \right\} > c_1 n^2 /3.
$$
Hence the polynomial $p$ is bounded on the interval $[-\delta, \delta]$ at $m+1 > c_1 n^2/3 + 1$ equidistant points. We deduce that
\be{\lb{p_i}
      \|p\|_{[-\delta,\delta],\infty} \le \beta^{3/c_1},
}
by the Coppersmith and Rivlin bound \R{coppersmith-rivlin}.

We now suppose that $k \le n \sqrt{(1-\delta^2)/2}$, so that the Markov-type inequality \rf[Markov] is applicable. First, consider partition of the interval $[-\delta,\delta]$ with the $m+1$ equispaced points 
$$
     x_i = - \delta  + 2 \delta i/m, \qquad i = 0,\ldots,m, \qquad 
    m \ge \max \left \{12 c_1 n k \frac{\delta}{\sqrt{1-\delta^2}}, c_1 k^2 \right\}.
$$
Take any subinterval $I = [x_r, x_s]  \subset [-\delta,\delta]$ containing $\lfloor c_1 k^2 \rfloor + 1$ points,  
so that 
\be{\lb{I}
     |I| = \frac{\lfloor c_1 k^2 \rfloor}{m} 2 \delta \le  \frac{c_1 k^2 }{m} 2\delta \le \frac{k \sqrt{1-\delta^2}}{6n},
}
where we used the inequality $m \ge 12 c_1 n k \frac{\delta}{\sqrt{1-\delta^2}}$. 

Now, with the same $k = \lfloor \log(1/\epsilon) \rfloor$, let $Q \in \bbP_k$ be the polynomial of best approximation to $p\in\bbP_n$ from $\bbP_k$ on $I$:
$$
     \|p - Q\|_{I,\infty} = \inf_{q \in \bbP_k} \|p -q\|_{I,\infty} =: E_k(p) \le E_{k-1}(p).
$$
By Lemma \ref{T}, 
$$
    \|p-Q\|_{I,\infty} \le \frac{2}{k!} \left(\frac{|I|}{4}\right)^k \|p^{(k)}\|_{I,\infty},
$$
and by Corollary \ref{1/2}
$$
     \|p^{(k)}\|_{I,\infty} 
\le \|p^{(k)}\|_{[-\delta,\delta],\infty} 
\le  1.251 n^k (1/\sqrt{1-\delta^2})^k  \|p\|_{[-1,1],\infty}.
$$
Hence, using the well-known estimate $k! \ge \sqrt{2\pi}\sqrt{k}  (k/\E)^k$, and the bound \rf[I] for $|I|$, we obtain 
\eas{
         \|p-Q\|_{I,\infty} 
&\le \frac{2.502}{\sqrt{2\pi} \sqrt{k}} \frac{\E^k}{k^k} \left( \frac{k \sqrt{1-\delta^2}}{4 \cdot 6n} \right)^k 
          n^k (1/\sqrt{1-\delta^2})^k \|p\|_{[-1,1],\infty} \\
&\le \rho^k \|p\|_{[-1,1],\infty} , \qquad \rho = \frac{\E}{24} < \frac{1}{\E^2}.
}
Now, recalling that $k = \lfloor \log (1/\epsilon) \rfloor$, whereas $\|p\|_{[-1,1],\infty} \le 1/\epsilon$, we conclude that
\be{\lb{pQ}
       \|p-Q\|_{I,\infty}  \le \E^{-2 \log(1/\epsilon) + 2} 1/\epsilon \le \E^2 \epsilon < 8 \epsilon.
}
 On the other hand, by construction the interval $I = [x_r, x_s]$ 
contains $\lfloor c_1 k^2 \rfloor +1$ equispaced points from the $(x_i)$.  By \R{coppersmith-rivlin},  
for $Q \in \bbP_k$, we have  
\be{\lb{Qi}
     \|Q\|_{I,\infty} \le  C_2 \max_{x_i \in I} |Q(x_i)|, \qquad C_2 = \beta^{c_3}, \qquad 
    c_3 = \frac{k^2}{\lfloor c_1 k^2 \rfloor} \le \frac{k^2}{c_1 k^2/2} \le \frac{2}{c_1}. 
}
Here, for the Copppersmith-Rivlin bound, the assumption $c_1k^2 \ge k$, i.e., boundedness of $Q \in \bbP_k$ 
on at least $k+1$ points, is required, and similarly we required $c_1 n^2/3 \ge n$ in obtaining \rf[p_i].  
This is where the theorem's condition 
$$
      c_1 \ge \max\{1/k,3/n\} = \max\{1/ \lfloor\log(1/\epsilon) \rfloor, 3/n\}
$$
came from. Also, since $c_1 k^2 \ge k \ge 1$, we have the inequality  $\lfloor c_1 k^2 \rfloor \ge c_1 k^2/2$ which
we used in \rf[Qi]

We now conclude the proof. Using \rf[pQ] and \rf[Qi], we obtain
\eas{
         \|p\|_{I,\infty} 
\le \|Q\|_{I,\infty}  + 8 \epsilon
\le  C_2 \max_{x_i \in I} |Q(x_i)| + 8 \epsilon
\le C_2 (\max_{x_i \in I} |p(x_i)| + 8 \epsilon) + 8 \epsilon
&\le C_2 (1+8\epsilon) + 8\epsilon,
}
as required.
}

\section{Proofs of the possibility and impossibility theorems}\label{s:main-thm-proof}

We now prove Theorems \ref{t:possibility-thm}--\ref{t:impossibility-extended}. For these, we first require the following result. 

\lem{
\label{poly-approx-bounds}
Suppose that $E_{\theta} \subset \bbC$ is the Bernstein ellipse \R{Bernstein-ellipse} with parameter $\theta > 1$ and $f \in B(E_{\theta})$. Then there exists a polynomial $p \in \bbP_n$ for which
\bes{
\nmu{f - p}_{[-1,1],\infty} \leq \frac{2}{\theta - 1} \nm{f}_{E_{\theta},\infty} \theta^{-n}.
}
Moreover, this polynomial satisfies
\bes{
\nm{p}_{E_{\tau},\infty} \leq \frac{2}{1-\tau/\theta} \nm{f}_{E_{\theta},\infty},\qquad 1 < \tau < \theta
}
and
\bes{
\nm{p}_{E_{\tau},\infty} \leq \left ( \frac{\tau}{\theta} \right )^n \frac{2\tau/\theta}{\tau/\theta-1}  \nm{f}_{E_{\theta},\infty},\qquad \tau > \theta.
}

}
\prf{
This result is essentially standard. We repeat it to obtain the explicit bounds for $p$.
Since $f \in B(E_{\theta})$ its Chebyshev expansion converges uniformly on $[-1,1]$, i.e.
\bes{
f(x) = \sum^{\infty}_{n=0} c_k \phi_k(x),\quad \phi_k(x) = \cos(n \arccos(x)),
}
and its coefficients satisfy $|c_k| \leq 2 \nm{f}_{E_{\theta}} \theta^{-k}$ \cite[Thm.\ 8.1]{trefethen2013approximation}. Let $p = \sum^{n}_{k=0} c_k \phi_k$. Then
\bes{
\nmu{f-p}_{[-1,1],\infty} \leq 2 \nm{f}_{E_{\theta}} \sum_{k > n}  \theta^{-k} .
}
Evaluating the sum gives the first result. 

For the other results, recall that the Bernstein ellipse is given by $E_{\tau} = \left \{ J(z) : z \in \bbC,\ 1 \leq |z | \leq \tau \right \}$, where $J(z) = \frac12(z+z^{-1})$ is the Joukowsky map, and the Chebyshev polynomials satisfy $\psi_k(J(z)) = \frac12 \left ( z^n + z^{-n} \right )$.
Hence $\nm{\psi_k }_{E_{\tau},\infty} \leq \tau^{n}$. Therefore
\bes{
\nm{p}_{E_{\tau},\infty} \leq 2\nm{f}_{E_{\theta}} \sum^{n}_{k=0} (\tau / \theta)^k,
}
which yields the result.
}

\prf{
[Proof of Theorem \ref{t:possibility-thm}]
We use Theorem \ref{t:acc-cond-poly-frame}. Note that the condition \R{m-n-relation} implies that \R{m-n-eps-poly-growth} holds. Therefore Theorem \ref{t:polynomial-inequality-main} implies that \R{C-cond-main-acc-stab} holds with $c$ as defined therein.

The desired bound for the condition number follows immediately. For the error bound, let $f \in B(E_{\theta})$, where $\theta > \gamma + \sqrt{\gamma^2-1}$, and let $p \in \bbP_n$ be as in the previous lemma. Since $[-\gamma,\gamma] \subset E_{\tau}$ where $\tau = \gamma + \sqrt{\gamma^2-1} < \theta$ by assumption, this lemma gives
\eas{
\nmu{f - p}_{[-1,1],\infty} + (n+1)\epsilon \nm{p}_{[-\gamma,\gamma],\infty} & \leq \left ( \frac{2}{\theta-1} \theta^{-n} + \frac{2(n+1)\epsilon}{1-(\gamma + \sqrt{\gamma^2-1}) / \theta} \right )  \nm{f}_{E_{\theta},\infty}
\\
& \leq g(\theta,\gamma) \left ( \theta^{-n} + n \epsilon  \right )  \nm{f}_{E_{\theta},\infty},
}
for some function $g(\theta,\gamma)$ depending on $\theta$ and $\gamma$ only. Here we also use the fact that $n \geq 1$, which follows from the fact that $n$ is an integer and $n \geq \log(1/\epsilon) > 0$ by assumption. Supposing that \R{C-cond-main-acc-stab} holds, Theorem \ref{t:acc-cond-poly-frame} now gives, up to a constant change in $g(\theta,\gamma)$,
\bes{
\nmu{f - \cP^{\epsilon',\gamma}_{m,n}(f) }_{[-1,1],\infty} \leq c \sqrt{m} g(\theta,\gamma) \left ( \theta^{-n} + n \epsilon \right ) \nm{f}_{E_{\theta},\infty},
}
which is the desired error bound with respect to $n$. To obtain the error bound with respect to $m$, we simply notice that \R{m-n-relation} implies that $m \leq 36 n \log(1/\epsilon) / \sqrt{\gamma^2-1} + 1 = n/c^*+1$. Hence $\theta^{-n} \leq \theta^{c^*(1-m)} = \rho^{1-m}$, as required.
}

\prf{
[Proof of Theorem \ref{t:possibility-slower-exp}]

The overall argument is similar to the previous proof. The only difference is the estimation of the term
\be{
\label{best-approx-term}
\nmu{f - p}_{[-1,1],\infty} + (n+1)\epsilon \nm{p}_{[-\gamma,\gamma],\infty}.
}
Suppose that $f$ is analytic in $E = E_{\theta}$ for some $\theta$ with $1 < \theta < \tau : = \gamma + \sqrt{\gamma^2-1}$.
In other words, the Bernstein ellipse $E_{\theta}$ does not contain the extended interval $[-\gamma,\gamma]$. Let $1 \leq k \leq n$ and $p \in \bbP_k$ be the polynomial guaranteed by Lemma \ref{poly-approx-bounds}. Then
\eas{
\nmu{f - p}_{[-1,1],\infty}  + (n+1)\epsilon \nm{p}_{[-\gamma,\gamma],\infty} & \leq \left ( \frac{2}{\theta-1} \theta^{-k} +  \left ( \frac{\tau}{\theta} \right )^k \frac{2 \tau / \theta}{\tau/\theta-1} (n+1)\epsilon \right )\nm{f}_{E_{\theta},\infty}
\\
& \leq g(\theta,\gamma) \left ( \theta^{-k} + \left ( \frac{\tau}{\theta} \right )^k n \epsilon \right ) \nm{f}_{E_{\theta},\infty}.
}
We now choose
\bes{
k = \min \left \{ n , \left \lfloor \frac{\log(1/\epsilon)}{\log(\tau)} \right \rfloor \right \},
}
so that
\bes{
 \left ( \frac{\tau}{\theta} \right )^k \epsilon \leq \epsilon^{1-\frac{\log(\tau/\theta) }{ \log(\tau) }} = \epsilon^{\frac{\log(\theta)}{\log(\tau)}  },
}
and
\bes{
\theta^{-k} \leq \theta^{-n} + \theta^{1-\frac{\log(1/\epsilon)}{\log(\tau)}} = \theta^{-n} + \theta \epsilon^{\frac{\log(\theta)}{\log(\tau)}}.
}
Since $n \geq 1$, we deduce that
\bes{
\nmu{f - p}_{[-1,1],\infty}  + (n+1)\epsilon \nm{p}_{[-\gamma,\gamma],\infty} \leq 2 g(\theta,\gamma) \left ( \theta^{-n} + n \epsilon^{\frac{\log(\theta)}{\log(\tau)}} \right )  \nm{f}_{E_{\theta},\infty},
}
as required.
}

\prf{
[Proof of Theorem \ref{t:possibility-algebraic}]

As in the previous proof, we only need to obtain the desired estimate for the term \R{best-approx-term}. Let $f \in C^{k}([-1,1])$. Then $f$ has a $C^k$-extension to the interval $[-\gamma,\gamma]$. Specifically, there is a function $\tilde{f} \in C^{k}([-\gamma,\gamma])$ satisfying $\tilde{f}(x) = f(x)$ for all $x \in [-1,1]$ and
\be{
\label{bounded-extension}
\nmu{\tilde{f}}_{C^{k}([-\gamma,\gamma])} \leq c_{k,\gamma} \nm{f}_{C^{k}([-1,1])},
}
for some constant $c_{k,\gamma} \geq 1$ depending on $k$ and $\gamma$ only. Since $\tilde{f} \in C^{k}([-\gamma,\gamma])$ a classical result (see, e.g., \cite[\S 4.6]{cheney1982introduction}) gives that
\bes{
\inf_{p \in \bbP_n} \nmu{\tilde{f} - p}_{[-\gamma,\gamma],\infty} \leq c'_{k,\gamma} n^{-k} \nmu{\tilde{f}^{(k)}}_{[-\gamma,\gamma],\infty},
}
for some $c'_{k,\gamma} > 0$.
Observe that
\bes{
\nmu{f - p}_{[-1,1],\infty} + (n+1)\epsilon \nm{p}_{[-\gamma,\gamma],\infty} \leq 2 \nmu{\tilde{f} - p}_{[-\gamma,\gamma],\infty} + (n+1) \epsilon \nmu{\tilde{f}}_{[-\gamma,\gamma],\infty}.
}
Therefore
\eas{
\inf_{p \in \bbP_n} \left \{ \nmu{f - p}_{[-1,1],\infty} + (n+1)\epsilon \nm{p}_{[-\gamma,\gamma],\infty}  \right \} & \leq 2 c'_{k,\gamma} n^{-k} \nmu{\tilde{f}^{(k)}}_{[-\gamma,\gamma],\infty} + (n+1) \epsilon \nmu{\tilde{f}}_{[-\gamma,\gamma],\infty}
\\
& \leq \left ( 2 c'_{k,\gamma} n^{-k} + (n+1) \epsilon \right ) \nmu{\tilde{f}}_{C^k([-\gamma,\gamma])}
\\
& \leq 2 \max \{c'_{k,\gamma} ,1 \} c_{k,\gamma} \left ( n^{-k} + n \epsilon \right ) \nmu{f}_{C^k([-1,1])},
}
where in the last step we used \R{bounded-extension} and the fact that $n \geq 1$.
This gives the desired result.
}

We conclude this section with the proof of Theorem \ref{t:impossibility-extended}.

\prf{
[Proof of Theorem \ref{t:impossibility-extended}]
The proof is similar to that of the impossibility theorem shown in \cite{platte2011impossibility}.  First observe that $\cR_m(0) = 0$ for any approximation procedure satisfying \R{extended-error-bound}. Now let $k \in \bbN_0$. Then the condition number \R{kappa-def} satisfies
\bes{
\kappa(\cR_m) \geq \lim_{\delta \rightarrow 0^+} \sup_{\substack{q \in \bbP_k \\ 0 < \nm{q}_{m,\infty} \leq \delta}} \frac{\nm{\cR_m(q)}_{[-1,1],\infty}}{\nm{q}_{m,\infty}}.
}
Let $p \in \bbP_k$ with $\nm{p}_{m,\infty} \neq 0$ and set $q = \delta p / \nm{p}_{m,\infty}$ so that $q \in \bbP_k$ with $\nm{q}_{m,\infty} = \delta$. Since polynomials are entire functions, \R{extended-error-bound} gives that
\eas{
\nmu{\cR_m(q)}_{[-1,1],\infty} & \geq \nm{q}_{[-1,1],\infty} - C \left ( \theta^{-c m^{\tau}} +  \epsilon \right ) \nm{q}_{E_{\theta},\infty}
 \geq \nm{q}_{[-1,1],\infty} \left ( 1 - C \left ( \theta^{-c m^{\tau}} +  \epsilon \right ) \theta^k \right ).
}
Here, in the second step, we used the classical inequality $\nm{q}_{E_{\theta},\infty} \leq \theta^k \nm{q}_{[-1,1],\infty}$, $\forall q \in \bbP_k$. Recalling the definition of $q$, we deduce that
\bes{
\kappa(\cR_m) \geq \left ( 1 - C \left ( \theta^{-c m^{\tau}} + \epsilon \right ) \theta^k \right ) \sup_{\substack{p \in \bbP_k \\ \nm{p}_{m,\infty} \neq 0 }}  \frac{\nm{p}_{[-1,1],\infty}}{\nm{p}_{m,\infty}}.
}
Now let
\bes{
k = \left \lfloor \min \left \{ \frac{\log(\frac{1}{4C \epsilon})}{\log(\theta)} , c m^{\tau} + \frac{\log(\frac{1}{4C})}{\log(\theta)} \right \} \right \rfloor.
}
This choice of $k$ gives $C \left ( \theta^{-c m^{\tau}} + \epsilon \right ) \theta^k \leq 1/2$,
and therefore
\bes{
\kappa(\cR_m) \geq \frac12 \alpha^{k^2/m},
}
where $\alpha > 1$ is as in \R{coppersmith-rivlin}. Observe that this holds for all $1 < \theta \leq \theta^*$. Now choose $\theta = \theta_m = \epsilon^{-\frac{1}{c m^{\tau}}}$,
and observe that $\theta_m \leq \theta^*$ for all large $m$. This value of $\theta$ and the fact that $0 < \epsilon \leq (4 C)^{-2}$ give
\bes{
k = \left \lfloor c m^{\tau} \left ( 1 + \frac{\log(\frac{1}{4C})}{\log(1/\epsilon)} \right )\right \rfloor \geq  \left \lfloor \frac12 c m^{\tau} \right \rfloor ,
}
and therefore $k \geq \frac13 c m^{\tau}$ for all sufficiently large $m$. We deduce that
\bes{
\kappa(\cR_m) \geq \frac12 a^{c^2 m^{2 \tau - 1}/9} \geq \sigma^{m^{2 \tau-1}}, 
}
for some $\sigma > 1$ and all large $m$, as required.
}

\section{Numerical examples}\label{s:numerical}

We conclude this paper with a series of experiments to examine the various theoretical results. Unless otherwise stated, we compute the discrete $L^{\infty}$-norm error of the approximation on a grid of 50,000 equispaced points in $[-1,1]$. Also, we consider the polynomial frame approximation threshold parameter $\epsilon$ rather than $\epsilon'$ (as used in the main theorems). Theoretically, this choice leads to a log-linear sample complexity, but in practice it appears to be adequate.

\begin{figure}[t]
\begin{small}
\begin{center}
\begin{tabular}{ccc}
\includegraphics[width = 0.3\textwidth]{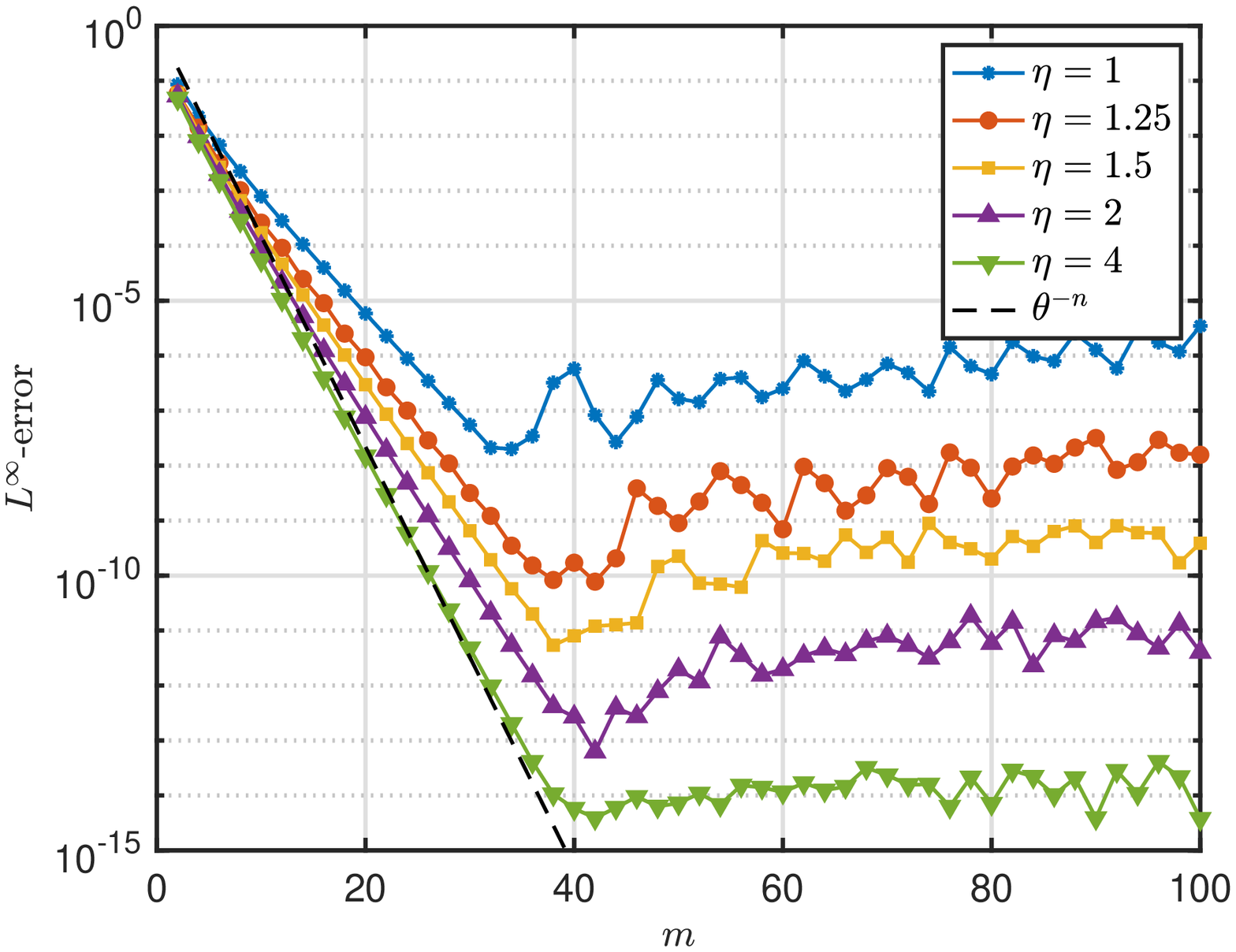}&
\includegraphics[width = 0.3\textwidth]{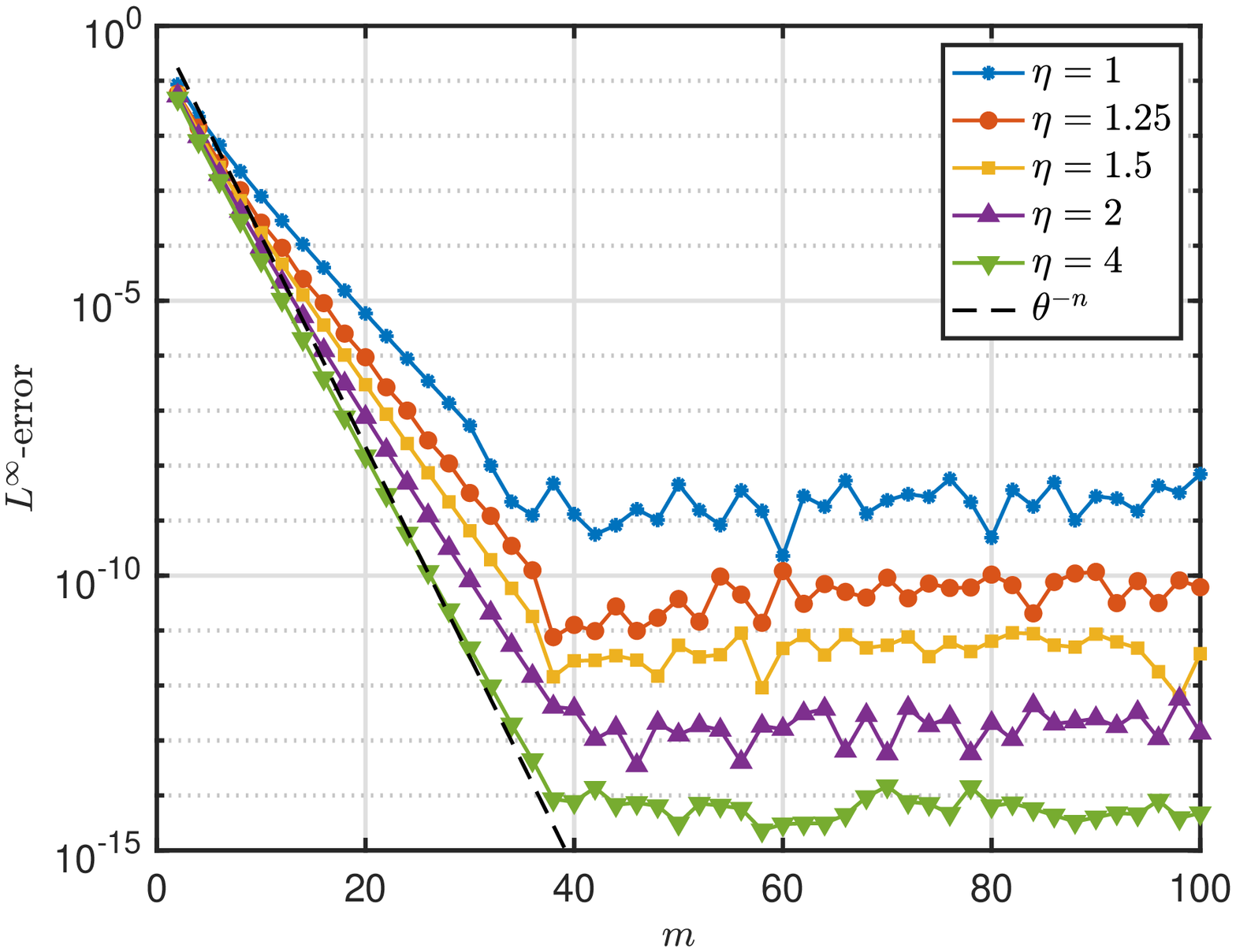}&
\includegraphics[width = 0.3\textwidth]{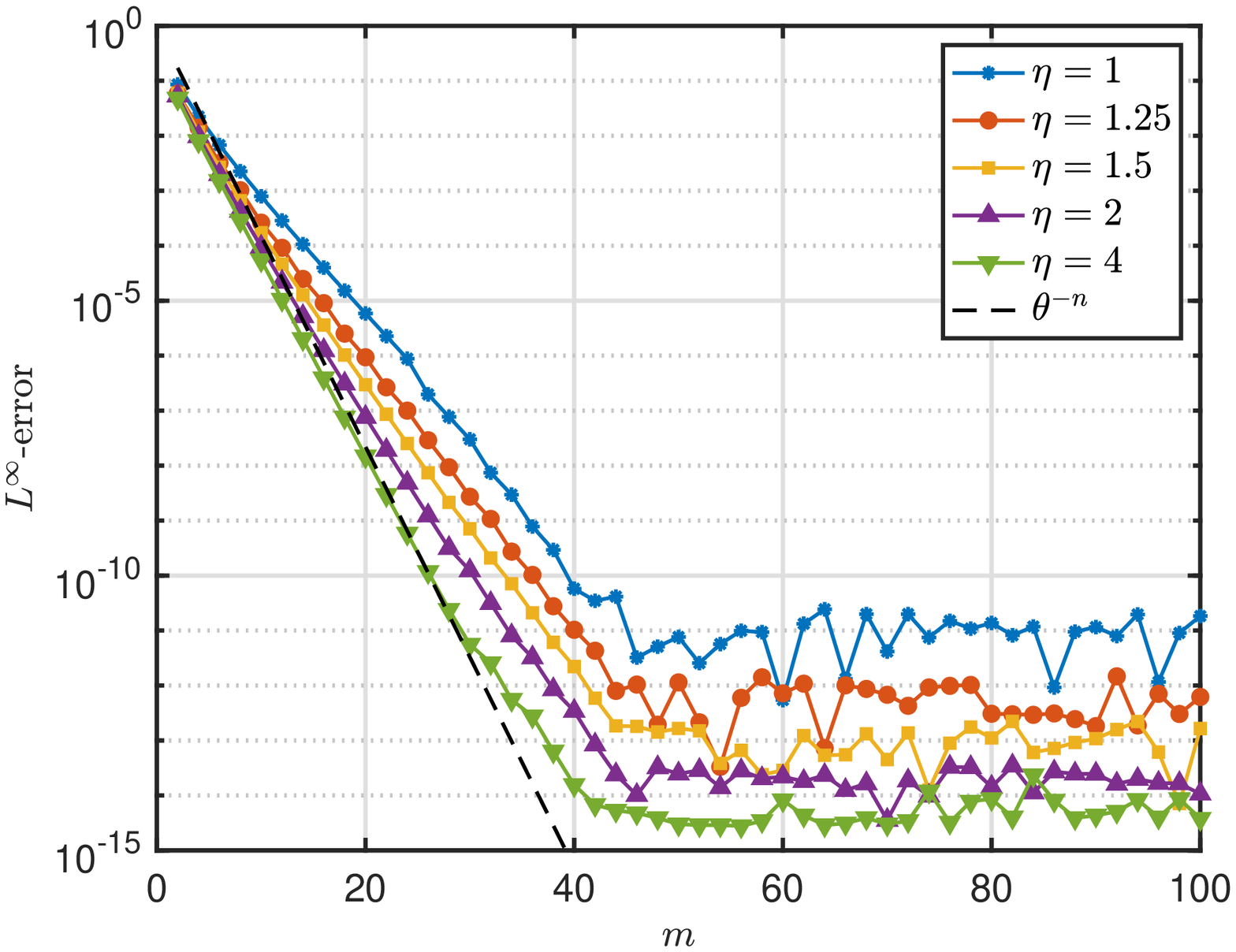}
\\
\includegraphics[width = 0.3\textwidth]{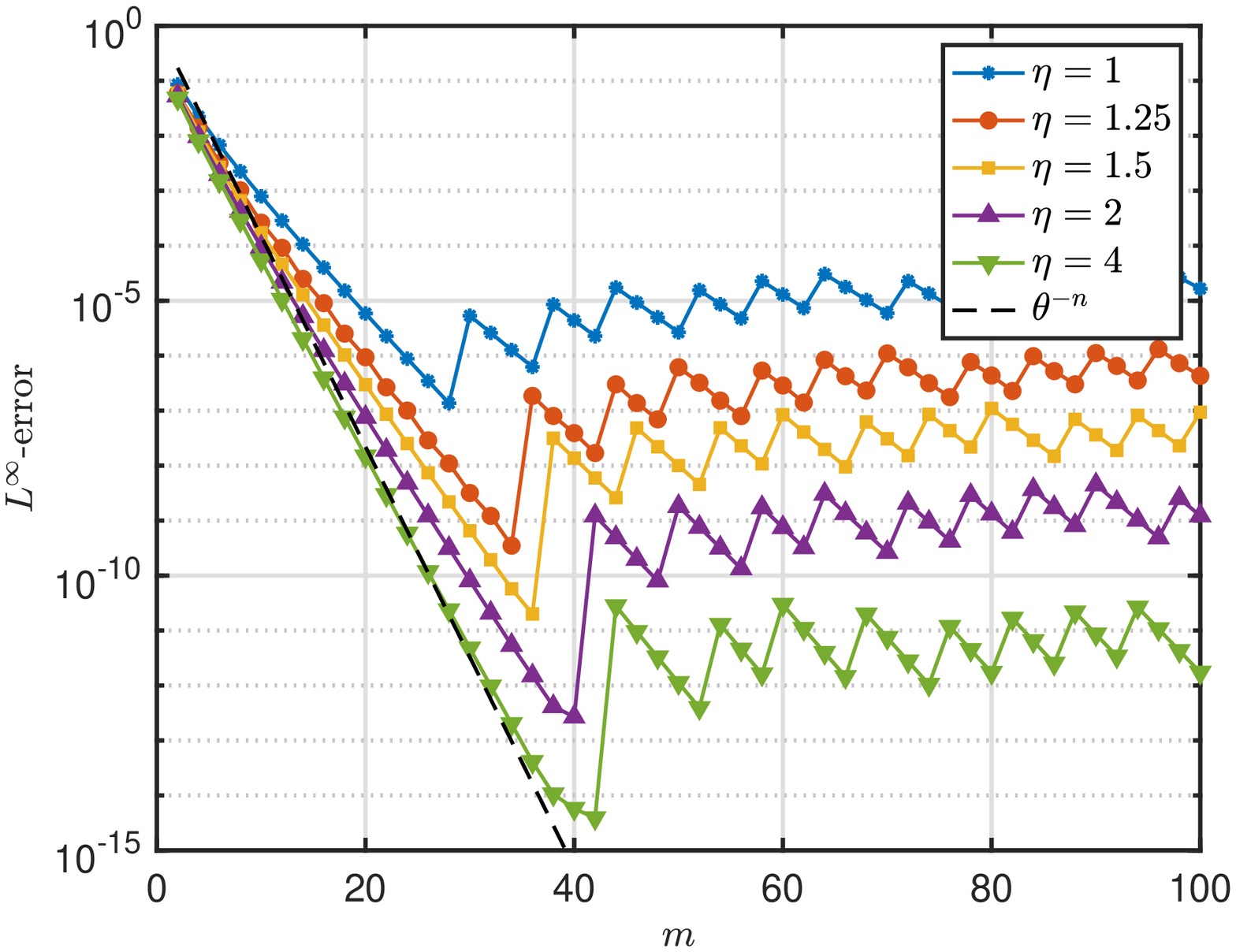}&
\includegraphics[width = 0.3\textwidth]{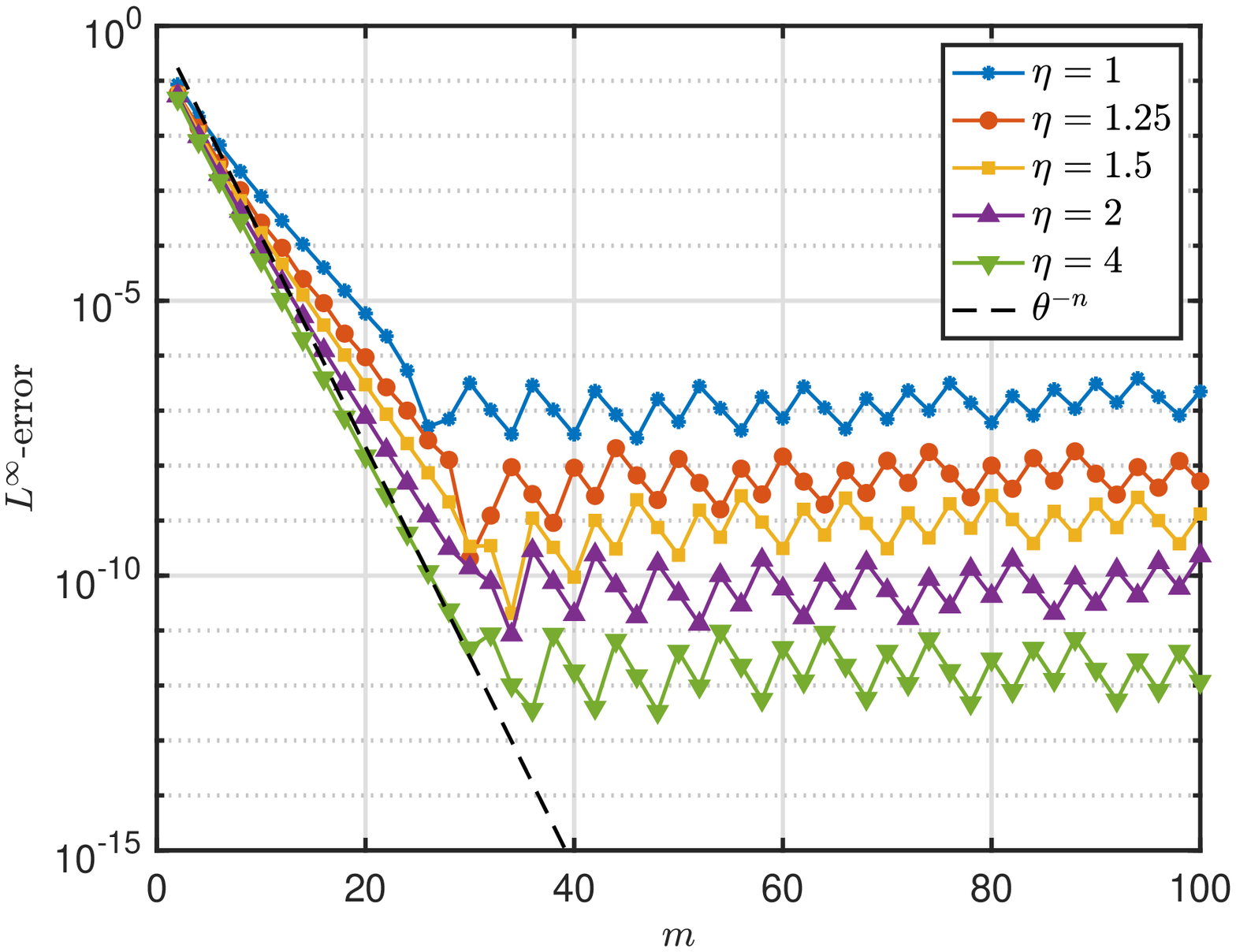}&
\includegraphics[width = 0.3\textwidth]{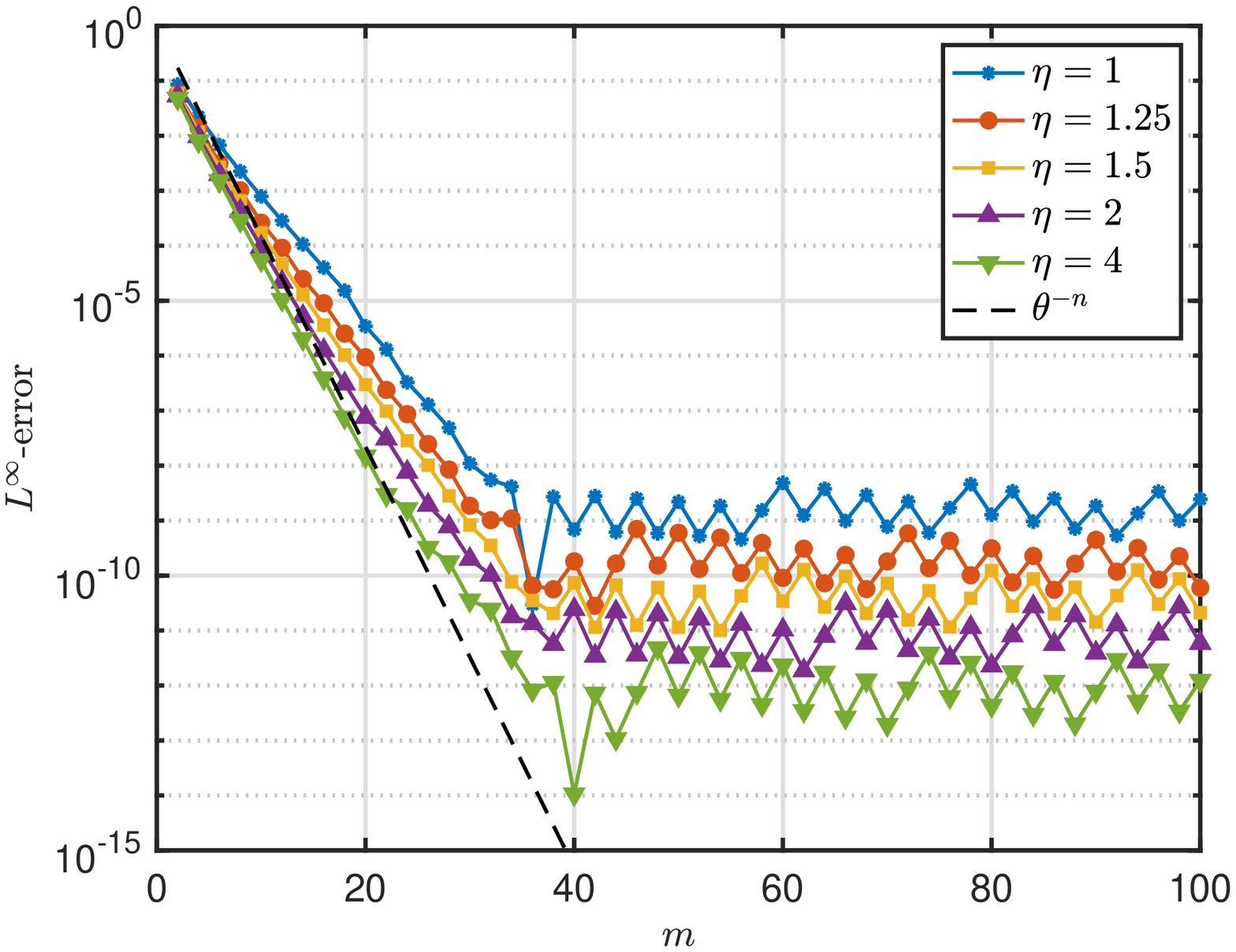}
\\
\includegraphics[width = 0.3\textwidth]{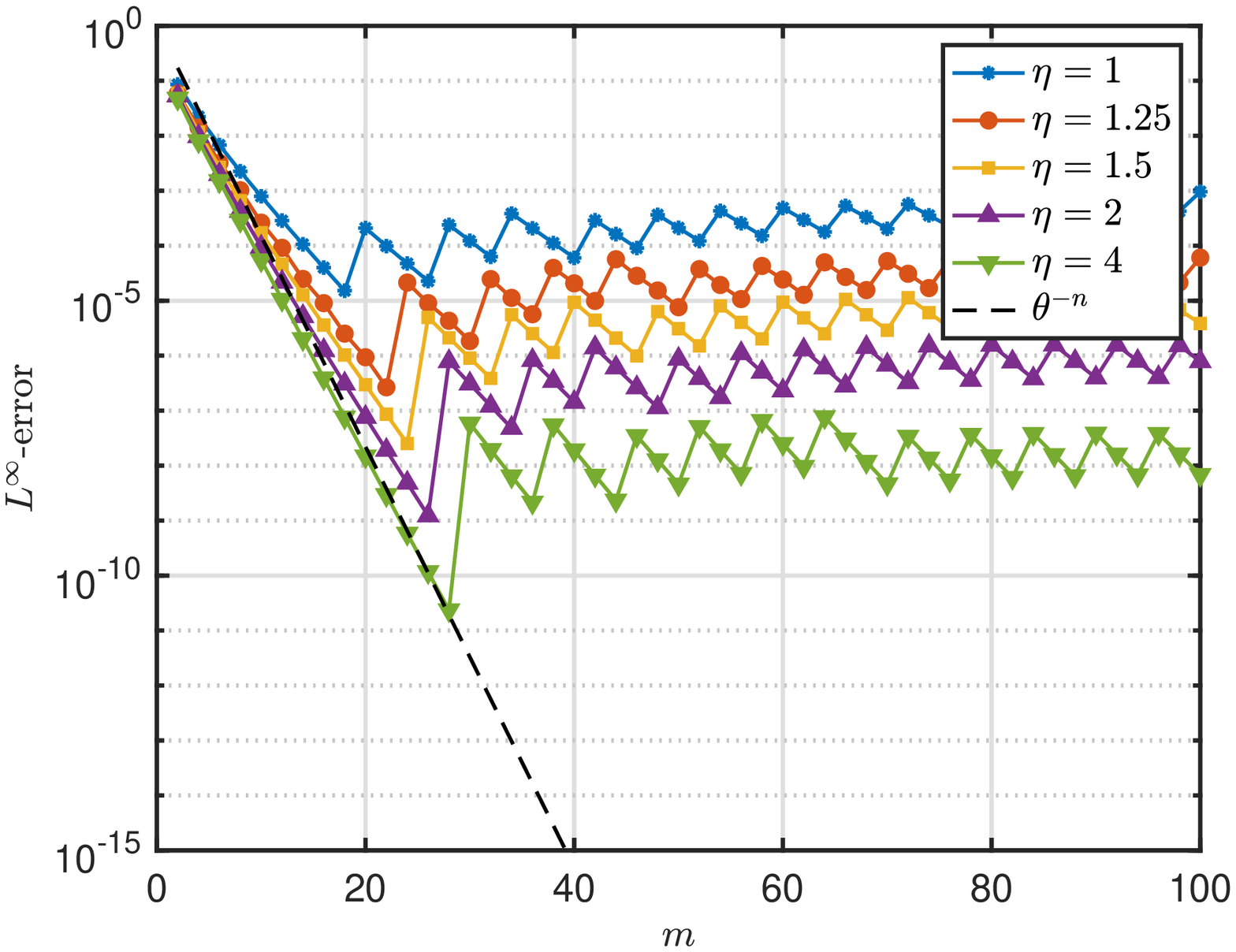}&
\includegraphics[width = 0.3\textwidth]{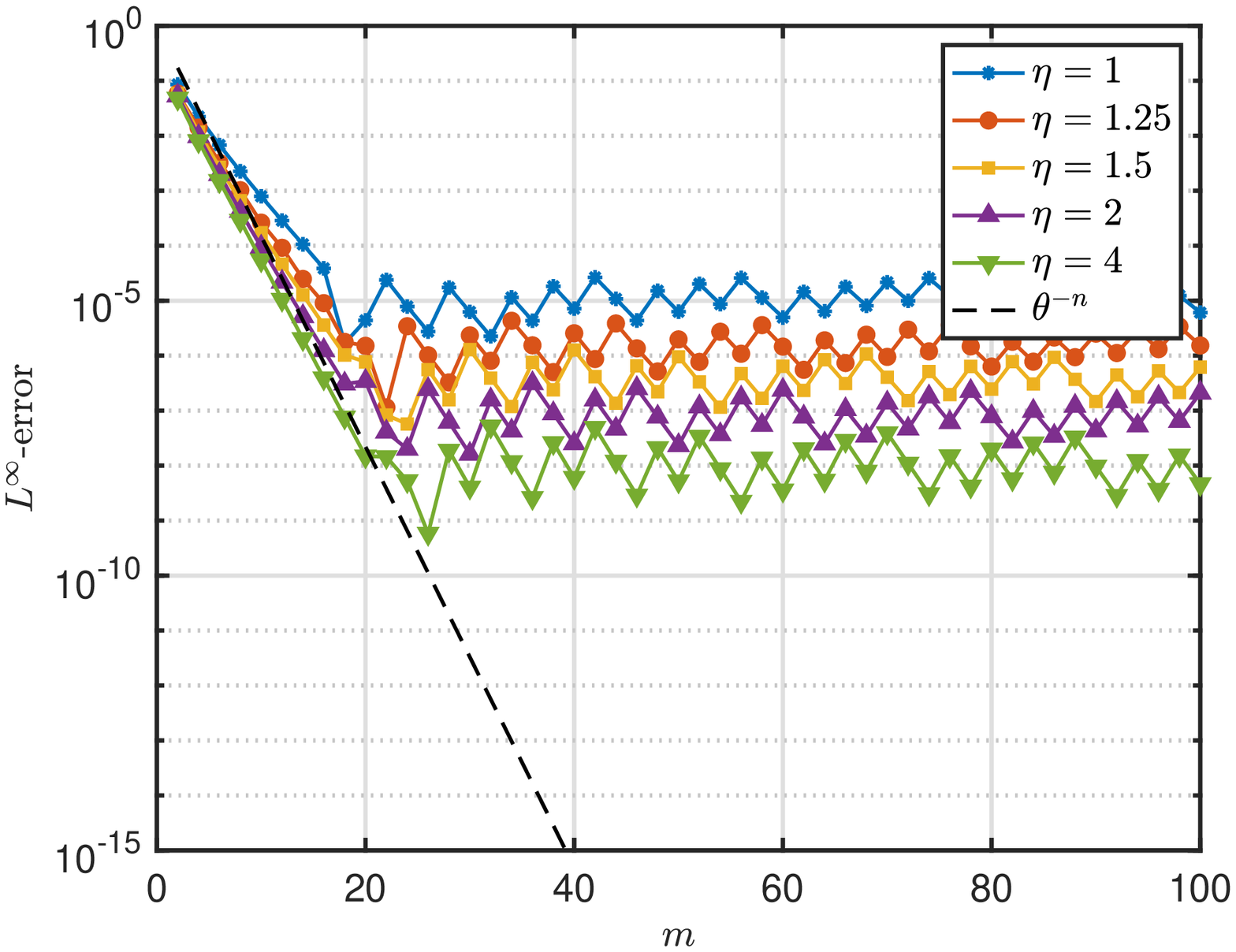}&
\includegraphics[width = 0.3\textwidth]{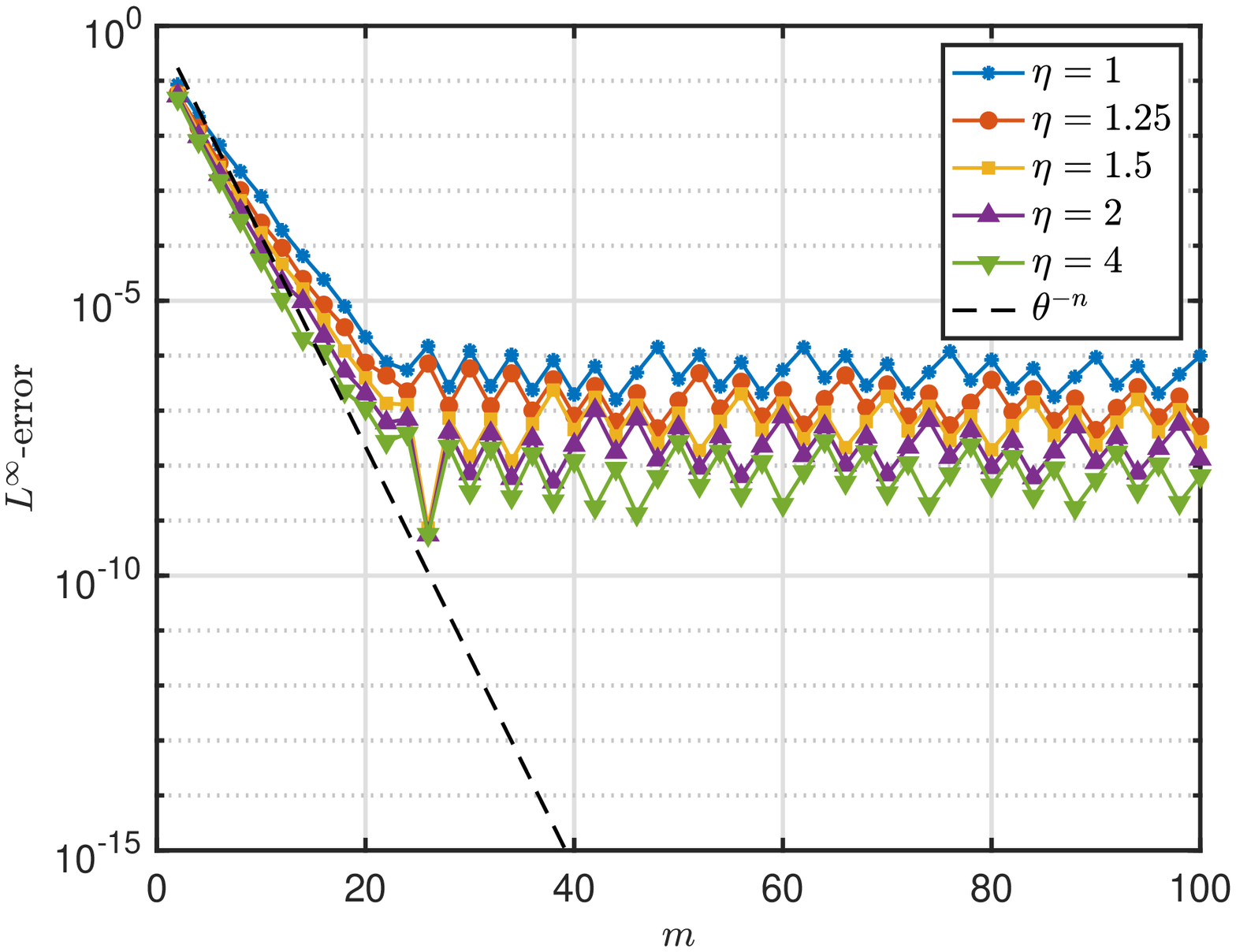}
\\
$\gamma = 1.2$ & $\gamma = 1.4$ & $\gamma = 1.8$
\end{tabular}
\end{center}
\end{small}
\caption{Approximation error versus $n$ for approximating the function $f(x) = \frac{1}{1+x^2}$ via $\cP^{\epsilon,\gamma}_{m,n}$, where $m/n = \eta$, using various different values of $\eta$, $\gamma$ and $\epsilon$. The values of $\epsilon$ used are  $\epsilon = 10^{-14}$ (top), $\epsilon = 10^{-10}$ (middle) and $\epsilon = 10^{-6}$ (bottom). The dashed line shows the quantity $\theta^{-n}$, where $\theta = \sqrt{2}+1$.
}
\label{f:fig1}
\end{figure}

In Fig.\ \ref{f:fig1} we plot the error versus $n$ for different values of the \textit{oversampling} parameter $\eta = m/n$. We compare several different values for the extended domain parameter $\gamma$, and several different values of $\epsilon$. Notice that in all cases, we witness exponential decrease of the error down to some fixed limiting accuracy. The limiting accuracy is related to the stability of the approximation and the parameter $\epsilon$. Observe that it gets smaller with increasing $\eta$, and for sufficiently large $\eta$ it closely tracks the value of $\epsilon$ used. Moreover, it is larger when $\gamma$ is smaller and smaller when $\gamma$ is larger. Both observations are intuitively true. Increasing the number of sample points (i.e.\ larger $\eta$) reduces the maximal growth of a polynomial on $[-1,1]$ relative to its values of the equispaced grid. Similarly, increasing $\gamma$ lengthens the region within which the polynomial cannot exceed the value $1/\epsilon$, and therefore it also cannot grow as large on $[-1,1]$. We also notice that increasing the oversampling parameter makes less difference to the limiting accuracy when $\epsilon$ is larger than it does when $\epsilon$ is smaller. Again, this is intuitively true, since larger $\epsilon$ means the polynomial cannot grow as large on the extended interval. These three observations are also supported by Theorem \ref{t:polynomial-inequality-main}. Here, the sufficient scaling between $m$ and $n$ depends on $\log(1/\epsilon)/\sqrt{\gamma^2-1}$, i.e., it is a decreasing function of both $\epsilon$ and $\gamma$.

In Fig.\ \ref{f:fig1} the function we consider has poles at $x = \pm \I$, meaning that it is analytic within any Bernstein ellipse $E_{\theta}$ for which $\frac12(\theta - \theta^{-1}) = 1$, i.e.\ $1 < \theta < \sqrt{2}+1 \approx 2.41$. Recall that the interval $[-\gamma,\gamma]$ is contained in the Bernstein ellipse $E_{\tau}$ with $\tau = \gamma + \sqrt{\gamma^2-1}$. In particular, $\tau < \sqrt{2}+1$ for $\gamma = 1.2$ and $\gamma = 1.4$. Our analysis in Theorem \ref{t:possibility-thm} therefore asserts exponential decay of the error with rate roughly $(\sqrt{2}+1)^{-n}$ for these two choices of $\gamma$. This is what we observe in practice in Fig.\ \ref{f:fig1}.

On the other hand, when $\gamma = 1.8$, $\tau = \gamma + \sqrt{\gamma^2-1} \approx 3.30 > 2.41 = \sqrt{2}+1$. In this case, our analysis in Theorem \ref{t:possibility-slower-exp} predicts exponential convergence with rate roughly $(\sqrt{2}+1)^{-n}$ down to roughly $\epsilon^{\frac{\log(\sqrt{2}+1)}{\log(\tau)}} \approx \epsilon^{0.74}$, and slower convergence below this level. This is again in agreement with the right column of Fig.\ \ref{f:fig1}.

\begin{figure}[t]
\begin{small}
\begin{center}
\begin{tabular}{ccc}
\includegraphics[width = 0.3\textwidth]{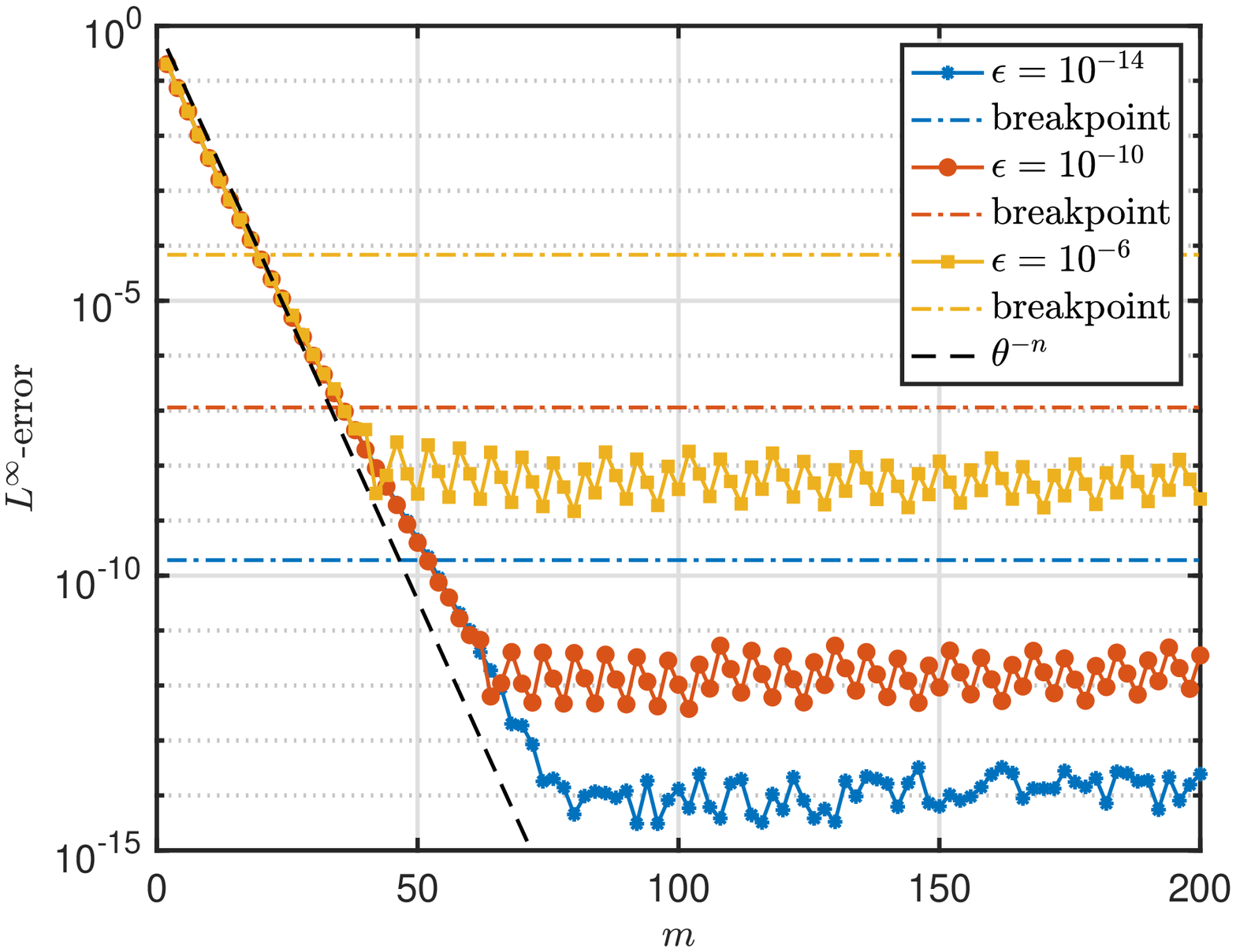}
&
\includegraphics[width = 0.3\textwidth]{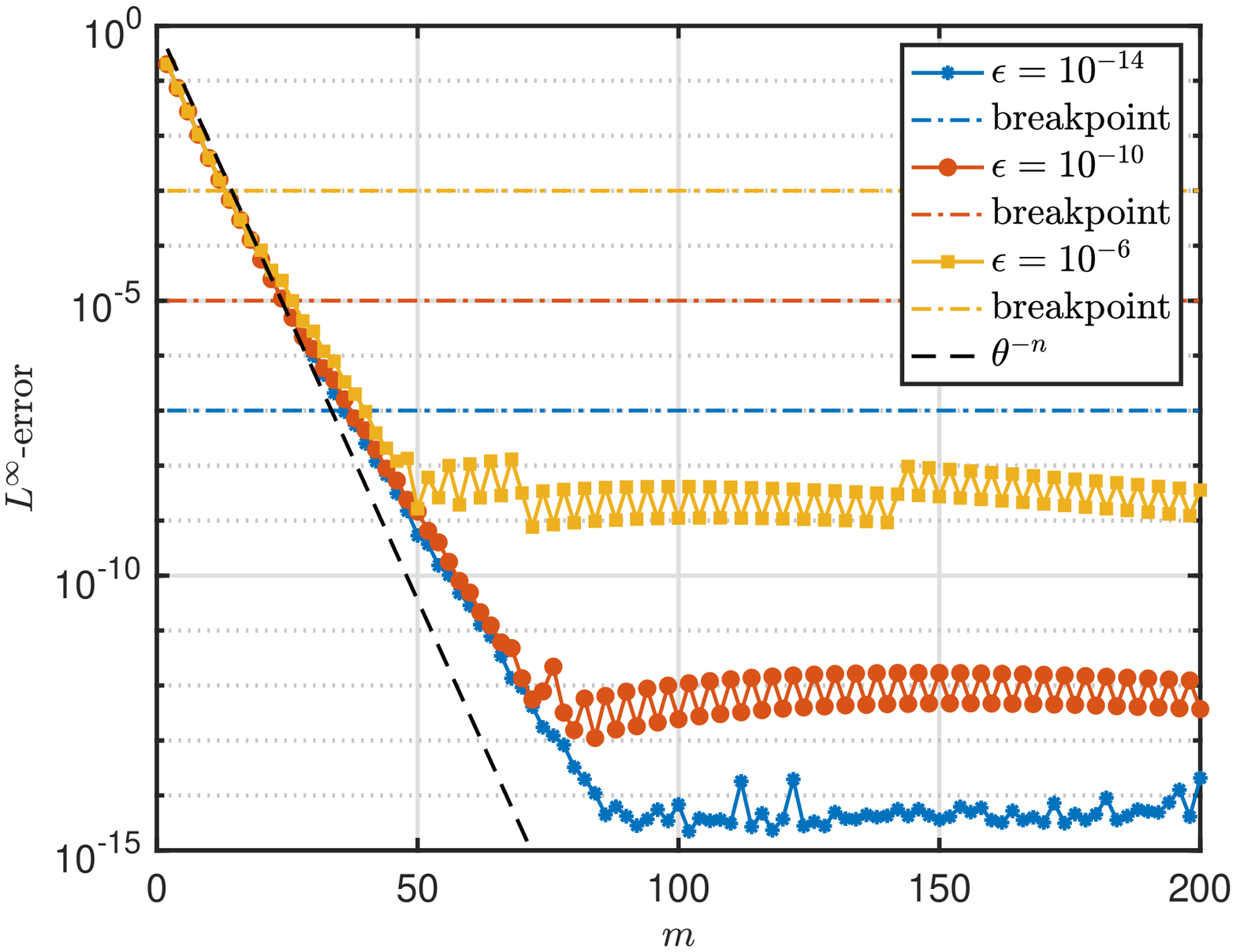}
&
\includegraphics[width = 0.3\textwidth]{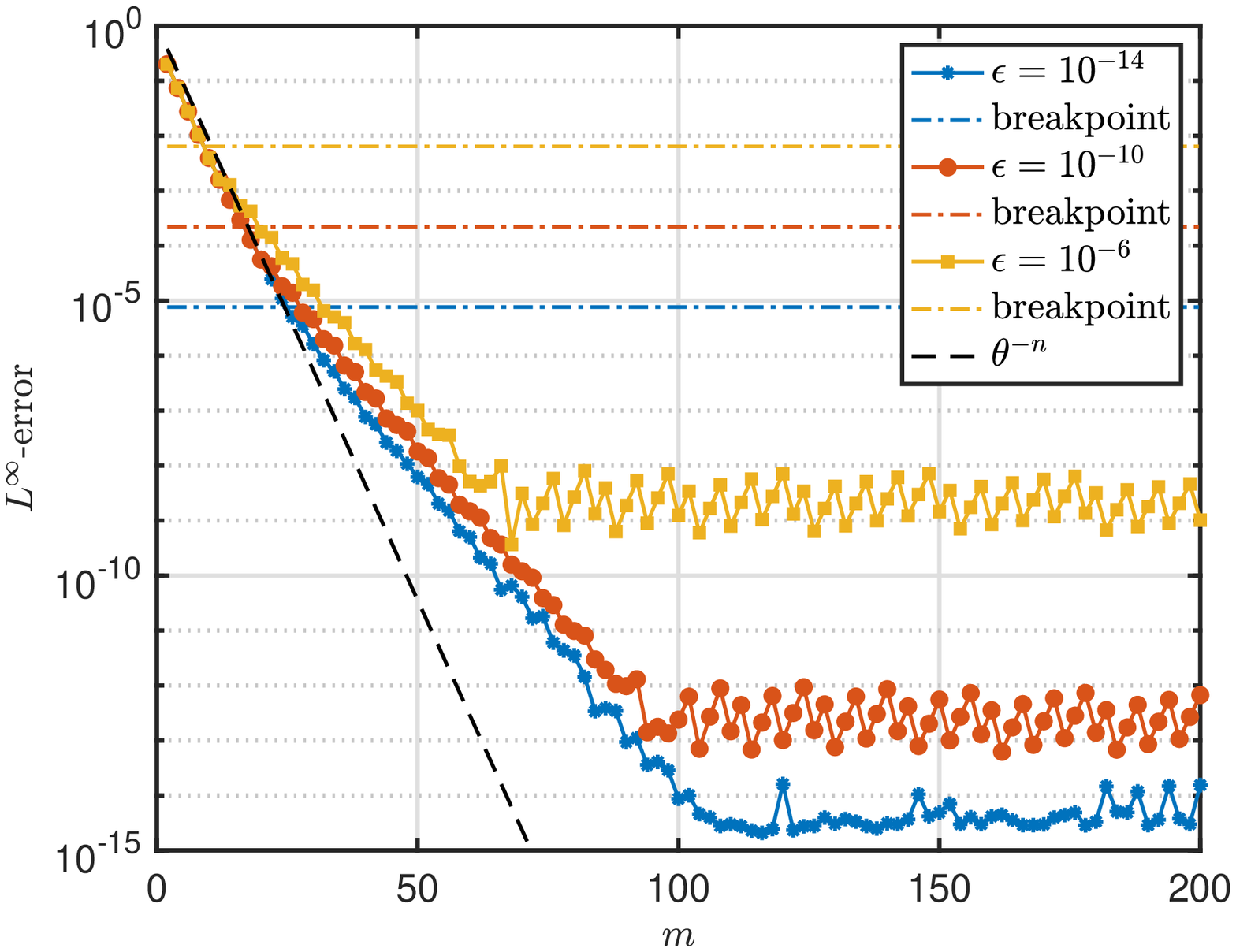}
\\
\includegraphics[width = 0.3\textwidth]{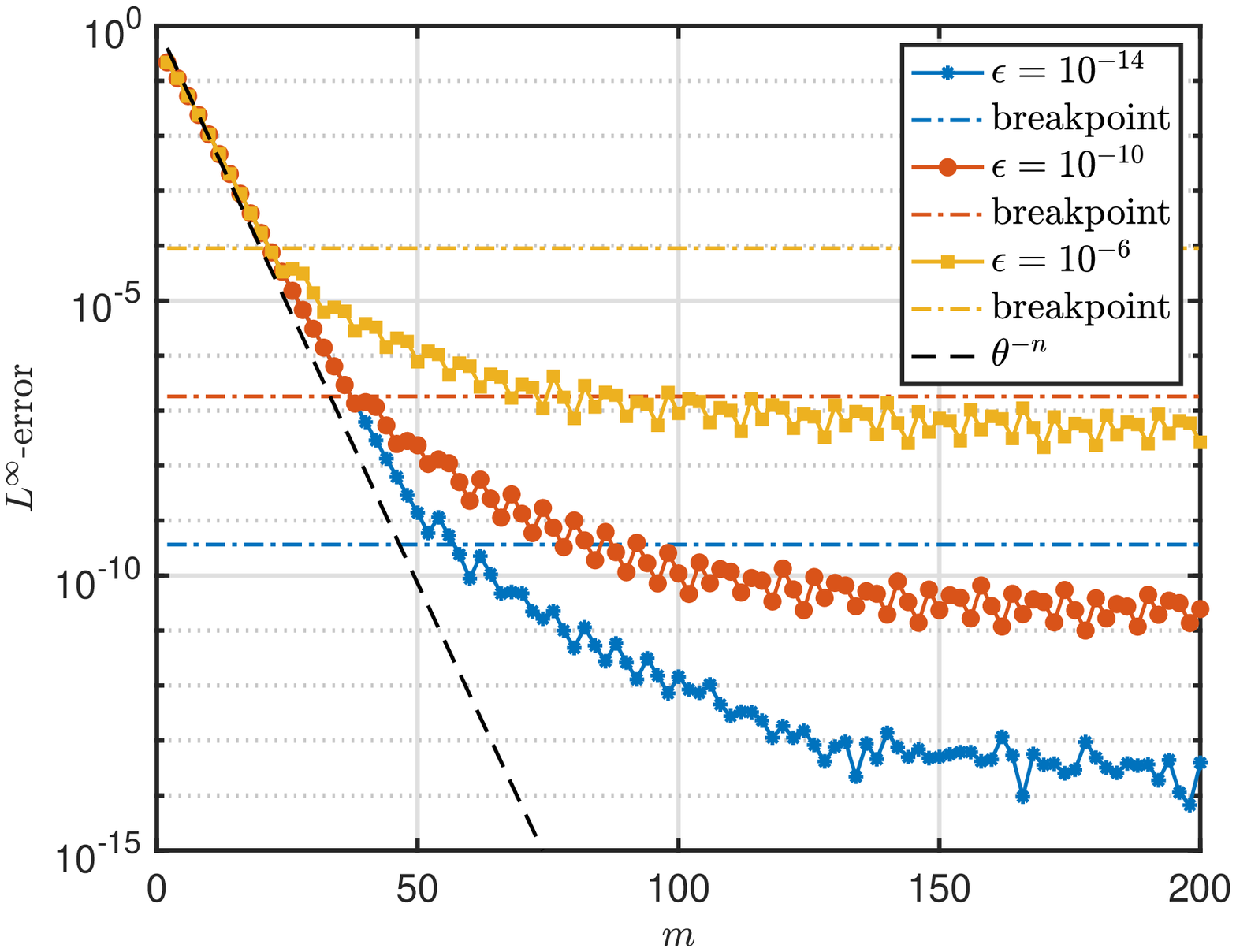}
&
\includegraphics[width = 0.3\textwidth]{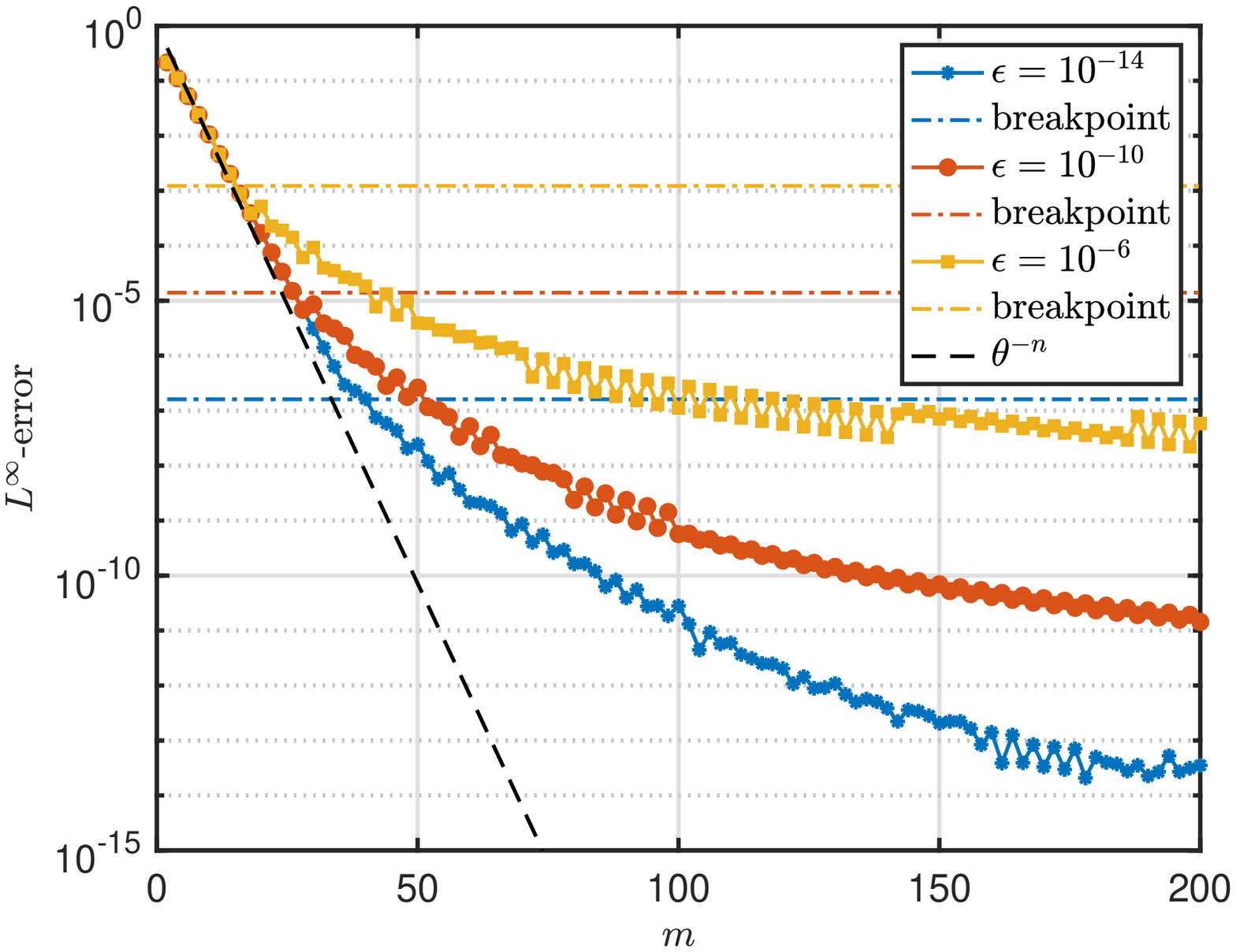}
&
\includegraphics[width = 0.3\textwidth]{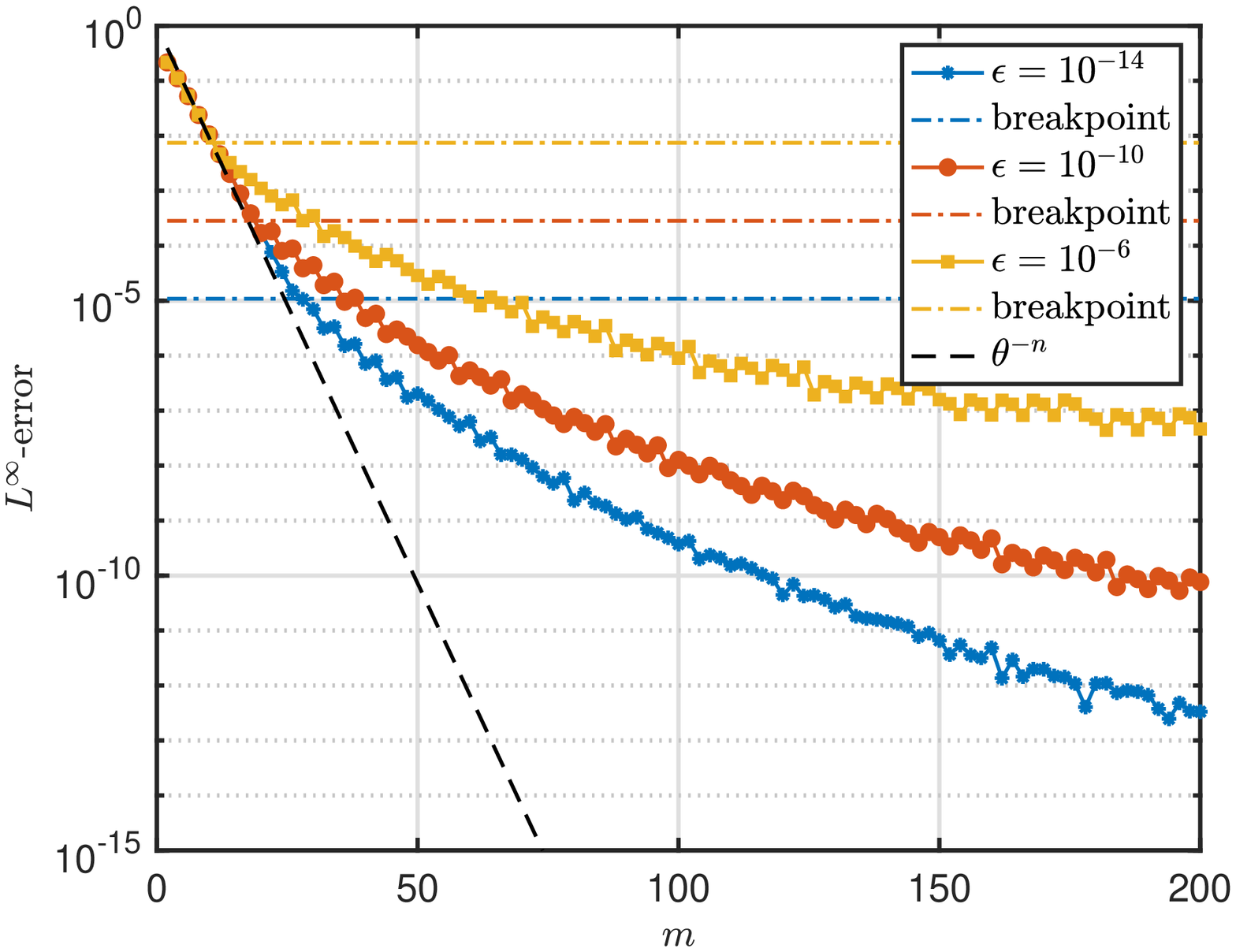}
\\
\includegraphics[width = 0.3\textwidth]{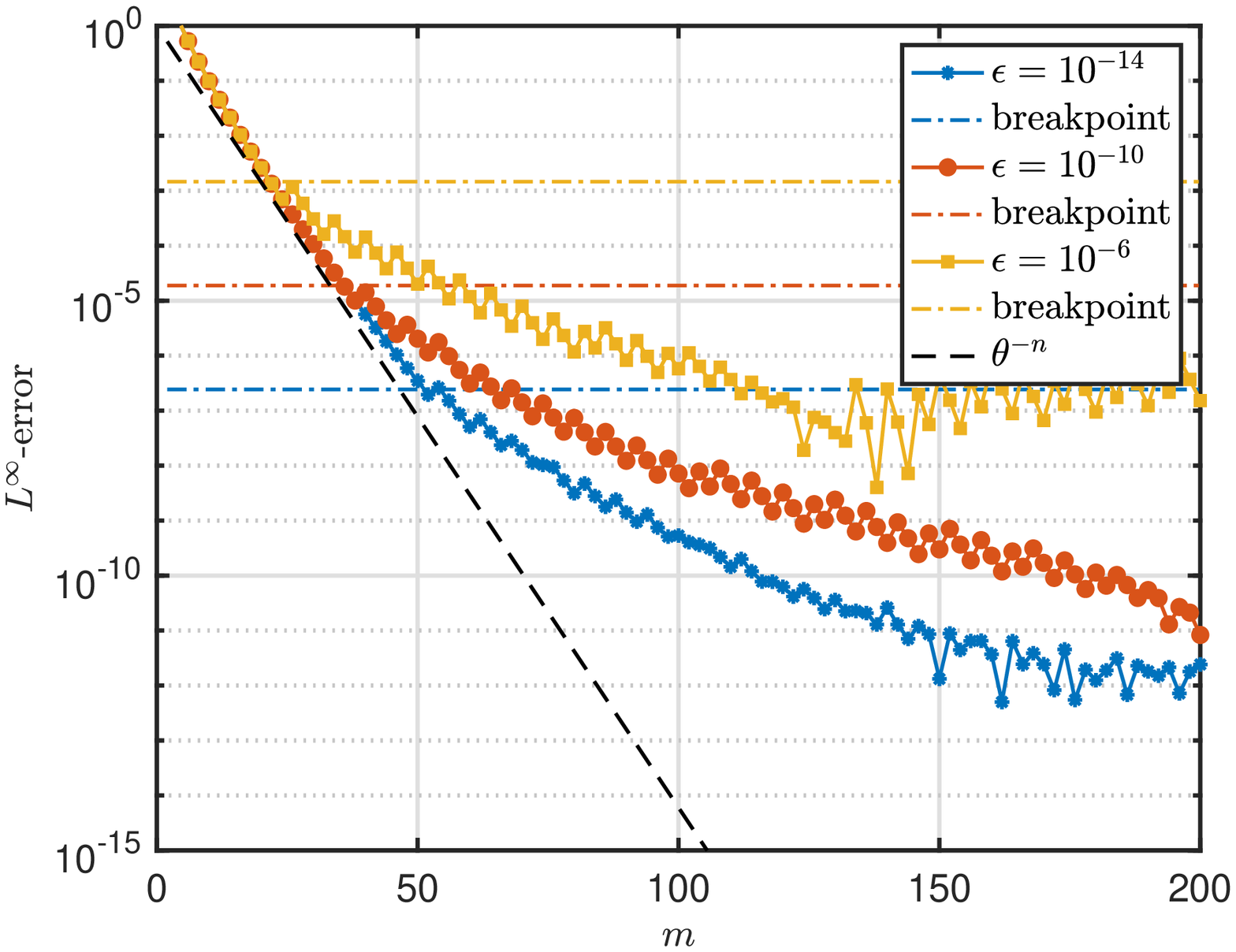}
&
\includegraphics[width = 0.3\textwidth]{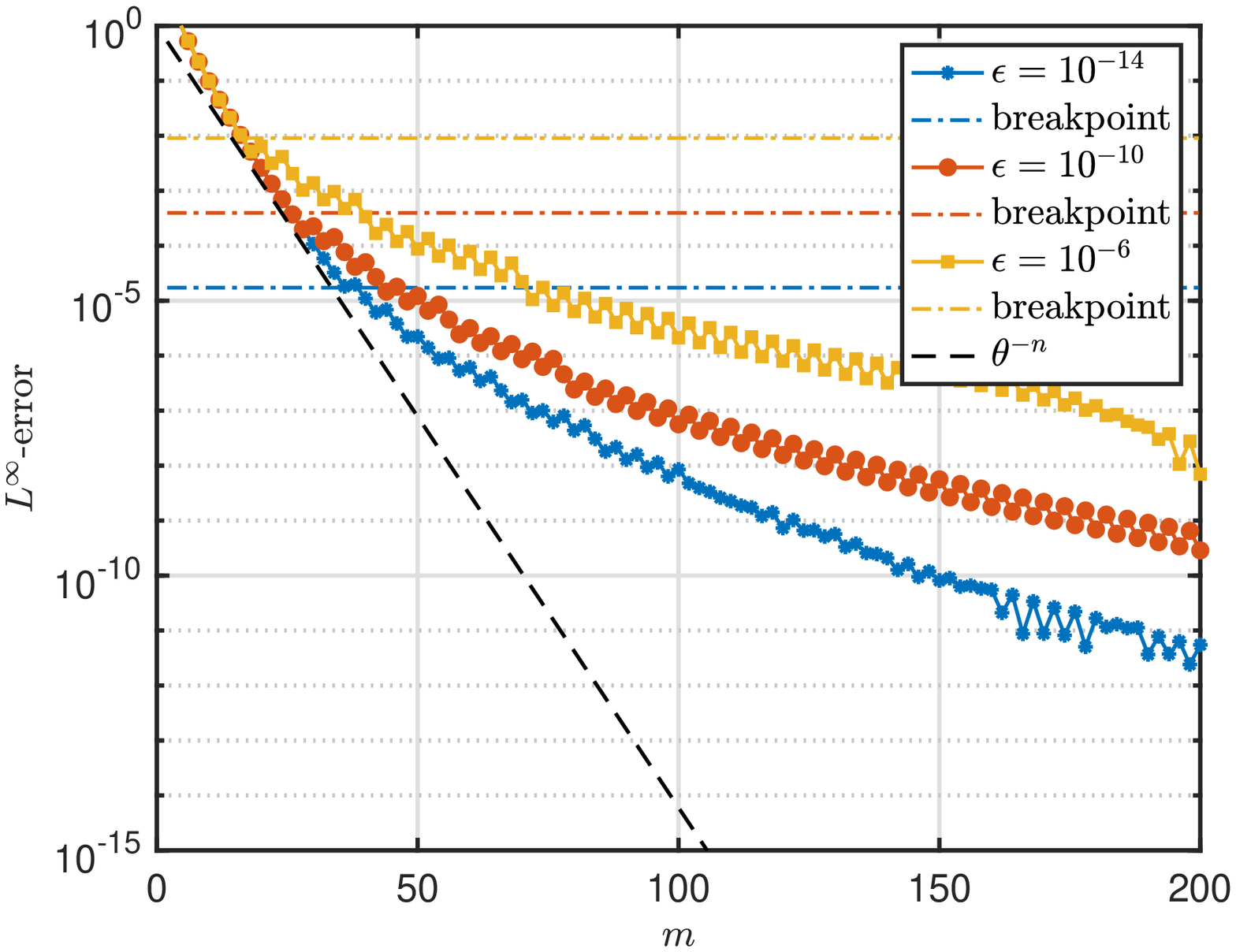}
&
\includegraphics[width = 0.3\textwidth]{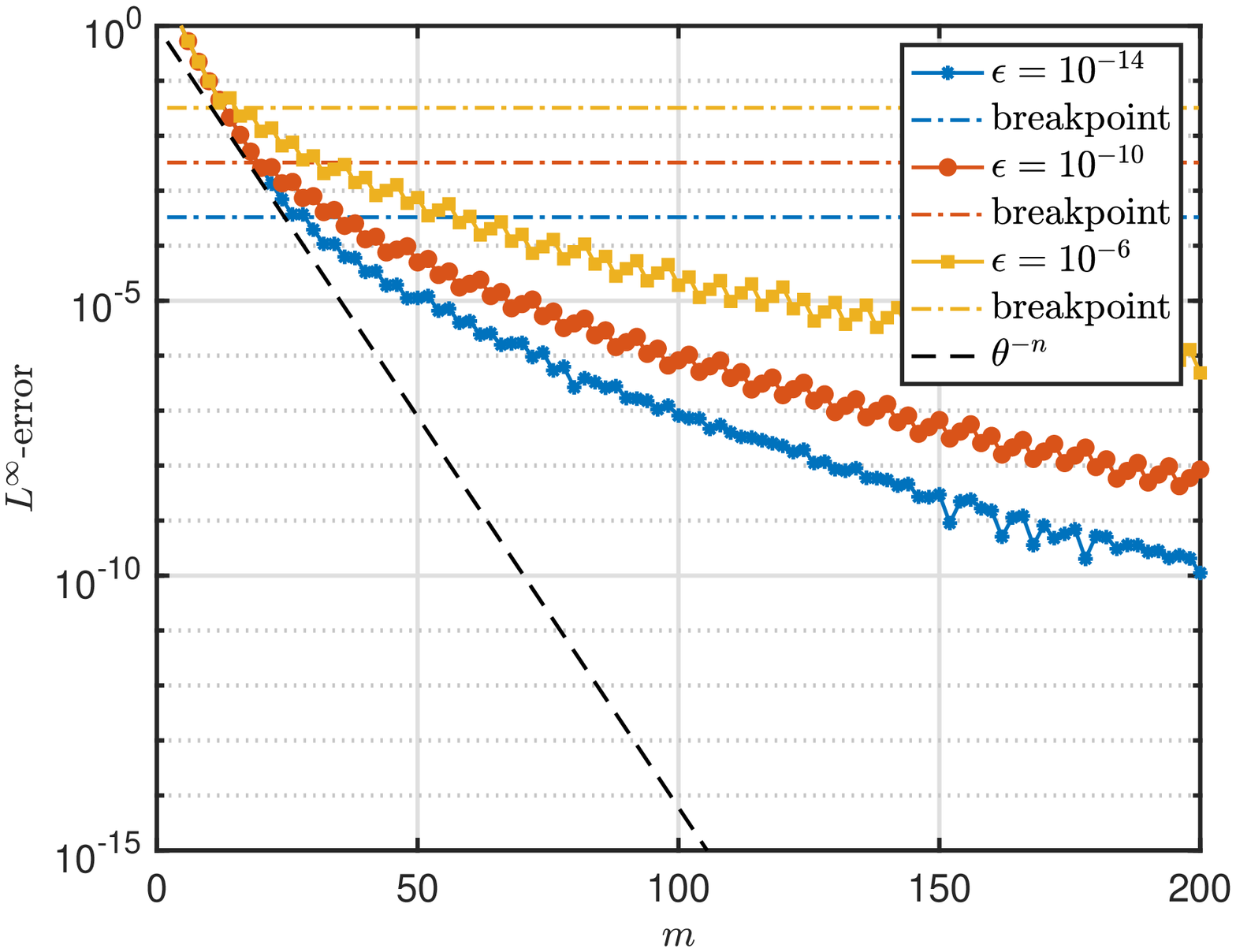}
\\
$\gamma = 1.25$ & $\gamma = 1.5$ & $\gamma = 2$
\end{tabular}
\end{center}
\end{small}
\caption{Approximation errors versus $n$ for approximating the functions $f_1(x) = \frac{1}{1+4x^2}$ (top), $f_2(x) = \frac{1}{10-9 x}$ (middle) and $f_3(x) = 25 \sqrt{9 x^2-10}$ (bottom) via $\cP^{\epsilon,\gamma}_{m,n}$, where $m/n = 4$, using various different values of $\gamma$ and $\epsilon$. The dot-dashed lines show the breakpoints \R{breakpoint} in each case and dashed line shows the quantity $\theta^{-n}$. In this experiment, the values of $\theta$ are $\theta = \frac12(1+\sqrt{5})$ (top), $\theta = \frac{1}{9}(10 + \sqrt{19})$ (middle) and $\theta = \sqrt{10/9} + 1/3$ (bottom).}
\label{f:fig2}
\end{figure}

To further investigate this effect of decreasing error decay for less regular functions, in Fig.\ \ref{f:fig2} we plot the error versus $n$ for several different functions and values of $\gamma$. We also plot the theoretical breakpoint described in Theorem \ref{t:possibility-slower-exp}, i.e.\ the value
\be{
\label{breakpoint}
\epsilon^{\frac{\log(\theta)}{\log(\tau)}},\qquad \tau = \gamma + \sqrt{\gamma^2-1}.
}
In all cases, we see reasonable agreement between the theoretical results and the empirical performance. First, the error decays with rate $\theta^{-n}$ down to the breakpoint, as in Theorem \ref{t:possibility-slower-exp}, before decreasing more slowly beyond it. This second phase is described by Theorem \ref{t:possibility-algebraic} (since analytic functions are infinitely differentiable).

It is notable that the error decay rate after the breakpoint is quite different for the functions considered. This effect has also been observed and discussed in the case of Fourier extensions \cite{adcock2014parameter,adcock2014resolution}. It can be understood through Theorem \ref{t:possibility-algebraic}. Recall that this theorem asserts that the error is bounded by
\be{
\label{limiting-accuracy-analysis}
c g(k,\gamma) \sqrt{m} \left ( n^{-k} + n \epsilon \right ) \nm{f}_{C^k([-1,1])},
}
for any $k \in \bbN$ (since all functions considered are infinitely smooth). The derivatives of the first function $f_1(x) = \frac{1}{1+4x^2}$ do not grow too large on $[-1,1]$ with increasing $k$. Hence the constants $ \nm{f}_{C^k([-1,1])}$ in the error term remain reasonably small and we see rapid decrease in $n$. On the other hand, the derivatves of the functions $f_2$ and $f_3$ grow rapidly with $k$, meaning the constants $ \nm{f}_{C^k([-1,1])}$ also grow rapidly with $k$. Thus, \R{limiting-accuracy-analysis} suggests that the error decays progressively more slowly the closer it gets to the limiting value $\epsilon$. This is exactly the effect we observe in practice.

\begin{figure}[t]
\begin{small}
\begin{center}
\begin{tabular}{ccc}
\includegraphics[width = 0.3\textwidth]{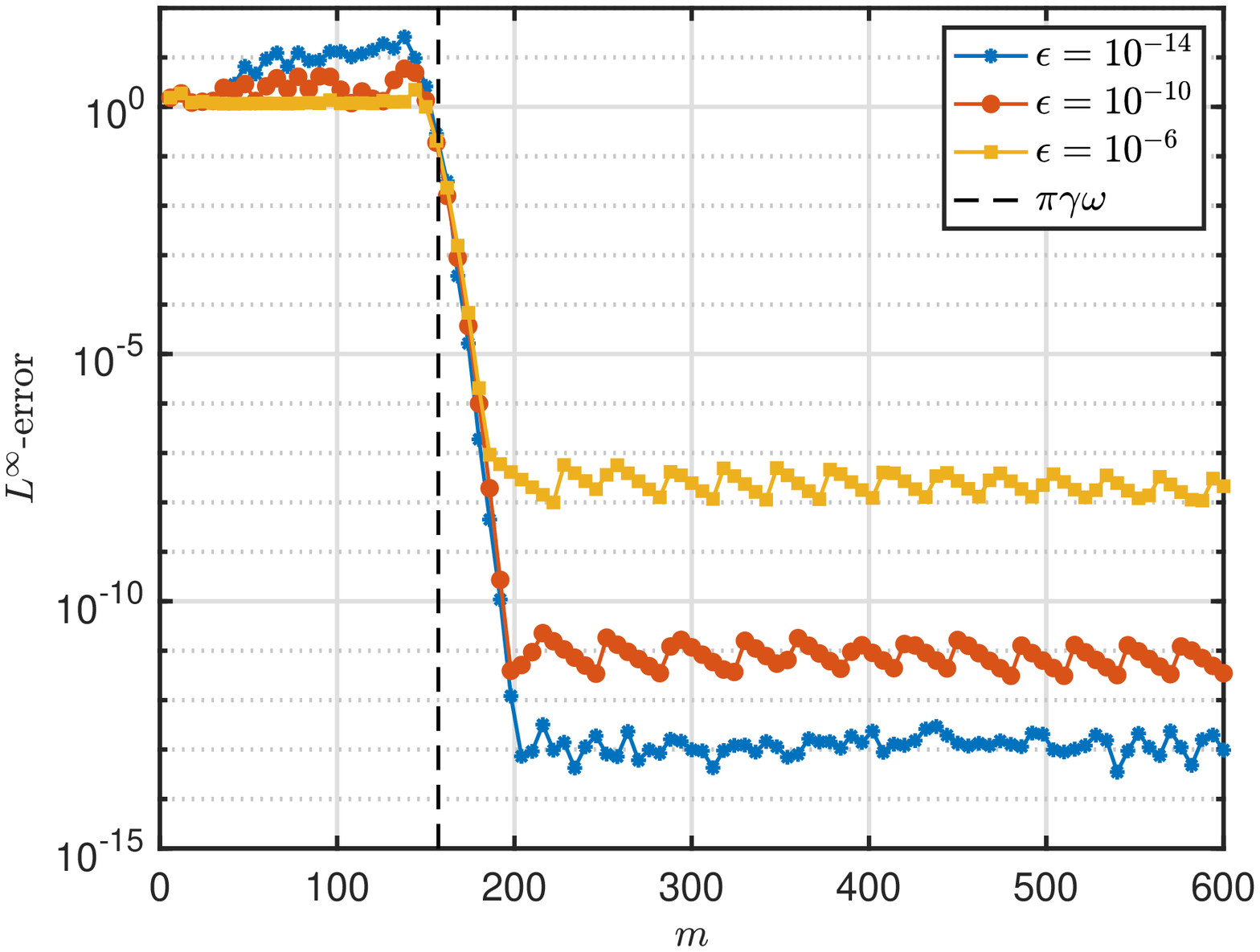}
&
\includegraphics[width = 0.3\textwidth]{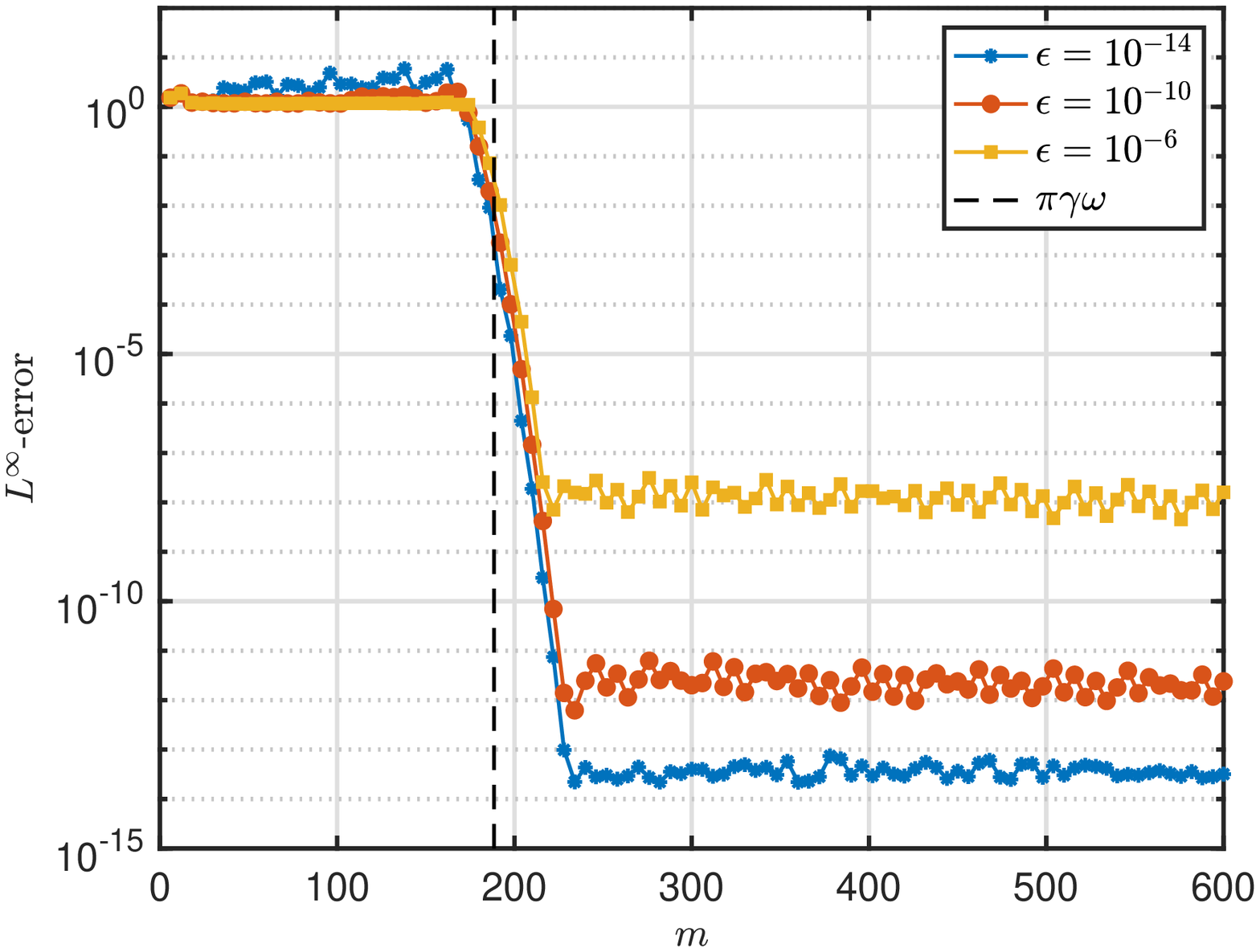}
&
\includegraphics[width = 0.3\textwidth]{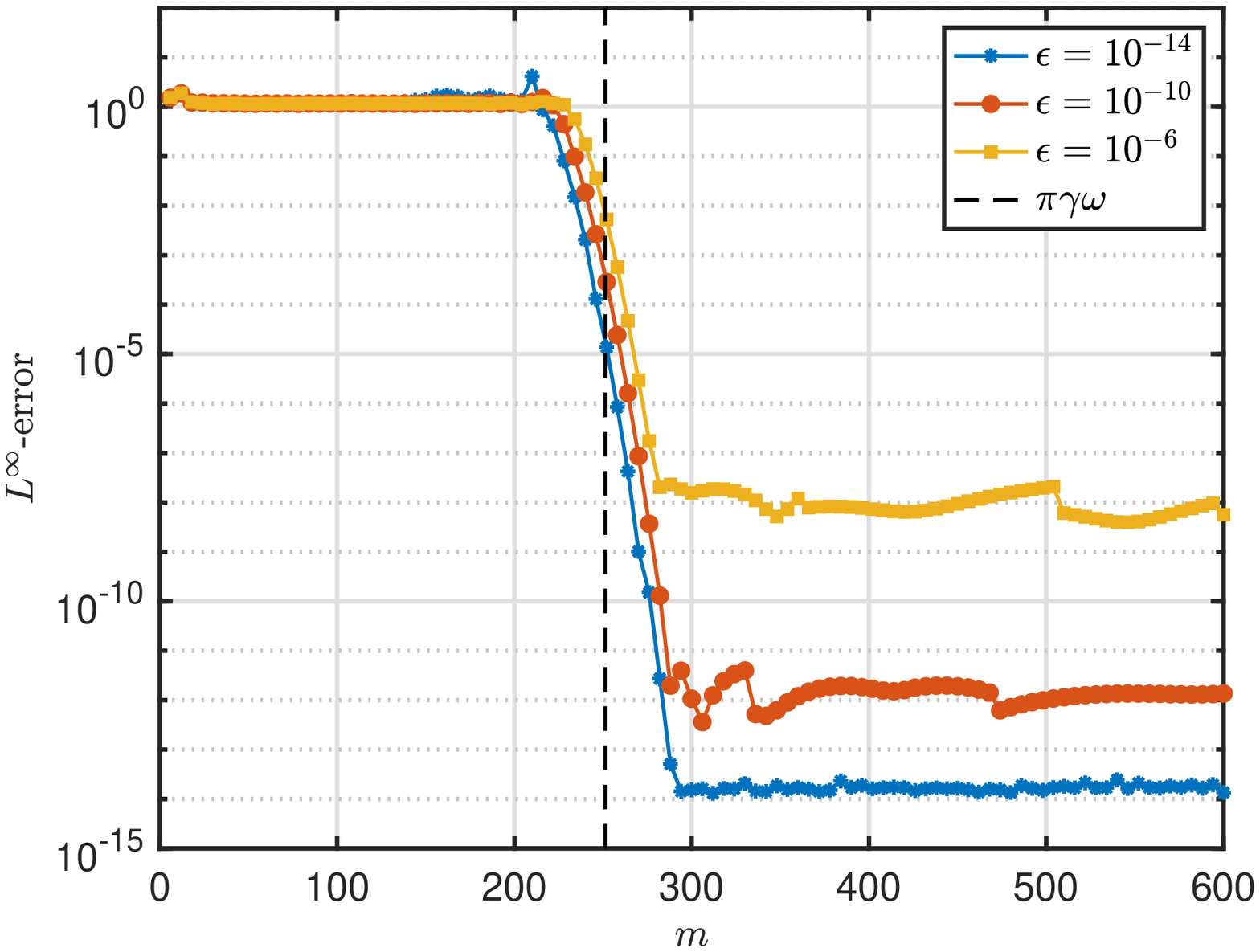}
\\
\includegraphics[width = 0.3\textwidth]{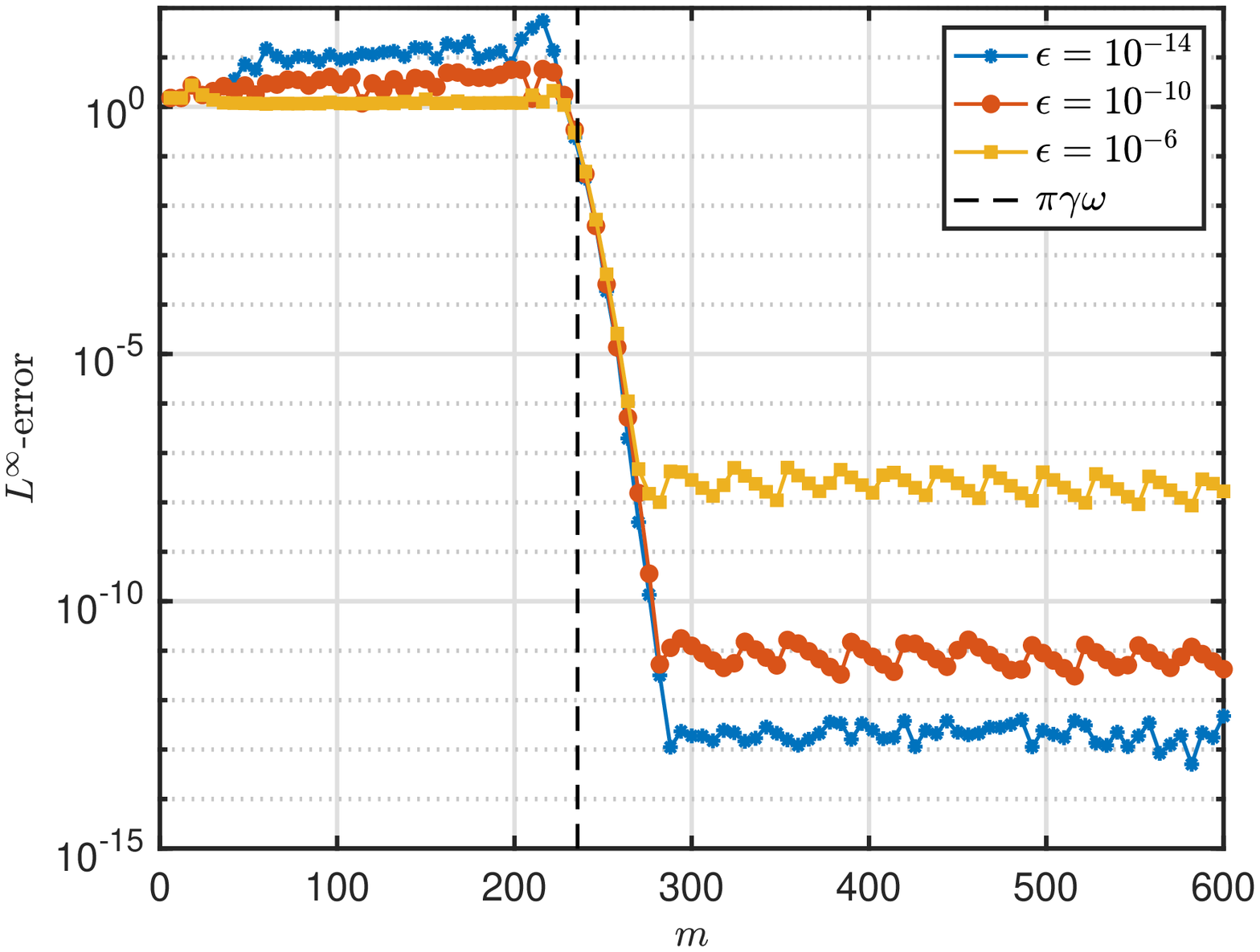}
&
\includegraphics[width = 0.3\textwidth]{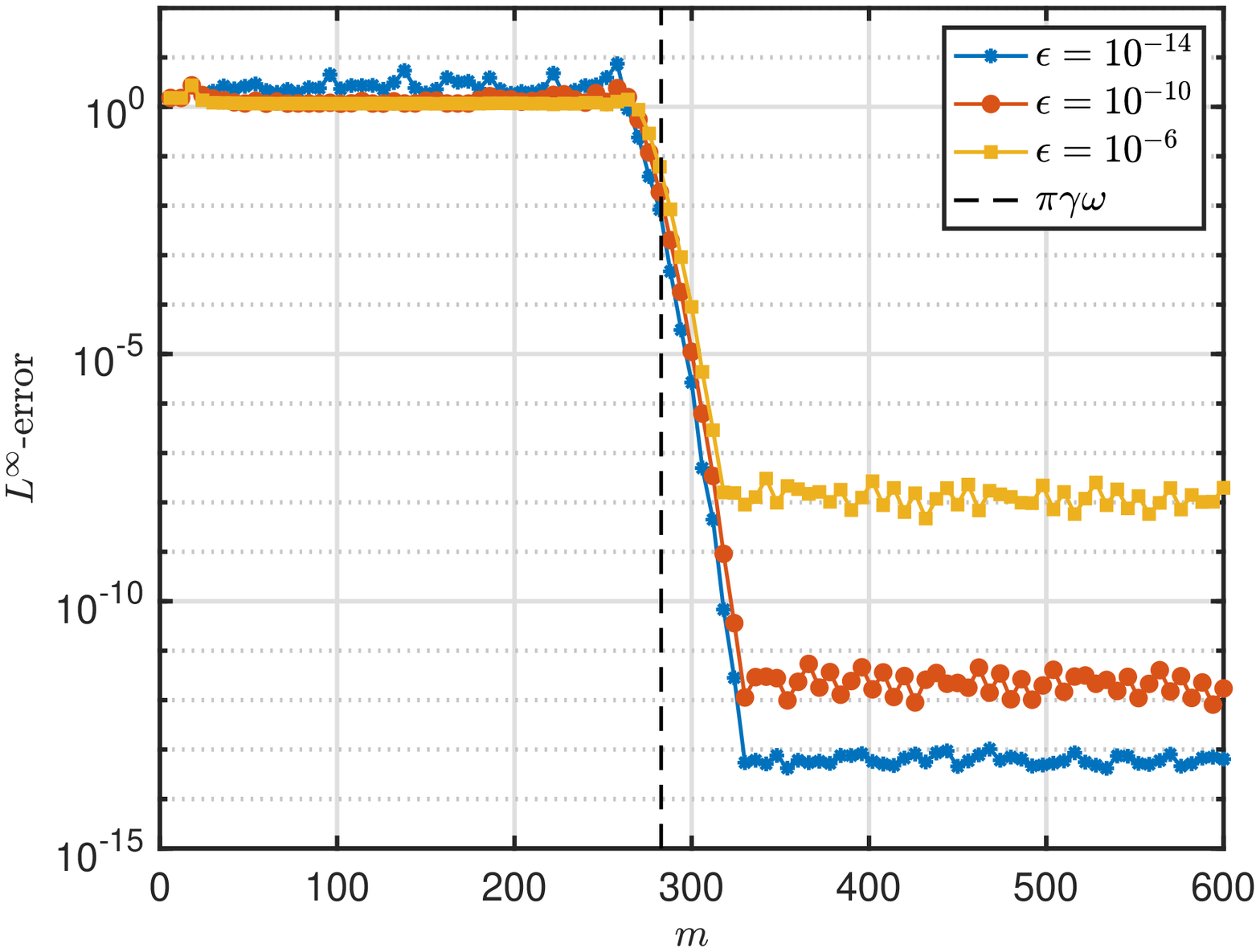}
&
\includegraphics[width = 0.3\textwidth]{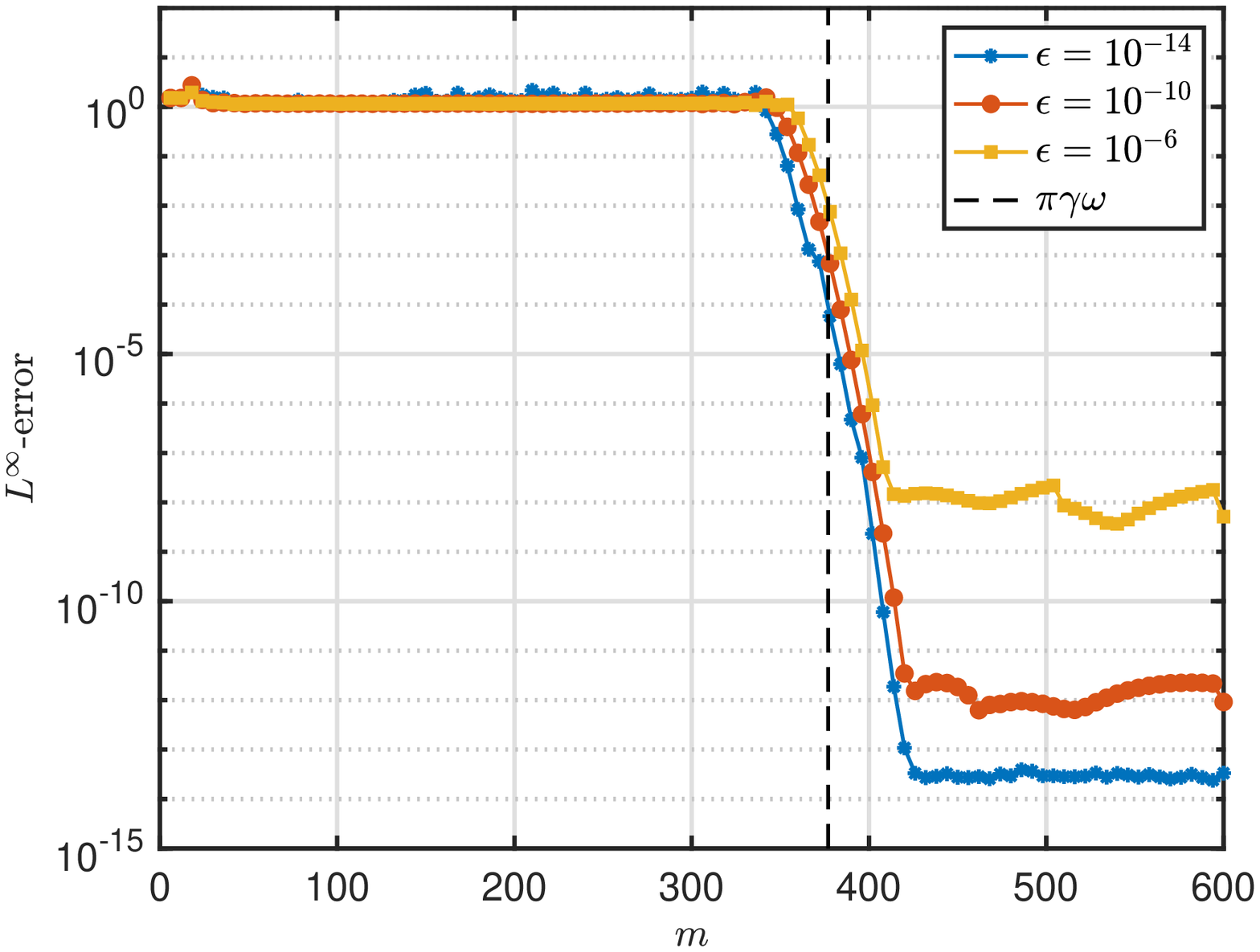}
\\
\includegraphics[width = 0.3\textwidth]{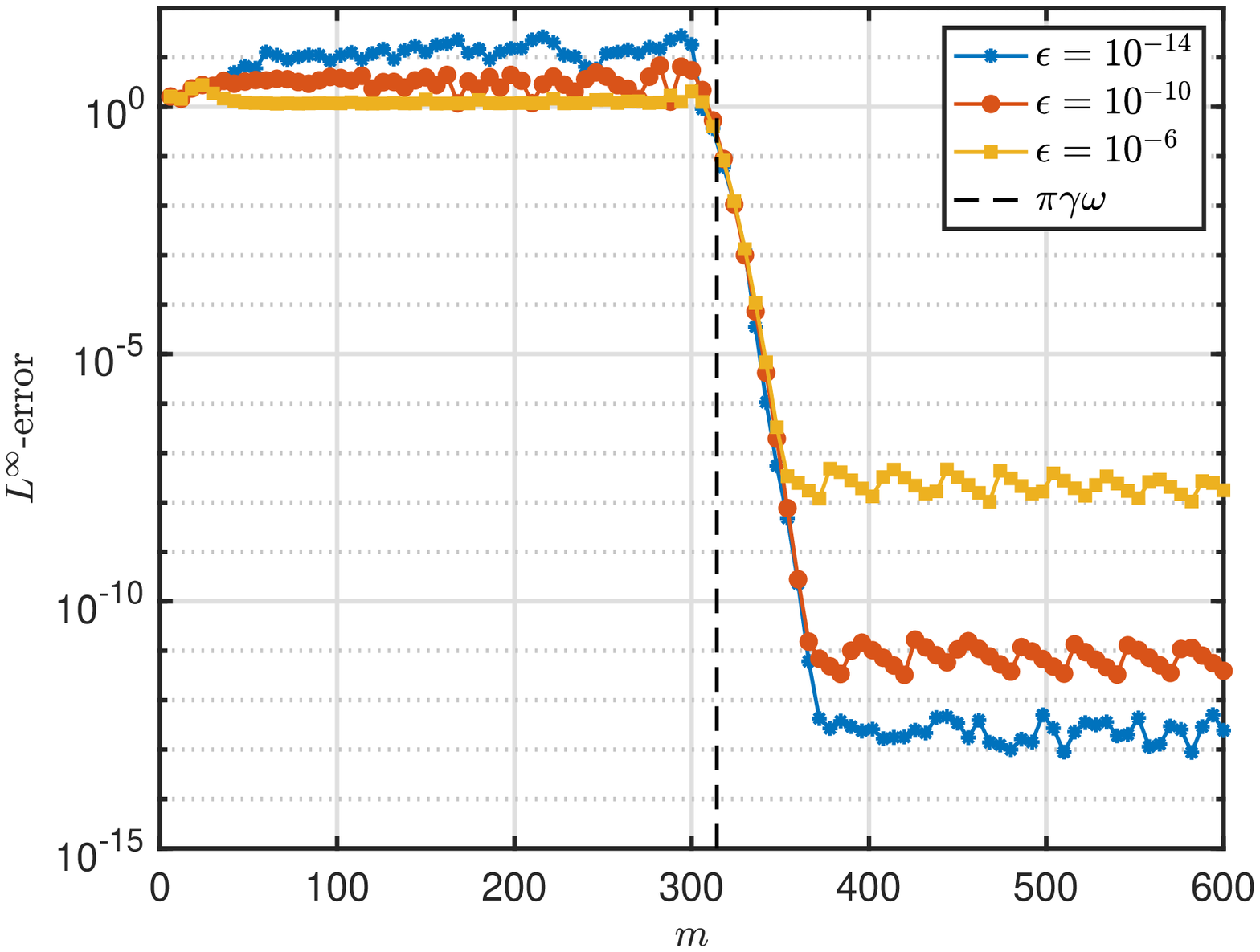}
&
\includegraphics[width = 0.3\textwidth]{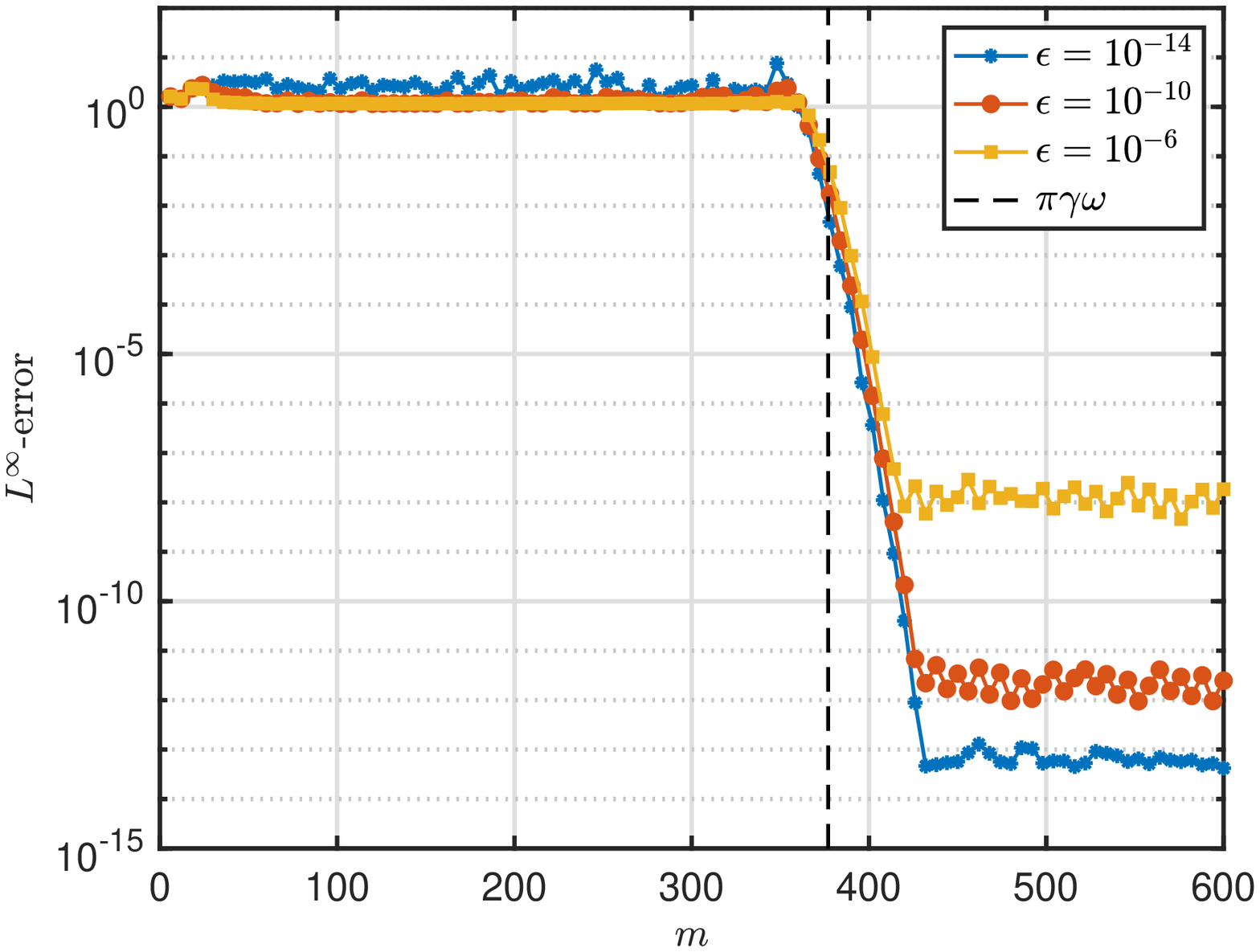}
&
\includegraphics[width = 0.3\textwidth]{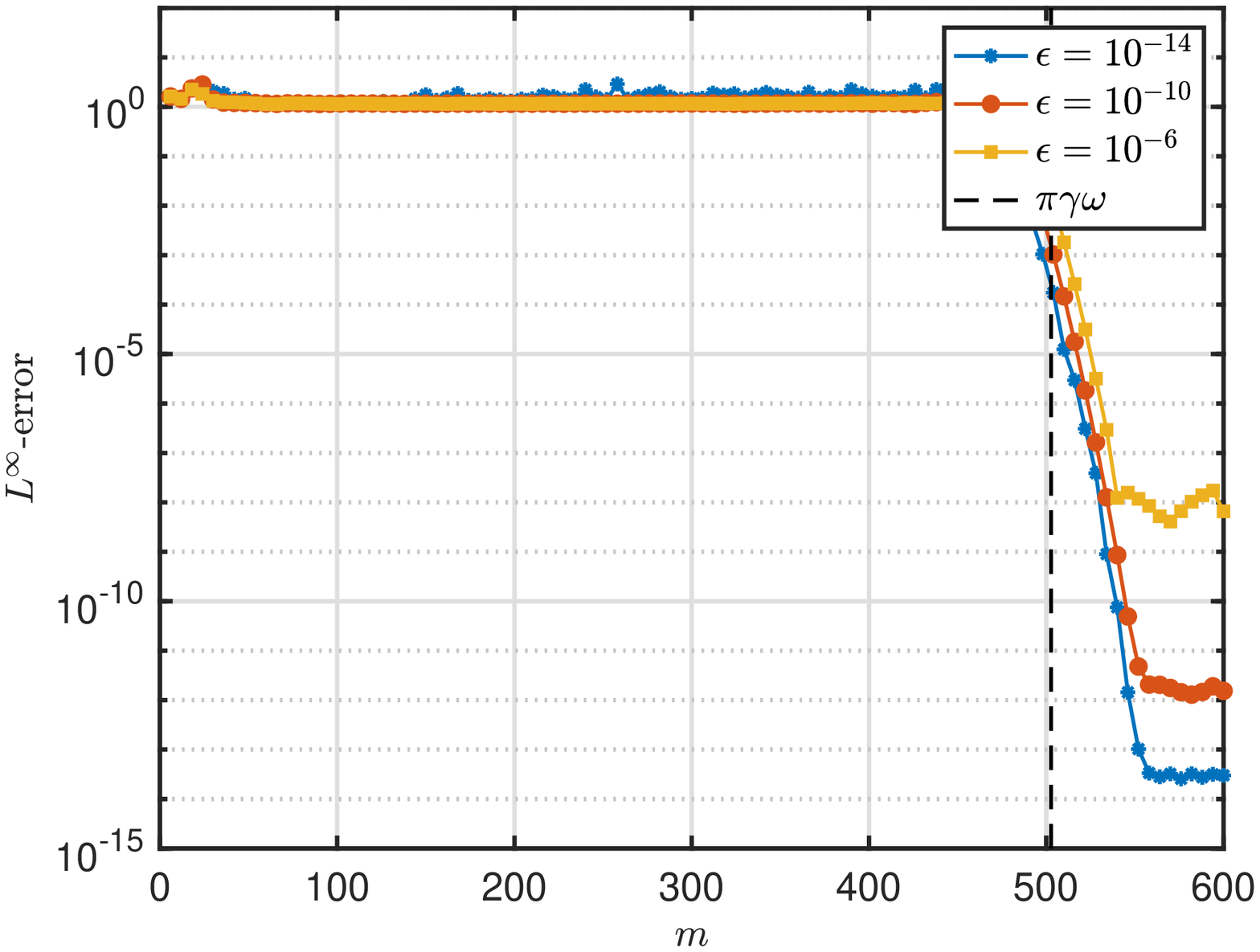}
\\
$\gamma = 1.25$ & $\gamma = 1.5$ & $\gamma = 2$
\end{tabular}
\end{center}
\end{small}
\caption{Approximation errors versus $n$ for approximating the functions $f(x) = \exp(\I \omega \pi x)$via $\cP^{\epsilon,\gamma}_{m,n}$, where $m/n = 4$, using various different values of $\gamma$ and $\epsilon$.  The values of $\omega$ are $\omega = 40$ (top), $\omega = 60$ (middle) and $\omega = 80$ (top).}
\label{f:fig3}
\end{figure} 

In Fig.\ \ref{f:fig3} we consider approximating the oscillatory function $f(x) = \E^{\I \omega \pi x}$ for various different values of $\omega$. Oscillatory functions are an interesting case study for approximation procedures from equispaced nodes. They are entire functions, yet they grow extremely rapidly along the imaginary axis for large $|\omega|$, meaning that the term $\nm{f}_{E_{\theta},\infty}$ is extremely large unless $\theta \approx 1$. As we see in Fig.\ \ref{f:fig3}, the approximation error is order one until a minimum value of $n = n_0(\omega)$ is met. After this point, the function begins to be resolved and the error decreases rapidly. Determining the behaviour of $n_0(\omega)$ allows us to examine the \textit{resolution power} of the approximation scheme, i.e.\ the number of points needed before decay of the error sets in.

For the polynomial frame approximation, we determine this point by recalling the error bound
\bes{
\inf_{p \in \bbP_n} \left \{ \nmu{f - p}_{[-1,1],\infty} + (n+1)\epsilon \nm{p}_{[-\gamma,\gamma],\infty} \right \}.
}
Since $f$ is entire and $|f(x)| = 1$ for real $x$, we can write
\bes{
\inf_{p \in \bbP_n} \left \{ \nmu{f - p}_{[-1,1],\infty} + (n+1)\epsilon \nm{p}_{[-\gamma,\gamma],\infty} \right \} \leq 2 \inf_{p \in \bbP_n}  \nmu{f - p}_{[-\gamma,\gamma],\infty} + (n+1)\epsilon.
}
Hence, the resolution power for polynomial approximation is related to the resolution power of best polynomial approximation on the extended interval $[-\gamma,\gamma]$. It is well known (see, e.g., \cite{hale2008new,boyd1989chebyshev,gottlieb1977numerical,adcock2014resolution}) that on the interval $[-1,1]$ an oscillatory function with frequency $\omega$ can be approximated by a polynomial once the degree $n$ exceeds the value $\pi \omega$. Since we consider the extended interval, this implies that $n_0(\omega) = \pi \gamma \omega$. 

This value is also shown in Fig.\ \ref{f:fig3}. It closely predicts the point at which the error begins to decrease. Note that this suggests choosing a small value of $\gamma$ so as to obtain a higher resolution power. Yet, as seen, this worsens the condition number of the approximation. Or, to put it another way, a smaller $\gamma$ necessitates a more severe oversampling ratio $\eta$ so as to maintain the same condition number, thus worsening the resolution power with respect to the number of equispaced samples $m$.

\begin{figure}[t]
\begin{small}
\begin{center}
\begin{tabular}{ccc}
\includegraphics[width = 0.3\textwidth]{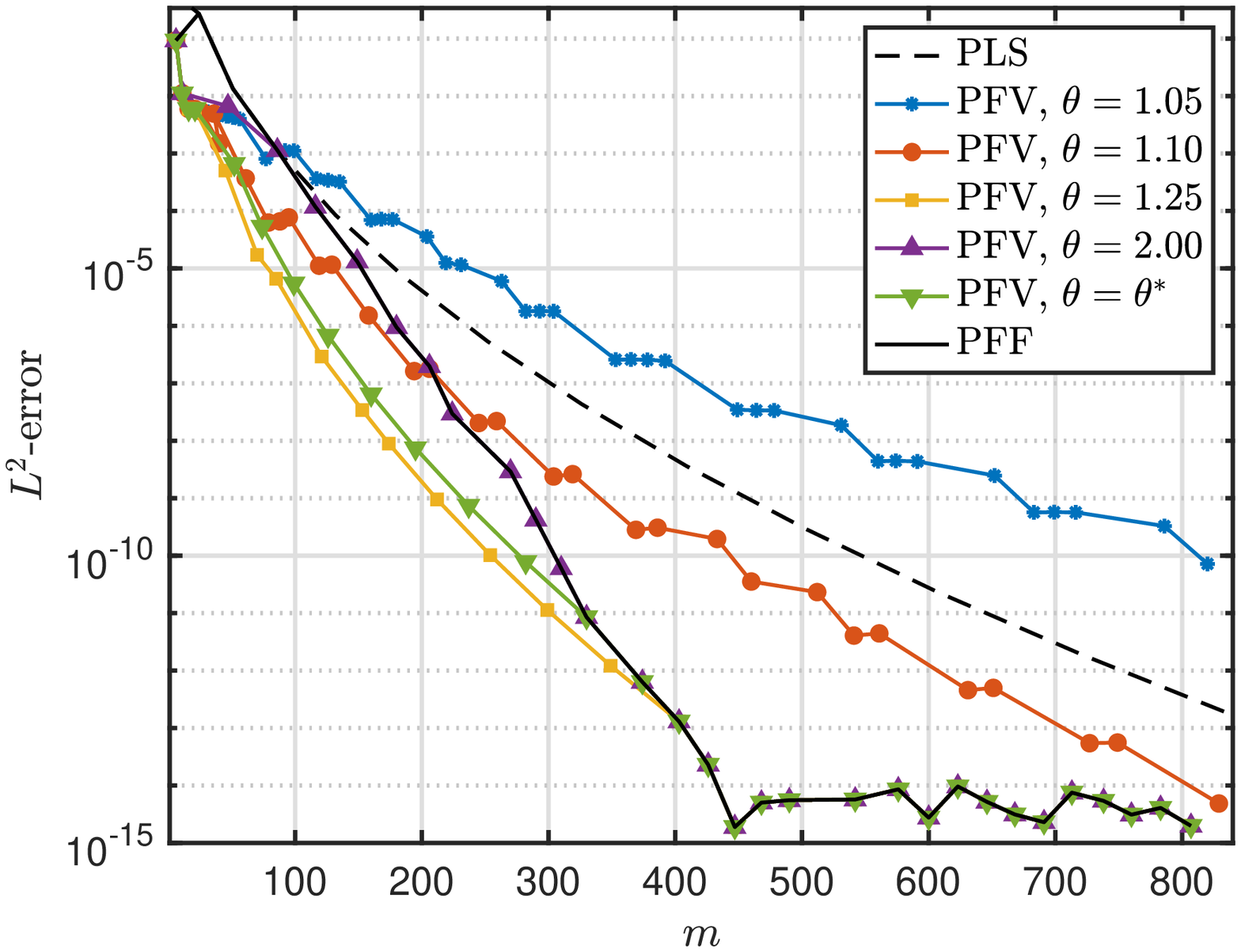}
&
\includegraphics[width = 0.3\textwidth]{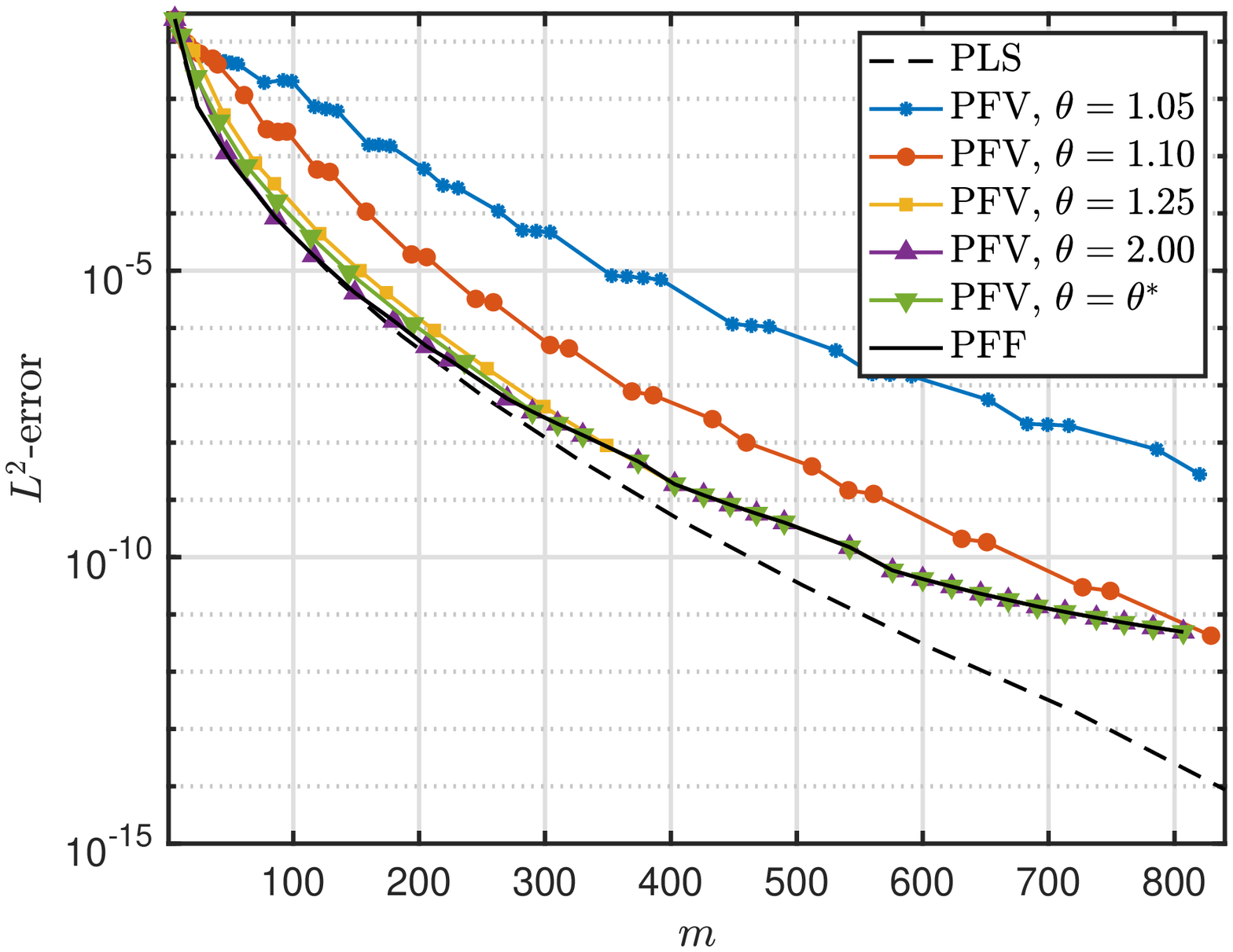}
&
\includegraphics[width = 0.3\textwidth]{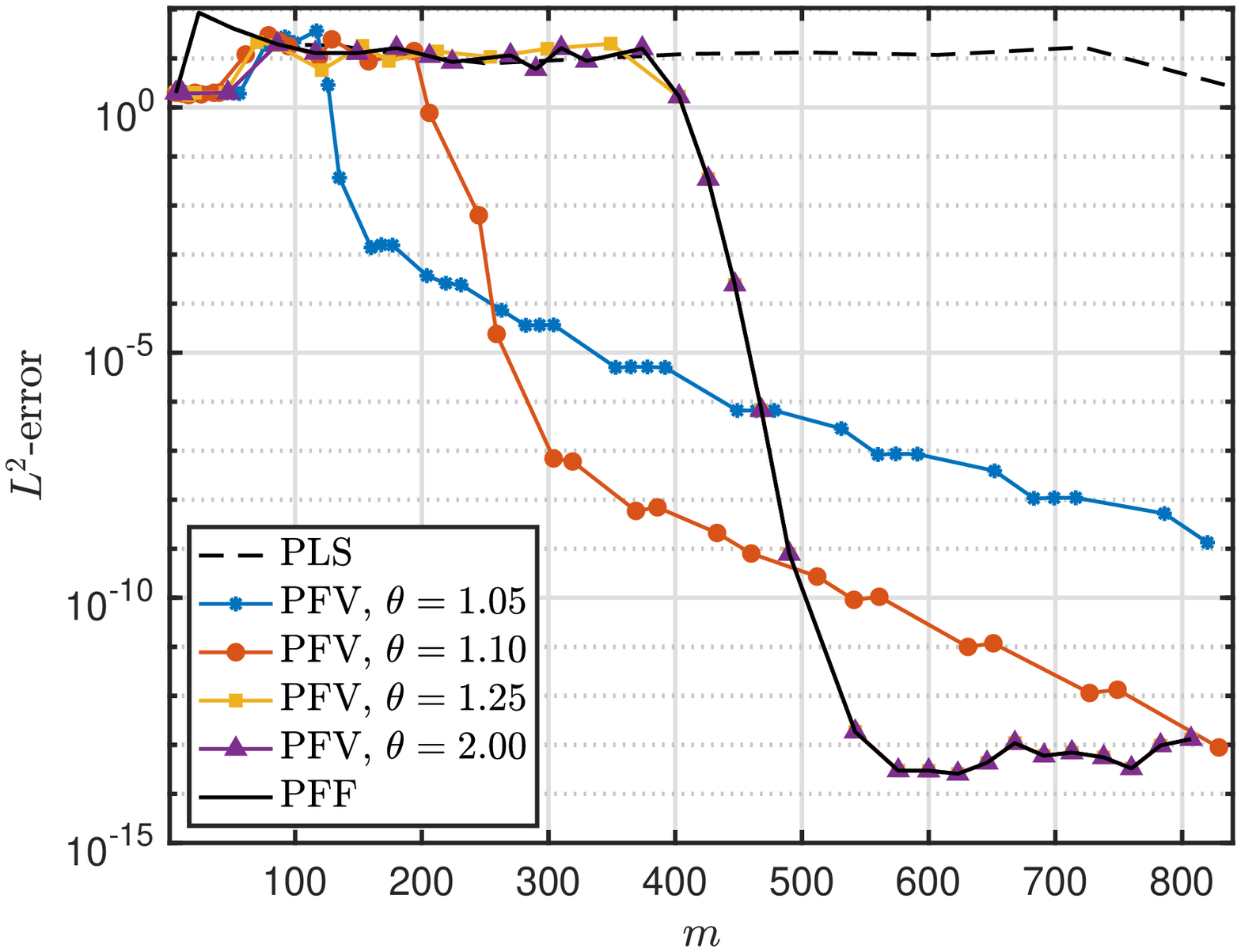}
\\
\includegraphics[width = 0.3\textwidth]{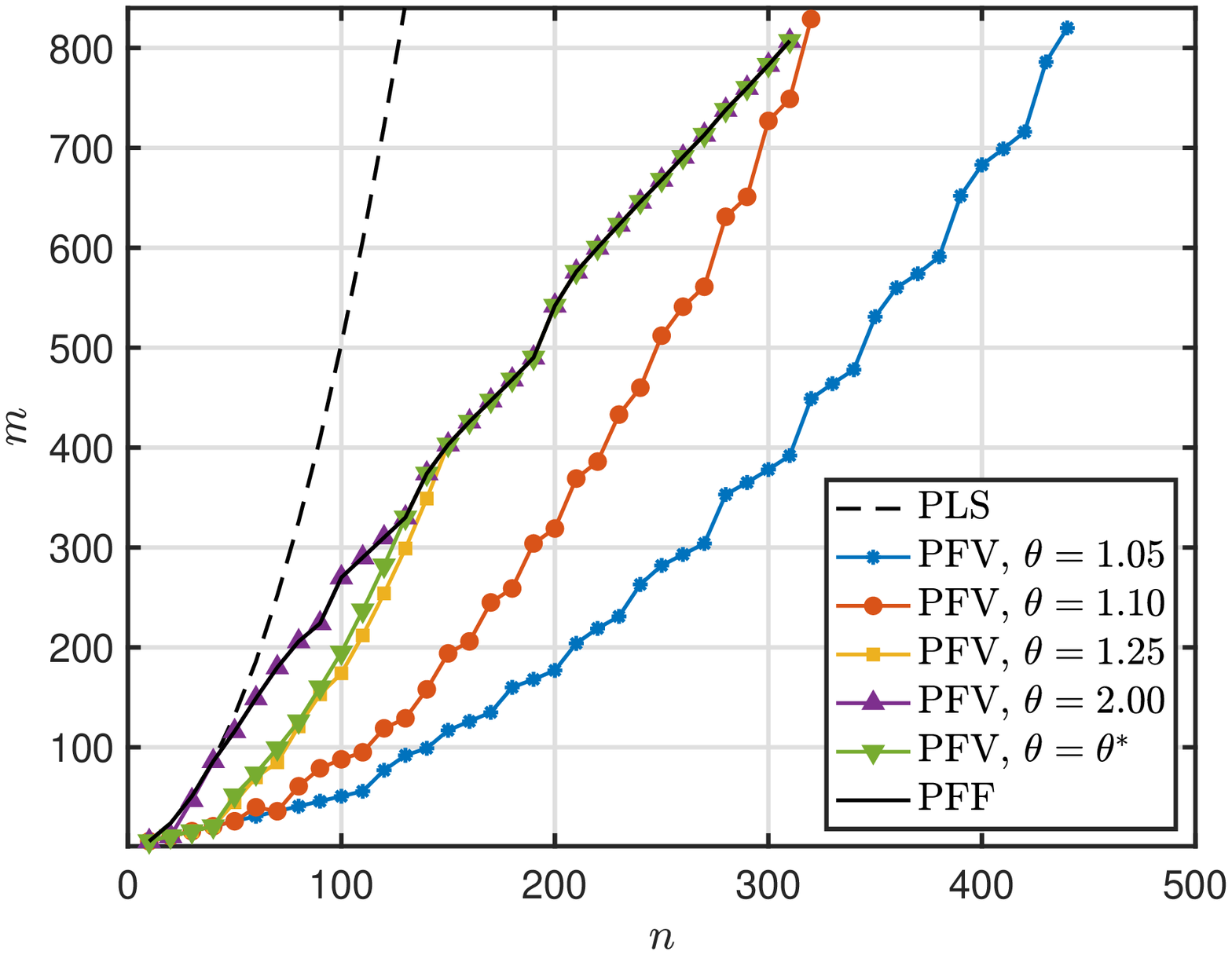}
&
\includegraphics[width = 0.3\textwidth]{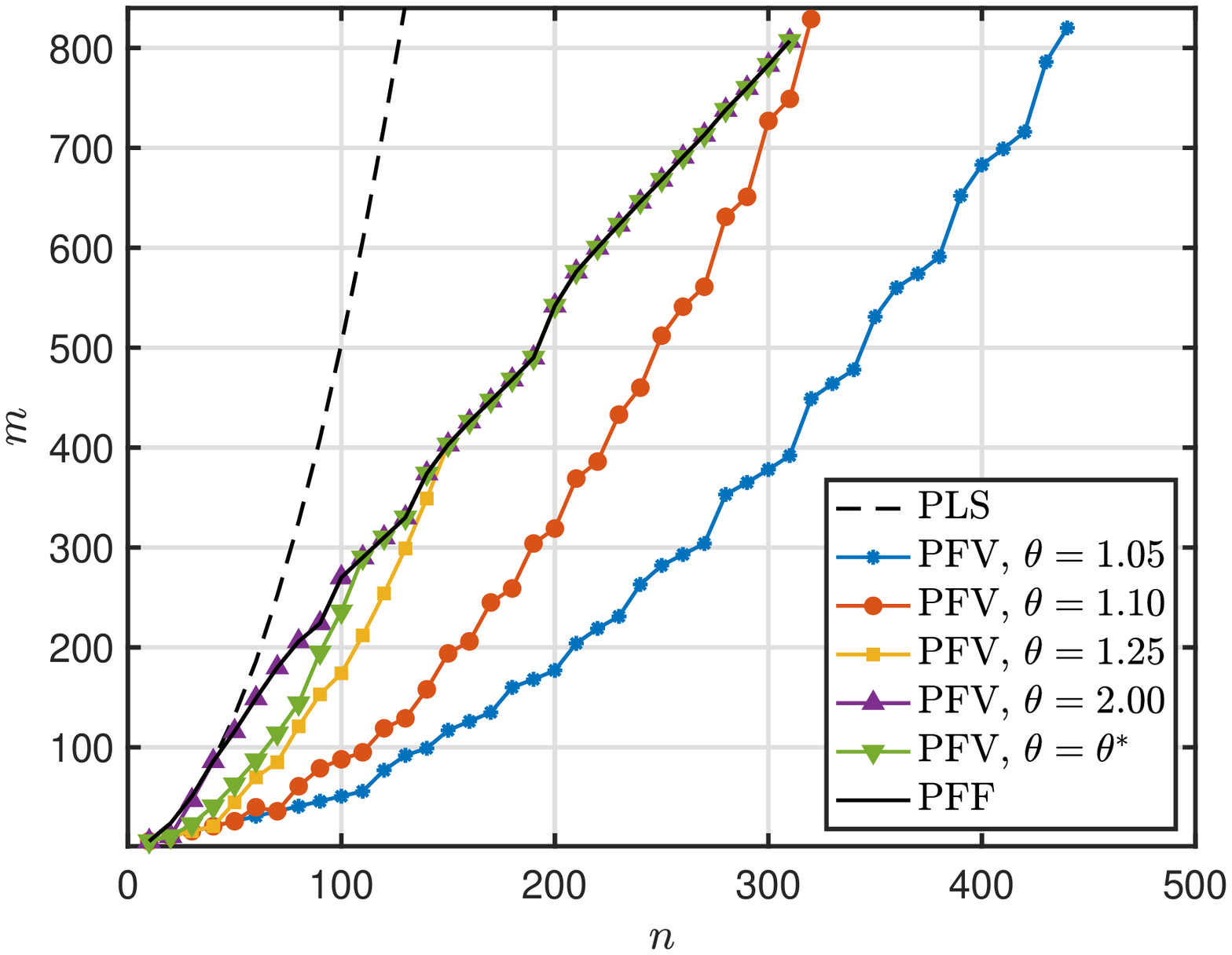}
&
\includegraphics[width = 0.3\textwidth]{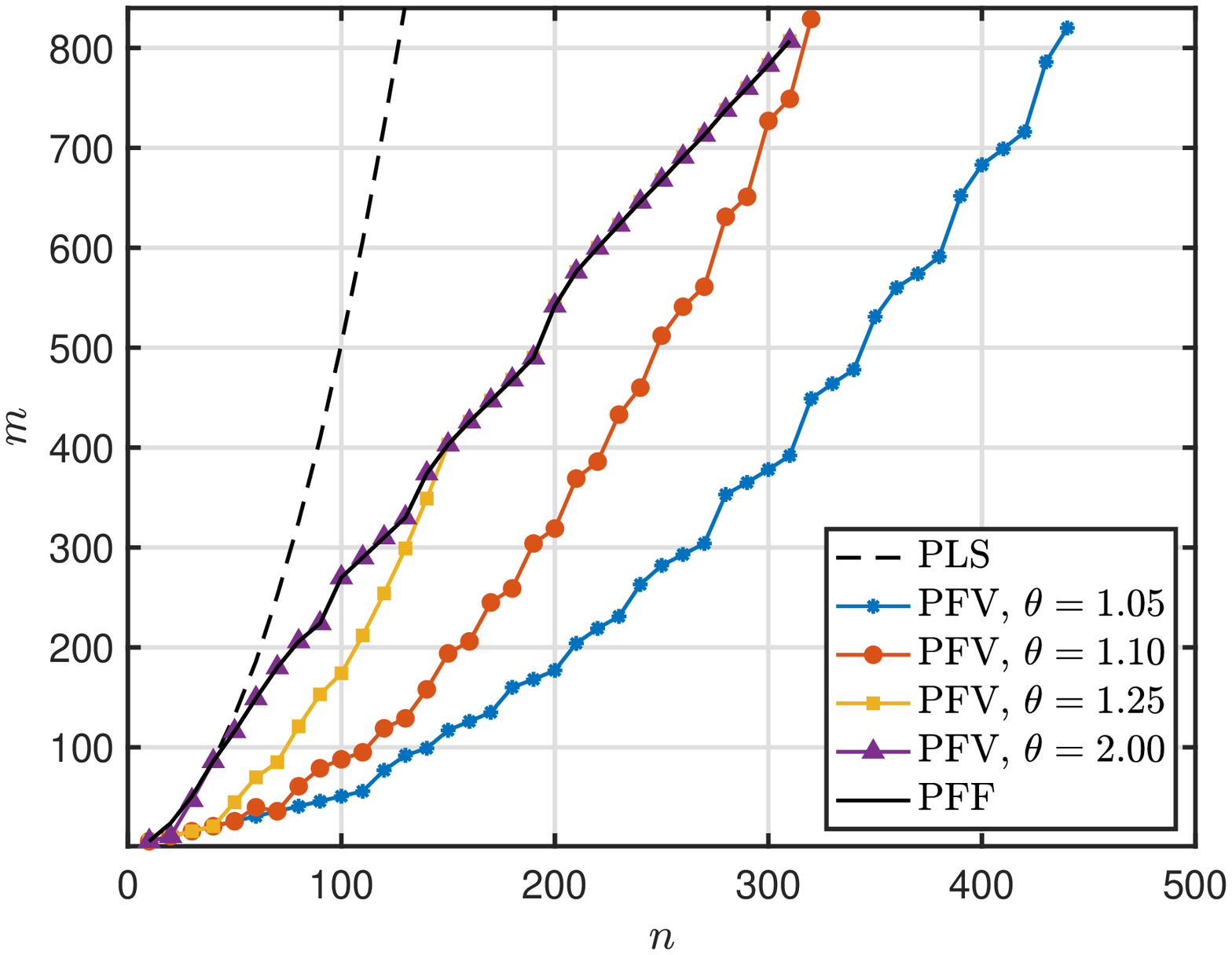}
\\
$f_1(x) = \frac{1}{1+16 x^2}$, $\theta^* = \frac{\sqrt{17}+1}{4}$ & $f_2(x) = \frac{1}{30-29 x}$,  $\theta^* = \frac{30 + \sqrt{59}}{29}$ & $f_3(x) = \E^{40 \pi \I  x}$, $\theta^* = \infty$
\end{tabular}
\end{center}
\end{small}
\caption{{Comparison of various schemes for approximating the functions $f_1$, $f_2$ and $f_3$. The top row shows the approximation error (measured in the discrete $L^2$-norm) versus $m$. The bottom row shows the relationship between $m$ and $n$ where, given $n$, $m$ is smallest integer such that the (discrete $L^2$-norm) condition number of each scheme is at most $100$. The schemes considered are: polynomial least-squares approximation (``PLS"), which uses orthonormal Legendre polynomials on $[-1,1]$ (see Remark \ref{r:why-not-unit-interval}); polynomial frame approximation with fixed $\epsilon = 10^{-14} $ (``PFF"); and polynomial frame approximation with varying $\epsilon = \max \{ \theta^{-n}, 10^{-14}  \}$ (``PFV"). Both PFF and PFV use the value $\gamma = 1.25$. The quantity $\theta^*$ is the `optimal' value of $\theta$, in the sense that it is the parameter of the largest Bernstein ellipse within which the given function is analytic. } }
\label{f:fig4}
\end{figure}

{ 
Finally, we conclude in Fig.\ \ref{f:fig4} with an experiment that compares fixed $\epsilon$ (the main setting in this paper) with varying $\epsilon$, the latter chosen as in Remark \ref{rem:varying-eps} to decay like $\theta^{-n}$ for a given $\theta > 1$. In order to make a fair comparison, we define a maximum allowable condition number $\kappa^* = 100$. Then, given $n$,  we find the largest value of $m$ for each scheme such that its condition number is at most $\kappa^*$. In order to do this, we use follow the approach of \cite[\S 8]{adcock2020approximating} and work in the discrete $L^2$-norm (as opposed to the discrete uniform norm considered in previous experiments) over a grid of 50,000 equispaced points in $[-1,1]$. Doing so means that the condition number can be computed as the norm of a certain matrix.

For the first function, $f_1$, we observe that varying $\epsilon$ can lead to a benefit whenever the parameter $\theta$ is chosen suitably: specifically, whenever it is chosen close to the value $\theta = \theta^*$, where $\theta^*$ is the largest Bernstein ellipse within which the function is analytic. If $\theta$ is chosen too large, then the scheme behaves similarly to the standard polynomial frame approximation with fixed $\epsilon$. On the other hand, if $\theta$ is chosen too small, then the scheme performs significantly worse.  In fact, it can even perform worse that standard polynomial least-squares approximation using orthonormal Legendre polynomials on $[-1,1]$ (see Remark \ref{r:why-not-unit-interval}). 

For the second function, $f_2$, varying $\epsilon$ conveys no benefit over fixed $\epsilon$, and, as before, it can lead to worse performance if $\theta$ is chosen too small. It is notable that polynomial least-squares approximation outperforms any of the polynomial frame approximations for this function. This is not surprising. The function grows rapidly near the endpoint $x = + 1$. In polynomial frame approximation the approximating polynomial is constrained to be of a finite size on $[-\gamma,\gamma]$. Hence it cannot capture this rapid growth as effectively as in polynomial least-squares approximation, where there is no such constraint.

Finally, the third function, $f_3$, is oscillatory and therefore entire. Observe that when $\theta$ is small, the polynomial frame approximation with varying $\epsilon$ initially resolves the function using fewer samples than the scheme with fixed $\epsilon$. Yet, as $m$ increases the error decays more slowly, and is eventually larger than the error for the latter method.

Finally, in the second row of Fig.\ \ref{f:fig4} we plot the scaling between $m$ and $n$. As expected, polynomial least-squares approximation exhibits a quadratic scaling, while polynomial frame approximation with fixed $\epsilon$ exhibits a linear scaling. Polynomial frame approximation with varying $\epsilon$ exhibits a quadratic scaling for small $\theta$. When $\theta$ is larger the scaling is at first quadratic and then linear.  This arises because $\epsilon$ is constrained to be no smaller than $10^{-14}$ in this experiments, this being done in order to avoid numerical effects in the thresholded SVD. Note that in this regime, the two polynomial frame approximations coincide.

The main conclusion from this experiment is that varying $\epsilon$ with $n$ can lead to some benefit (for small $m$), as long as the parameter $\theta$ is chosen carefully. Unfortunately, such a choice is function dependent (compare, for instance, $f_1$ versus $f_3$), and may require knowledge of the domain of analyticity of the (unknown) function.
}

\section{Conclusions and outlook}\label{s:conclusions}

In this work, we have shown a positive counterpart to the impossibility theorem of \cite{platte2011impossibility}. Namely, we have shown a possibility theorem (Theorem \ref{t:possibility-thm}) for polynomial frame approximation, which asserts stability and exponential convergence down to a finite, but user-controlled limiting accuracy. This holds for all functions that are analytic in a sufficiently large region -- a condition that must be the case for any scheme achieving this type of error decay, as shown by our extended impossibility theorem (Theorem \ref{t:impossibility-extended}). On the other hand, for insufficiently analytic functions, we have shown exponential decay down to some fractional power of $\epsilon$ (Theorem \ref{t:possibility-slower-exp}), and superalgebraic decay (Theorem \ref{t:possibility-algebraic}) beyond this point.

There are several avenues for further investigation. First, recall that our main error bounds involve (albeit small) algebraic factors in $m$ and $n$. It would be interesting to see if such factors could be removed by modifying the approximation scheme. Alternatively, one might instead work in the $L^2$-norm, which is arguably more natural for least-squares approximations.

Second, this work was inspired by so-called Fourier extensions \cite{adcock2014parameter,adcock2014numerical,matthysen2017function,lyon2012fast,huybrechs2010fourier,bruno2007accurate,boyd2002fourier,pasquetti1996spectral}, wherein a smooth, nonperiodic function on $[-1,1]$ is approximated by a Fourier series on $[-\gamma,\gamma]$. In practice, linear oversampling appears sufficient for accuracy and stability of $\epsilon$-regularized Fourier extensions, with exponential error decay down to $\epsilon$ \cite{adcock2014parameter,adcock2014numerical}. Proving a similar possibility theorem for this scheme is an open problem. Note that Fourier extension is equivalent to a polynomial approximation problem on an arc of the complex unit circle \cite{geronimo2020fourier,webb2020pointwise}. Fourier extension schemes have several advantages over the polynomial extension scheme studied herein. For example, they generally possess higher resolution power \cite{adcock2014resolution} (recall Fig.\ \ref{f:fig3} and the discussion in \S \ref{s:numerical}).

Third, we mention that equispaced points are not special. Similar impossibility theorems have been shown for scattered data, or more generally, any sample points that do not cluster quadratically near the endpoints $x = \pm 1$ \cite{adcock2019optimal}. Beyond pointwise samples, it is notable that an impossibility theorem also holds for reconstructing analytic functions from their Fourier samples \cite{adcock2012stable,adcock2014generalized,adcock2014stability}. The question of whether or not possibility theorems hold for these more general types of samples is also an open problem.

Fourth, notice that we have not concerned ourselves with fast computation of the polynomial frame approximation in this work. We anticipate, however, that a fast implementation may be possible, as it is with Fourier extensions \cite{matthysen2015fast,lyon2012fast,matthysen2017function}. One potential idea in this direction is the AZ algorithm \cite{coppe2020AZ}.

Fifth and finally, we note that the one-dimensional problem is, in some senses, a toy problem. Polynomial frame approximations we first formalized in \cite{adcock2020approximating} to accurately and stably approximate functions that are defined over general, compact domains in two or more dimensions. Here the domain is embedded in a hypercube, and a tensor-product orthogonal polynomial basis on the hypercube is used to construct the approximation. Such approximations are often used in practice in surrogate model construction problems in uncertainty quantification \cite{adcock2022sparse,adcock2020approximating}. Showing that linear oversampling is sufficient in two or more dimensions and a corresponding possibility theorem for analytic function approximation in arbitrary dimensions would be an interesting and practically-relevant extension.

\small
\bibliographystyle{plain}
\bibliography{polyleastsquaresbib}

\begin{thebibliography}{10}

\bibitem{adcock2022sparse}
B.~Adcock, S.~Brugiapaglia, and C.~G. {Webster}.
\newblock {\em Sparse Polynomial Approximation of High-Dimensional Functions}.
\newblock Society for Industrial and Applied Mathematics, Philadelphia, PA,
  2022.

\bibitem{adcock2012stable}
B.~Adcock and A.~C. Hansen.
\newblock Stable reconstructions in {H}ilbert spaces and the resolution of the
  {G}ibbs phenomenon.
\newblock {\em Appl. Comput. Harmon. Anal.}, 32(3):357--388, 2012.

\bibitem{adcock2014generalized}
B.~Adcock and A.~C. Hansen.
\newblock Generalized sampling and the stable and accurate reconstruction of
  piecewise analytic functions from their {F}ourier coefficients.
\newblock {\em Math. Comp.}, 84:237--270, 2014.

\bibitem{adcock2014stability}
B.~Adcock, A.~C. Hansen, and A.~Shadrin.
\newblock A stability barrier for reconstructions from {F}ourier samples.
\newblock {\em SIAM J. Numer. Anal.}, 52(1):125--139, 2014.

\bibitem{adcock2014resolution}
B.~Adcock and D.~Huybrechs.
\newblock On the resolution power of {F}ourier extensions for oscillatory
  functions.
\newblock {\em J. Comput. Appl. Math.}, 260:312--336, 2014.

\bibitem{adcock2019frames}
B.~Adcock and D.~Huybrechs.
\newblock Frames and numerical approximation.
\newblock {\em SIAM Rev.}, 61(3):443--473, 2019.

\bibitem{adcock2020approximating}
B.~Adcock and D.~Huybrechs.
\newblock Approximating smooth, multivariate functions on irregular domains.
\newblock {\em Forum Math. Sigma}, 8:e26, 2020.

\bibitem{adcock2020frames}
B.~Adcock and D.~Huybrechs.
\newblock Frames and numerical approximation {II}: generalized sampling.
\newblock {\em J. Fourier Anal. Appl.}, 26(6):87, 2020.

\bibitem{adcock2014numerical}
B.~Adcock, D.~Huybrechs, and J.~Mart{\'\i}n-Vaquero.
\newblock On the numerical stability of {F}ourier extensions.
\newblock {\em Found. Comput. Math.}, 14(4):635--687, 2014.

\bibitem{adcock2016mapped}
B.~Adcock and R.~Platte.
\newblock A mapped polynomial method for high-accuracy approximations on
  arbitrary grids.
\newblock {\em SIAM J. Numer. Anal.}, 54(4):2256--2281, 2016.

\bibitem{adcock2019optimal}
B.~Adcock, R.~Platte, and A.~Shadrin.
\newblock Optimal sampling rates for approximating analytic functions from
  pointwise samples.
\newblock {\em IMA J. Numer. Anal.}, 39(3):1360--1390, 2019.

\bibitem{adcock2014parameter}
B.~Adcock and J.~Ruan.
\newblock Parameter selection and numerical approximation properties of
  {F}ourier extensions from fixed data.
\newblock {\em J. Comput. Phys.}, 273:453--471, 2014.

\bibitem{boyd1989chebyshev}
J.~P. Boyd.
\newblock {\em Chebyshev and Fourier Spectral Methods}.
\newblock Springer--Verlag, 1989.

\bibitem{boyd2002fourier}
J.~P. Boyd.
\newblock A comparison of numerical algorithms for {F}ourier {E}xtension of the
  first, second, and third kinds.
\newblock {\em J. Comput. Phys.}, 178:118--160, 2002.

\bibitem{boyd2009exponentially}
J.~P. Boyd and J.~R. Ong.
\newblock Exponentially-convergent strategies for defeating the {R}unge
  phenomenon for the approximation of non-periodic functions. {I}.
  {S}ingle-interval schemes.
\newblock {\em Commun. Comput. Phys.}, 5(2--4):484--497, 2009.

\bibitem{bruno2007accurate}
O.~P. Bruno, Y.~Han, and M.~M. Pohlman.
\newblock Accurate, high-order representation of complex three-dimensional
  surfaces via {F}ourier continuation analysis.
\newblock {\em J. Comput. Phys.}, 227(2):1094--1125, 2007.

\bibitem{bullen2015dictionary}
P.~S. Bullen.
\newblock {\em Dictionary of Inequalities}.
\newblock Monographs and Research Notes in Mathematics. CRC Press, Boca Raton,
  FL, 2nd edition, 2015.

\bibitem{cheney1982introduction}
E.~W. Cheney.
\newblock {\em Introduction to Approximation Theory}, volume 317.
\newblock AMS Chelsea Publishing, Providence, RI, 2nd edition, 1982.

\bibitem{christensen2016introduction}
O.~Christensen.
\newblock {\em An Introduction to Frames and Riesz Bases}.
\newblock Applied and Numerical Harmonic Analysis. Birkh{\"a}user, Basel, 2nd
  edition, 2016.

\bibitem{coppe2020AZ}
V.~Copp{\'e}, D.~Huybrechs, R.~Matthysen, and M.~Webb.
\newblock {The AZ algorithm for least squares systems with a known incomplete
  generalized inverse}.
\newblock {\em SIAM J. Mat. Anal. Appl.}, 41(3):1237--1259, 2020.

\bibitem{coppersmith1992growth}
D.~Coppersmith and T.~J. Rivlin.
\newblock The growth of polynomials bounded at equally spaced points.
\newblock {\em SIAM J. Math. Anal.}, 23:970--983, 1992.

\bibitem{don1997accuracy}
W.-S. Don and A.~Solomonoff.
\newblock Accuracy enhancement for higher derivatives using {C}hebyshev
  collocation and a mapping technique.
\newblock {\em {SIAM} J. Sci. Comput.}, 18:1040--1055, 1997.

\bibitem{geronimo2020fourier}
J.~S. Geronimo and K.~Liechty.
\newblock The {F}ourier extension method and discrete orthogonal polynomials on
  an arc of the circle.
\newblock {\em Adv. Math.}, 365:107064, 2020.

\bibitem{gottlieb1977numerical}
D.~Gottlieb and S.~A. Orszag.
\newblock {\em Numerical Analysis of Spectral Methods: Theory and
  Applications}.
\newblock Society for Industrial and Applied Mathematics, Philadelphia, PA, 1st
  edition, 1977.

\bibitem{hale2008new}
N.~Hale and L.~N. Trefethen.
\newblock New quadrature formulas from conformal maps.
\newblock {\em SIAM J. Numer. Anal.}, 46(2):930--948, 2008.

\bibitem{huybrechs2010fourier}
D.~Huybrechs.
\newblock On the {F}ourier extension of non-periodic functions.
\newblock {\em SIAM J. Numer. Anal.}, 47(6):4326--4355, 2010.

\bibitem{konyagin2021stable}
S.~Konyagin and A.~Shadrin.
\newblock On stable reconstruction of analytic functions from selected
  {F}ourier coefficients.
\newblock {\em In preparation}, 2021.

\bibitem{kosloff1993modified}
D.~Kosloff and H.~Tal-Ezer.
\newblock A modified {C}hebyshev pseudospectral method with an
  {$\mathcal{O}(N^{-1})$} time step restriction.
\newblock {\em J. Comput. Phys.}, 104:457--469, 1993.

\bibitem{lyon2012fast}
M.~Lyon.
\newblock A fast algorithm for {F}ourier continuation.
\newblock {\em {SIAM} J. Sci. Comput.}, 33(6):3241--3260, 2012.

\bibitem{matthysen2015fast}
R.~Matthysen and D.~Huybrechs.
\newblock Fast algorithms for the computation of {F}ourier extensions of
  arbitrary length.
\newblock {\em SIAM J. Sci. Comput.}, 38(2):A899--A922, 2016.

\bibitem{matthysen2017function}
R.~Matthysen and D.~Huybrechs.
\newblock Function approximation on arbitrary domains using fourier frames.
\newblock {\em SIAM J. Numer. Anal.}, 56(3):1360--1385, 2018.

\bibitem{michelli1985lecture}
C.~A. Micchelli and T.~J. Rivlin.
\newblock Lectures on optimal recovery.
\newblock In P.~R. Turner, editor, {\em Numerical Analysis Lancaster 1984},
  volume 1129 of {\em Lecture Notes in Mathematics}. Springer, Berlin
  Heidelberg, 1985.

\bibitem{pasquetti1996spectral}
R.~Pasquetti and M.~Elghaoui.
\newblock A spectral embedding method applied to the advection--diffusion
  equation.
\newblock {\em J. Comput. Phys.}, 125:464--476, 1996.

\bibitem{platte2011impossibility}
R.~Platte, L.~N. Trefethen, and A.~Kuijlaars.
\newblock Impossibility of fast stable approximation of analytic functions from
  equispaced samples.
\newblock {\em SIAM Rev.}, 53(2):308--318, 2011.

\bibitem{schaeffer1938some}
A.~C. Schaeffer and R.~J. Duffin.
\newblock {On some inequalities of S.\ Bernstein and W.\ Markoff for
  derivatives of polynomials,}.
\newblock {\em Bull. Amer. Math. Soc.}, 44(4):289--297, 1938.

\bibitem{shadrin2004twelve}
A.~Shadrin.
\newblock Twelve proofs of the {M}arkov inequality.
\newblock In {\em Approximation Theory: A volume dedicated to Borislav
  Bojanov}, pages 233--298. Professor Marin Drinov Academic Publishing House,
  Sofia, 2004.

\bibitem{temlyakov2018multivariate}
V.~Temlyakov.
\newblock {\em Multivariate Approximation}, volume~32 of {\em Cambridge
  Monographs on Applied and Computational Mathematics}.
\newblock Cambridge University Press, Cambridge, UK, 2018.

\bibitem{temlyakov1993approximate}
V.~N. Temlyakov.
\newblock On approximate recovery of functions with bounded mixed derivative.
\newblock {\em J. Complexity}, 9(1):41--59, 1993.

\bibitem{trefethen2013approximation}
L.~N. Trefethen.
\newblock {\em Approximation Theory and Approximation Practice}.
\newblock Society for Industrial and Applied Mathematics, Philadelphia, PA,
  2013.

\bibitem{webb2020pointwise}
M.~Webb, V.~Copp{\'e}, and D.~Huybrechs.
\newblock Pointwise and uniform convergence of {F}ourier extensions.
\newblock {\em Constr. Approx.}, 52:139--175, 2020.

\end{thebibliography}

\end{document}